\documentclass[notitlepage,11pt,reqno]{amsart}
\usepackage{amssymb,amsmath,amscd,xspace,amsthm,tocvsec2,mathrsfs}
\usepackage{verbatim,color}
\usepackage[breaklinks=true]{hyperref}

\usepackage[letterpaper]{geometry}
\newcommand{\tangle}[1]{\ensuremath{\langle #1 \rangle}}
\newcommand{\bse}{\boldsymbol\varepsilon}
\newcommand{\bxi}{\boldsymbol\xi}
\newcommand{\bla}{\boldsymbol\lambda}
\renewcommand{\geq}{\geqslant}
\renewcommand{\leq}{\leqslant}
\DeclareMathOperator{\Res}{Res}
\DeclareMathOperator{\sgn}{sgn}
\DeclareMathOperator{\inte}{int}

\newtheorem{theorem}{Theorem}[section]
\newtheorem{prop}[theorem]{Proposition}
\newtheorem{lemma}[theorem]{Lemma}
\newtheorem{corollary}[theorem]{Corollary}

\newtheorem{utheorem}{\textrm{\textbf{Theorem}}}

\theoremstyle{definition}
\newtheorem{remark}[theorem]{Remark}
\newtheorem{defn}[theorem]{Definition}
\newtheorem{example}[theorem]{Example}
\newtheorem{qn}{Question}

\numberwithin{equation}{section}
\allowdisplaybreaks

\begin{document}


\title[Strong BMV property, commutation of operators, and
applications]{Strong Bessis--Moussa--Villani property, commutation of
operators, and applications}

\author[Sujit Sakharam Damase and Apoorva Khare]{Sujit Sakharam Damase
and Apoorva Khare\\ (with an Appendix by Vamsi P.~Pingali)}

\address[S.S.~Damase]{Department of Mathematics, Indian Institute of
Science, Bangalore, India;
Department of Mathematics, Massachusetts Institute of Technology,
Cambridge, USA;
PhD student}
\email{\tt sujits@iisc.ac.in; sujits@mit.edu}

\address[A.~Khare]{Department of Mathematics, Indian
Institute of Science, Bangalore, India}
\email{\tt khare@iisc.ac.in}

\address[V.P.~Pingali]{Department of Mathematics, Indian
Institute of Science, Bangalore, India}
\email{\tt vamsipingali@iisc.ac.in}

\date{\today}

\keywords{Bessis--Moussa--Villani conjecture, Hurwitz products,
commutativity of operators, Loewner order, strong BMV property;
positive operator, $C^*$-algebra, spectral measure, Weyl functional
calculus, Gibbs thermal weight, multivariate completely monotone
function; elementary symmetric polynomial, Schur polynomial, Jack
polynomial, Macdonald polynomial; maximum construction, Griffiths
negativity}

\subjclass[2020]{47B93, 
81R15, 
15A27, 
05E05, 
33D52, 
32L15, 
47Axx, 
15A15 (primary); 
44A10, 
46L05, 
47N50 
(secondary)}

\begin{abstract}
Fix an integer $l \geq 1$ and Hermitian matrices $B \geq 0$ and $A$ of
any (equal) size. We prove that the Hurwitz products $S_{m,l}(A,B) =
\left. \frac{1}{l!} \frac{d^l}{dt^l} (A+tB)^m \, \right|_{t=0}$ are
positive semidefinite for all $m \geq l$ only if $AB=BA$.
In fact, we prove a twofold strengthening:
(a)~$A,B$ now are self-adjoint elements of any $C^*$-algebra, with $B
\geq 0$ if $l \geq 2$; and
(b)~we accommodate negative powers $m$, provided $A>0$.
Moreover, (c)~if one only assumes $S_{m,l}(A,B) \geq 0$ for infinitely
many $m$ (possibly a zero density set) then one recovers ``partial
commutativity'': $A^2B = BA^2$. We show that even this need not follow if
$S_{m,l}(A,B) \geq 0$ for only finitely many $m$.

This result has numerous applications. Firstly, we resolve
(a)~a question from Lieb--Seiringer [\textit{J.\ Stat.\ Phys.}\ 2012], as
well as
(b)~a symmetrized variant of a 2002 question of Lieb.
These questions relate to the 1975 BMV conjecture (shown by Stahl
[\textit{Acta Math.}\ 2013] and Hein\"avaara [\textit{Invent.\ Math.}\
2024]), and ask when stronger conditions than in that conjecture hold.

Second, working with (non-)commuting operators $A$ and $B$ naturally
leads one to the functional calculus of Weyl and Wigner. Our results
recover -- under significantly weaker hypotheses --
the main theorem of Anderson [\textit{J.\ Funct.\ Anal.}\ 1969]: if the
Weyl distribution for self-adjoint operators $A_1, \dots, A_d$ is a
positive measure, then the $A_j$ commute pairwise.

Third, we return to the complete monotonicity in the BMV conjecture, and
characterize commuting of operators $B = H_1 \geq 0$ and $A = H_0 =
H_0^*$ via operator-monotonicity of the Gibbs thermal weight
$e^{-\beta(H_0 + t H_1)}$ for all (inverse) temperatures $\beta$.

Our fourth application is a strengthening of the BMV conjecture/theorem
itself. Given a positive definite matrix $A$, say $d_A$ is the
eigenspace-dimension for its largest eigenvalue. Let $g(\sigma(M))$
denote any monomial-positive symmetric polynomial in the eigenvalues of
$M = M^*$, whose constituent monomial functions $m_\lambda$ are indexed
by partitions with $\leq d_A$ parts. Alternately, for such $\lambda$ let
$g$ be the corresponding zonal, Schur, Jack, Hall--Littlewood, or even
Macdonald polynomial. We show that for every $l \geq 1$ and $B>0$ there
is a threshold $m_l = m_l(A,B,g) > 0$ such that
\[
g(\sigma(S_{m,l}(A,B))) > 0 \qquad \forall m \geq m_l,
\]
and similarly, $g(\sigma((-1)^l S_{-m,l}(A,B))) > 0 \ \forall m \geq m_l$
if we replace $d_A$ by $d_{A^{-1}}$.
(a)~The simplest case is $g = e_1$, aka the BMV conjecture.
We refine this for each individual $A$, from $g = e_1$ to $\{ e_1, \dots,
e_{d_A} \}$.
(b)~We show that the bounds of $d_{A^{\pm 1}}$ are best possible,
answering in detail the 2012~question arising from
Lieb--Seiringer's work.
(c)~We also show that unlike the BMV conjecture for $e_1$, the ``all
large $m \gg 0$'' conclusion cannot be improved to $m \geq l$ for $e_2,
\dots, e_{d_A}$. In particular, this shows that the original BMV
conjecture (aka Stahl's theorem) is remarkably special.

A fifth application is to operator algebras. We recover a host of results
on Loewner monotonicity of certain maps implying commutativity in
$C^*$-algebras -- by Ogasawara (1955), Pedersen (1979), Uchiyama (1993),
Wu (2001), and Moln\'ar (2017). We also resolve a question of
Ricker--V\"ath (2002).

Finally, an Appendix by Vamsi P.~Pingali contains a (negative)
application of our result for $l=1$, in geometry. Various natural
constructions of Hermitian metrics on holomorphic vector bundles are
explored, motivated by their analogues in the case of line bundles, in
the hope of developing a ``pluripotential theory'' of vector bundles. The
Appendix shows that none of these can work, by providing counterexamples
to Griffiths positivity/negativity for them all.
\end{abstract}

\maketitle

\tableofcontents
\settocdepth{section}

\section{Introduction and main results}\label{sec:Intro}

This work has two disparate motivations:
(a)~A question from Lieb--Seiringer's 2012 work, as well as another by
Lieb in $\sim$2002. These questions ask when two strengthenings of the
condition in the Bessis--Moussa--Villani (BMV) conjecture hold. This 1975
conjecture, originating in quantum statistical
mechanics and resolved in 2013, has led to much subsequent activity.
This work leads to a novel, BMV-type characterization of
when two self-adjoint Hilbert space operators (or even $2 \times 2$ or $n
\times n$ matrices) commute -- a question going back to the foundations
of quantum mechanics.

We were also motivated by
(b)~the question of Griffiths negativity for certain natural
constructions of Hermitian metrics on holomorphic complex bundles, aiming
towards a pluripotential theory. This question also has ties to mirror
symmetry. Indeed, an (as of now hypothetical) vector bundle version of the
deformed Hermitian--Yang--Mills equation under a degeneration would
produce singular metrics satisfying some positivity condition on
curvature. Griffiths positivity\footnote{The dual metric of a Griffiths
positive metric is Griffiths negative, so both notions will be used
interchangeably via dualization.}
is the simplest such condition.

In this section, we provide a detailed account of our motivations,
interweaving our results and applications into it -- with the
geometry motivation and results sketched here and detailed in
Appendix~\ref{sec:geom} by Vamsi P.~Pingali.

\subsection{Motivation from quantum physics; Lieb--Seiringer's and Lieb's
questions in matrix analysis}\label{Sstart}

In 1975, Bessis--Moussa--Villani (BMV) formulated a conjecture while
studying quantum partition functions
$Z(t) := {\rm Tr}(e^{-\beta (H_0 + t H_1)})$. Here,
\begin{itemize}
\item $H_0$ is an unperturbed Hamiltonian (e.g.\ kinetic energy),
\item $H_1$ is an interaction Hamiltonian (e.g.\ a repulsion potential),
\item $t$ is the coupling constant / perturbation strength,
\item multiplying by $\beta = 1/(k_{Bol} T)$ makes the terms
dimensionless (where $k_{Bol}$ is the Boltzmann constant and $T$ the
temperature),
\item and one works with $A = \beta H_0$ and $B = \beta H_1$.
\end{itemize}
One takes the operator $B$ to be positive here (as will we in much of
this paper) -- i.e.\ nonnegative interaction energy -- in order to avoid
in certain physical settings a thermodynamically unstable system without
ground state (as $t \to \infty$).

While the quantum operators $H_0, H_1$ (or $A,B$) do not commute in
general, Bessis--Moussa--Villani \textbf{conjectured}~\cite{BMV} that the
$l$th order perturbations of $Z(t)$ alternate in sign. That is, $t
\mapsto Z(t) = {\rm Tr}(e^{-(A+tB)})$ is a completely monotone function
(see Section~\ref{S1agniva} for the definition) if $A,B$ are
positive.\footnote{BMV originally worked with $B$ positive and any
self-adjoint $A$ -- which is equivalent since perturbing $A$ by a scalar
$c\, {\rm Id}$ in the exponent simply rescales the partition function by
$e^{-c}$.}
If true, this conjecture would help prove convergence of Pad\'e
approximants used in quantum field models. It would also show that for
all integers $m \geq l \geq0$, the total statistical weight of all closed
quantum paths of length $m$ that scatter exactly $l$ times off $H_1$, is
always nonnegative regardless of any destructive quantum interference.

Bessis--Moussa--Villani (BMV) proved their conjecture in some special
cases, including the following two situations:
(a)~when $A,B$ are $2 \times 2$ matrices; and
(b)~when $A,B$ are commuting operators over a Hilbert space, which is the
case in classical statistical mechanics (and which also leads us to the
Hilbert space setting, below).

In~2002, the BMV conjecture was reformulated into the following
equivalent -- and seemingly simple-sounding -- assertion:

\begin{theorem}[Lieb--Seiringer, \cite{LS}]\label{TLieb}
The BMV conjecture~\cite{BMV} is equivalent to the following claim:

For all positive integers $m,n$ and all $n \times n$ positive
semidefinite (henceforth termed \textbf{psd})\footnote{Recall that a real
symmetric (or complex
Hermitian) $n \times n$ matrix $A$ is positive semidefinite, henceforth
abbreviated as \textit{psd}, if $h^* A h \geq 0$ for all real (or
complex) vectors $h \neq 0$. It is positive definite if this inequality
is always strict. It is folklore that $A$ is psd if and only if all its
eigenvalues are in $[0,\infty)$.}
Hermitian matrices $A,B$, the trace-polynomial $t \mapsto {\rm Tr}((A + t
B)^m)$ has all coefficients nonnegative.
\end{theorem}

Over the last five decades, the BMV conjecture has resulted in tremendous
activity, including various special cases proved early on. It was
eventually resolved in~{2013} in a celebrated work~\cite{Stahl} by Stahl,
and strengthened even further in~2024 using tracial joint spectral
measures by Hein\"avaara~\cite{Heina}.\medskip

We now come to the present work, which characterizes stronger conditions
than the BMV trace-positivity. In this we were inspired by
Bessis--Moussa--Villani themselves, who explored at length a conjecture
in~\cite[Section~II]{BMV} stronger than the ``named'' BMV-conjecture.

The first problem we consider was formulated by Lieb and Seiringer (2012)
as a \textbf{more precise question} than the BMV conjecture. Note that
$(A+tB)^m$ is a polynomial in $t$ whose coefficients are sums of
\textit{matrix words}, i.e.\ non-commutative strings of the form $A B A A
B \dots$ with $m$ letters, which are evaluated via matrix multiplication.
We now make this precise by writing the coefficient of each monomial
$t^l$.

\begin{defn}\label{DSml}
Given integers $m, l \geq 0$ and matrices $A,B$ of the same size,
define the corresponding \textit{Hurwitz product} $S_{m,l}(A,B)$ to be
the sum of all words in $A,B$ of length $m$ with exactly $l$-many $B$s,
all evaluated under multiplication. We set $S_{0,0}(A,B) := {\rm Id}$.
\end{defn}

Some of the ``edge cases'' are:
\[
S_{m,0}(A,B) = A^m, \qquad S_{m,m}(A,B)= B^m, \qquad
S_{m,l}(A,B) = 0, \qquad \forall l > m \geq 0;
\]
and our initial object of interest below is:
\begin{equation}
S_{m,1}(A,B) := A^{m-1} B + A^{m-2} B A + \cdots + A B A^{m-2} + B
A^{m-1}.
\end{equation}

Note that if $A,B$ are symmetric or Hermitian, then so are all matrices
$S_{m,l}(A,B)$. Moreover, the Bessis--Moussa--Villani polynomial
$(A+tB)^m$ is the generating function of these Hurwitz products, or
equivalently, they are its Taylor coefficients:
\begin{equation}\label{Esml}
(A+tB)^m = \sum_{l=0}^\infty t^l S_{m,l}(A,B), \qquad
S_{m,l}(A,B) = \left. \frac{1}{l!} \, \frac{d^l}{dt^l} (A+tB)^m \,
\right|_{t=0}, \qquad \forall m,l \geq 0.
\end{equation}

Now the BMV conjecture says that ${\rm Tr} (S_{m,l}(A,B)) \geq 0$
for all $0 \leq l \leq m$ if $A,B$ are psd (respectively, positive if
$A,B$ are positive definite). Lieb and Seiringer discussed a
\textit{stronger phenomenon} than the BMV conjecture: when are all of the
coefficient-matrices $S_{m,l}(A,B)$ themselves positive (semi)definite --
see the discussion in \cite[Section~3]{LS2} (and \cite[Remark~6]{LS}).
For the matrices $A,B$ for which this holds, taking traces would
immediately yield the BMV conjecture for $A,B$.

We provide a bit more detail. In \cite[Section~3]{LS2}, Lieb--Seiringer
note that by Ballantine's work~\cite{Ball69}, there exist positive
definite matrices $A,B \in \mathbb{C}^{n \times n}$ such that
$\det(S_{2,1}(A,B)) = \det(AB+BA) < 0$, even though individually $AB, BA$
have positive eigenvalues.\footnote{We nontrivially generalize this below
-- see Example~\ref{Exeventual}.}
They add later that:
\begin{itemize}
\item If $W_1, \dots, W_{\binom{m}{l}}$ are the words in $A,B$ that occur
in $S_{m,l}(A,B)$ then the BMV conjecture implies
$\sum_{i=1}^{\binom{m}{l}} e_k(W_i) \geq 0$ for all $1 \leq k \leq n$
(and conversely). Here $e_k(W_i)$ is the $k$th elementary symmetric
polynomial in the eigenvalues of $W_i$.
\item Lieb--Seiringer further mention: it is false in general that
$e_k(\sum_{i \geq 1} W_i) \geq 0$, except for $k=1$ when the trace is
additive (so we reduce to the previous point).
\end{itemize}

The first point is now known affirmatively, given that the BMV conjecture
is resolved. One goal of this paper is to address the second point; and
indeed, we completely classify the pairs of positive matrices $(A,B)$ for
which the latter happens.\medskip

A \textbf{second question} in this regard was posed even earlier by Lieb.
In/before 2002, he asked Hillar and Johnson if for every (length $m$)
matrix-word in positive definite matrices $A,B$, the eigenvalues are all
positive. Note that some of these words have all eigenvalues positive,
whence so is their trace -- for instance, if $A,B$ are positive definite
then $B^l A^{m-l}$ has the same eigenvalues as $B^{l/2} A^{m-l} B^{l/2}$,
which is positive definite. While exploring Lieb's question, Johnson and
Hillar found~\cite{HJ} that while it holds for $2 \times 2$ matrices (see
\cite[Fact~6]{HJ}), it is not always so for larger matrices -- for
instance, there exist $3 \times 3$ positive definite matrices $A,B$ such
that $BABAAB$ has negative trace. (However, the sum of traces of all
words is indeed positive, by the BMV conjecture / Stahl's theorem.)

One problem with the preceding question is that e.g.\ the word $BABAAB$
is not Hermitian -- and it may be possible for general $m,l$ to exist a
word whose eigenvalues are not all real (and then asking if their sum is
real and positive is a later question). A close variant of Lieb's
question, therefore, would be to pair each word with its ``reverse''
(i.e., adjoint), and to ask if/when all such symmetrized (or
``Hermitianized'') word-pair-sums have positive eigenvalues. Thus, we
have the following questions, and we now state them slightly more
generally for psd matrices:

\begin{qn}\label{Qlieb}
Fix an integer $n \geq 1$. Given integers $m \geq l \geq 1$, fix an
enumeration
$\{ W_i = W_{i,m,l}(A,B) \, : \, 1 \leq i \leq \binom{m}{l} \}$ of the
$t$-degree $l$ matrix-words that occur in $(A+tB)^m$ in two
(noncommuting) square matrices $A,B$.
\begin{enumerate}
\item (Weak Lieb--Seiringer, 2012).
Classify the pairs of psd matrices $(A,B)$ such that for all $1 \leq k
\leq n$, we have $e_k \left( \sum_{i=1}^{\binom{m}{l}} W_i \right) \geq
0$ for all integers $m \geq l \geq 1$.

\item (Strong Lieb--Seiringer, 2012).
A refinement, inspired by the two bullet points above: given positive
definite $A \in \mathbb{C}^{n \times n}$, classify those $k \in [1,n]$
such that for every positive definite $B$ and integer $l \geq 1$, we have
$e_k \left( \sum_{i=1}^{\binom{m}{l}} W_i \right) \geq 0$ for all $m$
large enough. (We will answer this in Theorem~\ref{Tej}; the other two
questions here are answered in Theorem~\ref{Tcommute1}.)

\item (``Symmetrized Lieb'', $\sim$2002).
Classify the pairs of psd matrices $(A,B)$ such that for all integers $m
\geq l \geq 1$ and all $1 \leq i \leq \binom{m}{l}$, the word-sum $W_i +
W_i^*$ is psd. (Taking sums and trace would yield the BMV property.)
\end{enumerate}
\end{qn}

It appears \textit{a priori} that the third condition here is much more
restrictive than the first (all word-pairs vs.\ the sum of all words). In
fact they are equivalent, as our first result shows. It also reveals
that, remarkably, they are the same as the matrices commuting:

\begin{utheorem}\label{Tcommute1}
Fix an integer $n \geq 1$ and psd matrices $A,B \in \mathbb{C}^{n \times
n}$. The following are equivalent:
\begin{enumerate}
\item The ``weak Lieb--Seiringer'' condition in Question~\ref{Qlieb}(1)
holds.

\item The ``symmetrized Lieb'' condition in Question~\ref{Qlieb}(3) holds.

\item Every matrix $S_{m,l}(A,B)$ is positive semidefinite, for all $m
\geq l \geq 1$.

\item $AB = BA$.
\end{enumerate}
\end{utheorem}

Thus, from the BMV conjecture, we obtain multiple characterizations of
when two positive matrices commute. This suggests to call the additional
condition~(3) of the positivity of all Hurwitz products $S_{m,l}(A,B)$,
as the \textbf{strong BMV property}. We distill this property into
Definition~\ref{DstrongBMV} below, from a far stronger result (see
Theorem~\ref{Thilbert}).\footnote{{A word about the ``letter-labeled
results'': in a sense, Theorems~\ref{Thilbert} and~\ref{Tnxn} are our
``main results'', and Theorems~\ref{Tcommute1}
and~\ref{Tanderson}--\ref{Tvamsi}} are our main applications.}

We also add that a decade before Lieb, Uchiyama~\cite{Uchiyama} had
considered a question related to the ``symmetrized Lieb'' condition. In
Section~\ref{Suchi1}, we vastly generalize Uchiyama's result.

\subsection{Stronger results: fewer powers, negative powers, extension to Hilbert spaces}

Theorem~\ref{Tcommute1} is in fact one of numerous applications of our
two main results in this paper (so we have started with an application
and now get to our two results -- by successively strengthening
Theorem~\ref{Tcommute1}.)
As a first strengthening of Theorem~\ref{Tcommute1}(3): one only needs
the linearized BMV polynomials -- aka the Hurwitz products $S_{m,1}$ --
and not even all of them:

\begin{theorem}\label{Tcommute2}
The following are equivalent for $n \geq 1$ and $A,B \in \mathbb{C}^{n
\times n}$ positive semidefinite:
\begin{enumerate}
\item $S_{m,1}(A,B)$ is positive semidefinite for infinitely many $m \geq
1$.

\item $AB = BA$.
\end{enumerate}

If moreover $A,B$ are positive definite, then these are also
equivalent to the assertion:
\begin{itemize}
\item[$(1')$] $S_{m,1}(A,B)$ is positive definite for infinitely many
$m \geq 1$.
\end{itemize}

Here, the use of ``infinitely many $m$'' is best possible, in that if
$S_{m,1}(A,B)$ is psd only for finitely many $m$ then $A,B$ need not
commute.
\end{theorem}

\noindent In particular, the set of such $m$ required to deduce $AB=BA$
can be chosen to have density zero, and we also fix $l=1$. While these
strengthen the Lieb--Seiringer questions above by significantly weakening
their assumptions, the answer is the same: $AB=BA$.

The precision with which we show the final sentence of
Theorem~\ref{Tcommute2} is one reason we write Theorem~\ref{Tcommute2}
separately from both the previous result and the next ones. Another
reason is that we originally sought precisely this
Theorem~\ref{Tcommute2} -- involving the linearizations $S_{m,1}$, and
for matrices -- due to a non-BMV, geometric motivation which is described
in Appendix~\ref{sec:geom}. That question naturally leads to the terms
$S_{m,1}(A,B)$ for psd $A,B$.

We now strengthen the preceding results in numerous ways.
As mentioned above, Bessis--Moussa--Villani had already settled their
conjecture over a Hilbert space, for $A,B$ two commuting operators.
Recently, Pradhan and Skripka also worked out the BMV conjecture for
operators in Hilbert spaces \cite{PS}, using Hein\"avaara's machinery in
the finite-dimensional case~\cite{Heina}.

In this work, we too prove the above characterizations of the strong BMV
property in this setting. Thus, we first set this notation.

\begin{defn}
Throughout this paper, $(\mathcal{H}, \tangle{\cdot, \cdot})$ denotes a
nonzero complex Hilbert space, with the \textbf{global conventions}:
(a)~the inner product is $\mathbb{C}$-linear in the second argument; and
(b)~$A \geq 0$ means $A = A^* \in \mathcal{B}(\mathcal{H})$ is a
self-adjoint, positive\footnote{Recall that an operator $A$ is
self-adjoint if $\tangle{Ag, h} = \tangle{g, Ah}\ \forall g,h \in
\mathcal{H}$. Moreover, the self-adjointness of the (everywhere defined)
operators $A,B$ means that they are bounded, by Hellinger--Toeplitz.
Also, $A$ is positive if $A = B^* B$ for some $B \in
\mathcal{B}(\mathcal{H})$, and by the continuous functional calculus this
is equivalent to: $A = A^*$ and $\tangle{Ah,h} \geq 0\ \forall h \in
\mathcal{H}$.}
operator.

Now given self-adjoint $A,B \in \mathcal{B}(\mathcal{H})$, define for
integers $l \geq 0$ and $m$ (following~\eqref{Esml}):
\begin{equation}
S_{m,l}(A,B) := \left. \frac{1}{l!} \frac{d^l}{dt^l} (A+tB)^m \,
\right|_{t \to 0^+},
\end{equation}
with the \textbf{convention}:
(c)~if $m<0$ then $A$ is \textbf{positive definite} (i.e., positive and
invertible).\footnote{By open mapping / Banach--Schauder, $A^{-1}$ is
then automatically bounded.}
\end{defn}

With this, we can state our two main results, which imply commutativity
in a very general setting:
(1)~We work over any Hilbert space.
(2)~As above, we only require $S_{m,l}(A,B) \geq 0$ for fixed $l$ and
infinitely many $m$.
(3)~We allow \textit{negative} powers $m$ as well.

But a fourth, important strengthening is:
(4)~Akin to the original BMV formulation, we do not require $A$ to be
positive.
Indeed, the assumption $A,B \geq 0$ is needed in order to
deduce $AB = BA \implies S_{m,l}(A,B) \geq 0\ \forall m \geq l \geq 1$.
E.g., if $B = - {\rm id}_{\mathcal{H}}$ and $A \gneq 0$, then
$S_{m,1}(A,B) = - m A^{m-1}$, which is not positive for any $m>0$.
Similarly, if $A = - {\rm id}_{\mathcal{H}}$ and $B^l \geq 0$, then
$S_{m,l}(A,B) = -\binom{m}{l} B^l \leq 0$ for all $m$ of the opposite
parity to $l$.

It is natural to wonder if the \textit{reverse} implication also depends
on the positivity of $A,B$. Surprisingly, we reveal that it does not!
But this comes at a price: one obtains not $AB=BA$, but something weaker
-- which is, in a sense, the ``true'' underlying commutativity:

\begin{utheorem}\label{Thilbert}
Suppose $\mathcal{H}$ is a Hilbert space and $A,B \in
\mathcal{B}(\mathcal{H})$ are self-adjoint.
\begin{enumerate}
\item Suppose $l = 1$. If $S_{m,l}(A,B) \geq 0$ for infinitely many $m
\geq l$ (resp.\ ${\rm sgn}(m)^l \cdot S_{m,l}(A,B) \geq 0$ for infinitely
many $m \in \mathbb{Z}$ if $A>0$), then $A^2B = BA^2$.

\item If $l \geq 2$, then the same assertion holds if $B \geq 0$ and the
spectrum $\sigma(A) \setminus \{ 0 \}$ is
discrete.\footnote{It is standard that if $\sigma(A) \setminus \{ 0 \}$
is discrete, then $A = A^* \in \mathcal{B}(\mathcal{H})$ has a complete
orthonormal eigenbasis for $\mathcal{H}$ -- see Lemma~\ref{Leigen}.
Below, we will use this eigenbasis without further mention.}

\item If $A^2B = BA^2$ and $\sigma(A) \cap \sigma(-A) \subseteq \{ 0 \}$
(e.g.\ if $\pm A \geq 0$), then $AB = BA$.
\end{enumerate}
\end{utheorem}

As in the matrix case, working with finitely many $S_{m,l}$ does not
yield commutativity:

\begin{theorem}\label{Tfinite}
Fix a complex Hilbert space $\mathcal{H}$ of dimension in $[2,\infty]$,
and finitely many integer pairs $(m_1, l_1), \dots, (m_n, l_n)$ with all
$l_k \geq 0$.
Then there exist bounded self-adjoint positive operators $A,B \in
\mathcal{B}(\mathcal{H})$ such that $A^2B \neq BA^2$ (whence $AB \neq
BA$), and yet for every word $W$ occurring in any $\sgn(m_k)^l \cdot
S_{m_k,l_k}(A,B)$, the symmetrization $W+W^*$ is positive (whence so is
each $\sgn(m_k)^l \cdot S_{m_k,l_k}(A,B)$).
\end{theorem}

These results subsume Theorems~\ref{Tcommute1} and~\ref{Tcommute2}
(except for Theorem~\ref{Tcommute1}(1), which does not have an
infinite-dimensional counterpart in general). In particular, one can also
work with negative $m$ for positive definite matrices $A,B$ above.

\begin{remark}
We provide a significant strengthening of Theorem~\ref{Thilbert}(1)
below, valid for any self-adjoint $A,B$ over any Hilbert space. See
Theorem~\ref{Tadmissible}.
\end{remark}

Our next main result obtains the original commutation $AB=BA$. As in
Theorem~\ref{Thilbert}(1), we use ${\rm sgn}(m)^l$ to unify the
statements for $m \geq l$ and $m<0$.

\begin{utheorem}\label{Tnxn}
Suppose $\mathcal{H}$ is a Hilbert space and $A,B \in
\mathcal{B}(\mathcal{H})$, with $B \geq 0$, $A = A^*$, and $\sigma(A)
\setminus \{ 0 \}$ discrete. Then the following are equivalent for any $l
\geq 1$:
\begin{enumerate}
\item ${\rm sgn}(m)^l \cdot S_{m,l}(A,B) \geq 0$ for all $m \geq l$
(resp.\ all $m \in \mathbb{Z}$ if $A>0$).

\item ${\rm sgn}(m)^l \cdot S_{m,l}(A,B) \geq 0$ for infinitely many $m
\geq l$ with $m-l$ odd (resp.\ infinitely many $m < 0$ with $m$ even, if
$A>0$).

\item $AB=BA$ and $P_{\mathcal{H}_-} B P_{\mathcal{H}_-} = 0$
(equivalently, $B P_{\mathcal{H}_-} = B|_{\mathcal{H}_-} = 0$). Here,
$\mathcal{H}_-$ is the (closed) orthogonal sum of all
negative-eigenspaces $\ker(A - \lambda \, {\rm id}_{\mathcal{H}})$ for $0
> \lambda \in \sigma(A)$.
\end{enumerate}
\end{utheorem}

\begin{remark}
For $l=1,2$, we explain in Sections~\ref{S1tight} and~\ref{S2tight}
respectively, how our results are ``tight'' in ensuring $AB = BA$.
\end{remark}

\subsection{$C^*$-algebras}

Many of our results hold over an arbitrary $C^*$-algebra $\mathcal{A}$
(the prototypical example of this is $\mathcal{B}(\mathcal{H})$). Thus,
we recast the above results in this setting:

\begin{theorem}\label{Trecast}
Suppose $\mathcal{A}$ is a complex $C^*$-algebra. Fix $A = A^*, B = B^*$
in $\mathcal{A}$.
\begin{enumerate}
\item Theorem~\ref{Thilbert} holds verbatim in $\mathcal{A}$.

\item If $B \geq 0$ and $\sigma(A) \setminus \{ 0 \}$ is discrete, then
Theorem~\ref{Tnxn}(1) and Theorem~\ref{Tnxn}(2) are equivalent and imply
$AB = BA$.
\end{enumerate}
\end{theorem}

\begin{proof}
Indeed, by Gelfand--Naimark--Segal (GNS) we have an isometric
$*$-preserving algebra map $\pi : \mathcal{A} \to
\mathcal{B}(\mathcal{H})$ (hence injective). So $0 \leq \pi(S_{m,l}(A,B))
= S_{m,l}(\pi(A), \pi(B))$. Now e.g.\ for~(1), Theorem~\ref{Thilbert}
yields $0 = [\pi(A)^2, \pi(B)] = \pi([A^2,B])$, whence $0 = [A^2,B]$.
Similarly, use Theorem~\ref{Tnxn} to show~(2).
\end{proof}

\begin{remark}
In the preceding proof, and below, $[p,q]$ denotes the commutator $pq-qp$
in any associative ring (in fact $C^*$-algebra).
\end{remark}

We now provide a ``first application'' of our previous two main theorems:
our original Theorem~\ref{Tcommute1} for matrices is a special case,
which holds far more generally:

\begin{corollary}\label{Ccstar}
Given any $C^*$-algebra $\mathcal{A}$ and any $0 \leq a,b \in
\mathcal{A}$ with $\sigma(a) \setminus \{0\}$ discrete, the following are
equivalent:
\begin{enumerate}
\item For some fixed $l_0 \geq 1$, the Hurwitz product ${\rm sgn}(m)^l
\cdot S_{m,l_0}(a,b) \geq 0$ for infinitely many $m \geq l_0$ (or for
infinitely many $m < 0$ if $a \geq 0$ is invertible).

\item $W_i(a,b) + W_i(a,b)^* \geq 0$ for all $W_i$, for all $m \geq l
\geq 1$ (or for all $m\in \mathbb{Z}$ if $a \geq 0$ is invertible), where
$W_i = W_{i,m,l}$ are the words occurring in ${\rm sgn}(m)^l \cdot
S_{m,l}$.

\item ${\rm sgn}(m)^l \cdot S_{m,l}(a,b) \geq 0$ for all $m \geq l \geq
1$ (or for all $m\in \mathbb{Z}$ if $a \geq 0$ is invertible).

\item $ab=ba$ in $\mathcal{A}$.
\end{enumerate}
\end{corollary}

\begin{proof}
Even without assuming $a,b \geq 0$, one clearly has $(2) \implies (3)
\implies (1)$. If next~(1) holds, then by Theorem~\ref{Thilbert} we get
$a^2b = ba^2$.

Now choose a sequence of polynomials $p_k(x)$ converging
uniformly on $\sigma(a^2) \subseteq [0,\|a^2\|]$ to $\sqrt{x}$. Thus,
$p_k(a^2) b = b p_k(a^2)$ for all $k$. Taking the limit, the continuous
functional calculus -- see Remark~\ref{Rfnalcalc} in
Appendix~\ref{Sspectral} -- implies $ab=ba$ (one can either work in
$\mathcal{A}$, or embed $\mathcal{A} \hookrightarrow
\mathcal{B}(\mathcal{H})$ and work over $\mathcal{H}$).

Similarly, if $ab=ba$ for $a,b \geq 0$ then every word $W_i(a,b)$ inside
$S_{m,l}(a,b)$ equals $a^{m-l}b^l$ for $m \geq 0$. Since $a$ commutes
with $b$, so does $a^{1/2}$ by a similar argument to above. Hence $W_i +
W^*_i = 2 a^{m-l} b^l = w^* w \geq 0$, where $w = \sqrt{2} \, b^{l/2}
a^{(m-l)/2} \geq 0$. If $m<0$, we use Proposition~\ref{Psml} below.
\end{proof}

Theorem~\ref{Trecast} (following Theorem~\ref{Thilbert}) motivates the
following notion:

\begin{defn}\label{DstrongBMV}
Let $a,b$ be self-adjoint operators in a complex $C^*$-algebra
$\mathcal{A}$. We say the ordered pair $(a,b)$ has the \textbf{strong BMV
property of order $l \geq 1$} if $S_{m,l}(a,b) \geq 0$ for infinitely
many $m \geq l$.
\end{defn}

\noindent This implies $[a^2,b]=0$ if $l=1$, and also for $l \geq 2$ if
say $b \geq 0$ and $a$ has discrete nonzero spectrum, by
Theorem~\ref{Thilbert}.

\begin{remark}
While we impose a discreteness hypothesis on $\sigma(A)$ if $l \geq 2$,
this setting nevertheless subsumes all self-adjoint compact operators and
all projection operators.
\end{remark}

\begin{remark}\label{Rdichotomy}
The equivalence of $S_{m,l} \geq 0$ ``eventually'' and for infinitely
many $m \geq l$, reveals a sharp dichotomy. Namely, suppose $l \geq 1$
and $0 \leq A,B \in \mathcal{B}(\mathcal{H})$ (with $A$ as in
Theorem~\ref{Tnxn}). Then the sequence of Hurwitz products $S_{l,l}(A,B),
\ S_{l+1,l}(A,B), \ \dots$ is eventually always positive or eventually
not positive (but still with nonnegative trace when $A,B$ are matrices).

Also, as Bessis--Moussa--Villani remark on \cite[p.~2319]{BMV}: commuting
operators correspond to classical mechanics, and so
Theorems~\ref{Thilbert} and~\ref{Tnxn} are ``negative quantum results''.
\end{remark}

It is worth isolating the finite-dimensional version of these results:

\begin{corollary}
Theorems~\ref{Thilbert} and~\ref{Tnxn} hold for $n \times n$ matrices
$A,B$, for every $n \geq 1$.
\end{corollary}

For completeness, we also address the real case. The results for real psd
matrices follow immediately from the above results for Hermitian
matrices. As for the Hilbert space results:

\begin{theorem}\label{Treal}
Theorems~\ref{Thilbert} and~\ref{Tnxn} also hold for real Hilbert spaces
-- and in particular, for real psd matrices. So do the other theorems
above.
\end{theorem}

\noindent As this is standard, we defer a quick proof to
Appendix~\ref{Appendix3}.

A quick note on the proof-techniques. As we explain below, the methods
used to show Theorems~\ref{Tcommute1} and~\ref{Tcommute2} for matrices do
not work to show the Hilbert space theorems~\ref{Thilbert}(2)
and~\ref{Tnxn}, not even for $\mathcal{H} = \mathbb{C}^2$. Thus, we
provide alternate arguments that are akin to the ``analytic
combinatorics'' found in e.g.\ \cite[Chapters IV--VI]{FS} -- and also
similar to an extent with some arguments in Campanino and Klein's
work~\cite{CK}, wherein they prove Simon and Taylor's
conjecture~\cite{ST}. There are additional challenges that do not allow
us to simply adapt the $l=1$ arguments; we describe these in
Section~\ref{Shilbert2}. Moreover, as $B$ need not have eigenvalues, we
use spectral projections instead. (This is also done in proving
Theorem~\ref{Thilbert}(1),(3).)

\subsection{Anderson's result on Weyl quantization}\label{Sweylintro}

We next present applications of our results and techniques (in addition
to answering Lieb and Seiringer's questions above). Here we discuss the
Weyl functional calculus; this was introduced by Weyl in his 1928
monograph~\cite{Weyl} as a way to act by a multivariate function on a
tuple of \textit{non-commuting} operators on a Hilbert space. A special
case is the \textit{Wigner distribution}~\cite{Wigner}, extensively used
in quantum optics, quantum field dynamics, statistical mechanics, and
other areas. The original, famous instance of noncommuting operators are
of course the position and momentum operators in quantum mechanics:
$[p,q] \neq 0$, in which case this is Weyl--Wigner quantization; see
also~\cite{McCoy}.

Briefly, let $\mathcal{H}$ be a real or complex Hilbert space, and $A,B
\in \mathcal{B}(\mathcal{H})$ self-adjoint operators. If $AB=BA$ then one
can apply a (continuous/holomorphic/etc.) bivariate function $f$ to them
via the Spectral theorem, to obtain $f(A,B)$. But for general $(A,B)$, a
workaround is required, and a natural choice is the Weyl calculus. This
is an operator valued-distribution $\mathcal{W}_{A,B} :
C^\infty(\mathbb{R}^2) \to \mathcal{B}(\mathcal{H})$, which is uniquely
determined by being linear, ``continuous'', and sending $x^{m-l} y^l$ to
$S_{m,l}(A,B) / \binom{m}{l}$ for all integers $m \geq l \geq 0$. If
$A,B$ commute then $\mathcal{W}_{A,B}$ is indeed their joint spectral
measure. Some more details -- for a $d$-tuple of operators -- are
provided in Appendix~\ref{Sweylappendix}; we also cite the comprehensive
monograph~\cite{Jeff}.

Now Anderson's result~\cite[Theorem~3.7]{And1}, states:
\textit{If for self-adjoint operators $A_1, \dots, A_d \in
\mathcal{B}(\mathcal{H})$ the Weyl distribution $\mathcal{W}_{A_1, \dots,
A_d}$ is a positive measure, then the $A_j$ commute pairwise.}

Notice that the positivity condition on $\mathcal{W}_{A_1, \dots, A_d}$
is very strong. For instance, if e.g.\ $d=2$ and $A = A_1, B = A_2$ are
positive, it contains the information that
\[
0 \leq \mathcal{W}_{A,B} (x^{m-l} y^l) = S_{m,l}(A,B) / \textstyle
\binom{m}{l}, \quad \forall m \geq l \geq 0,
\]
and from this information one deduces $AB=BA$. By contrast,
Theorem~\ref{Tnxn} deduces the same equality from merely a single $l \geq
1$ (i.e., a one-dimensional set of hypotheses instead of 2-dim.) -- and
not all $m$ but only infinitely many of them. (As the reader will know,
the positivity of infinitely many moments of a measure is weaker than
that of all moments.)

We thus show the following result, whose hypotheses are significantly
weaker than Anderson's but give the same conclusion:

\begin{utheorem}[Weyl calculus]\label{Tanderson}
Let $d \geq 2$, $\mathcal{H}$ be a complex Hilbert space, and $A_1,
\dots, A_d \in \mathcal{B}(\mathcal{H})$ be self-adjoint, and suppose for
all $1 \leq i < j \leq d$ that $\mathcal{W}_{A_1, \dots,
A_d}(p_i(x_i)(x_j+c_j))$ is a positive operator, for all polynomials $p_j
\in \mathbb{R}[t]$ with $p_j(A_j), A_j + c_j \, {\rm id}_{\mathcal{H}}
\geq 0\ \forall j \geq 1$. Then the $A_j$ commute pairwise.
\end{utheorem}

As our proof reveals, we only need one value of $c_j$ for each $j$: any
scalar in $[\|A_j\|,\infty)$.

\subsection{Commuting via operator complete monotonicity, via Gibbs
thermal weights, and Moln\'ar--Wu}\label{S1agniva}

This part returns to the original formulation of the BMV conjecture, but
on the level of operators. Recall that a \textit{completely monotone
function} is a continuous map $f : [0,\infty) \to \mathbb{R}$ that is
smooth on $(0,\infty)$ and satisfies: $(-1)^m f^{(m)} \geq 0$ on
$(0,\infty)$ for all $m \geq 0$. A celebrated paper of
Bernstein~\cite{Bernstein} (see also~\cite{Widder}) shows these are
precisely the functions that are Laplace transforms of positive measures,
aka a ``mixture of exponentials.''

In this subsection we present three results. The first comes from recent
work~\cite{Agniva} of Bhowmik--Chatterjee--Putinar, where they studied
completely monotone functions in several variables, defined on cones and
applied to Hermitian matrices $A_1, \dots, A_d$ via the Loewner
order.\footnote{Recall the Loewner order, in fact on self-adjoint
operators in $\mathcal{B}(\mathcal{H})$ for any $\mathcal{H}$: $A \geq B$
means $A-B$ is positive.} In other words, they studied the complete
monotonicity of the BMV partition function -- or more precisely,
$\widetilde{Z}(s,t) := e^{-(sA+tB)}$ -- without taking traces.

\begin{defn}
Fix $d \geq 1$. A \textit{cone} $\Gamma$ in $\mathbb{R}^d$ is a subset
closed under rescaling by $\mathbb{R}_+ = (0,\infty)$. We say a closed
cone $\Gamma = \overline{\Gamma}$ is \textit{acute} if $\tangle{\bxi,x}
\geq 0$ for all $\bxi,x \in \Gamma$, and \textit{solid} if its span is
$\mathbb{R}^d$.

The \textit{dual cone} to $\Gamma$ is $\Gamma^* := \{ y \in \mathbb{R}^d
: \tangle{\bxi,y} \geq 0\ \forall \bxi \in \Gamma \}$.
Given $\bxi = (\xi_1, \dots, \xi_d) \in \Gamma$, define the directional
derivative $D_{\bxi} := \sum_{j=1}^d \xi_j \frac{\partial}{\partial
x_j}$. 
\end{defn}

\begin{defn}\label{DmatrixCM}
A map $\varphi : \Gamma \to \mathbb{C}^{N \times N}$ that is
smooth\footnote{In~\cite{Agniva} the authors worked with smooth maps into
$\mathbb{R}^{N \times N}$, but the results do not change here.}
on $\inte \Gamma$ (with all derivatives extended to the boundary by
continuity), is said to be \textit{matrix-valued completely monotone on
$\Gamma$} if
\begin{equation}\label{EmatrixCM}
(-1)^m (D_{\bxi_1} \cdots D_{\bxi_m} \varphi)(x) \geq {\bf 0}_{N \times
N}, \qquad \forall m \geq 0, \ x \in \inte \Gamma, \ \bxi_1, \dots,
\bxi_m \in \Gamma.
\end{equation}
\end{defn}

In~\cite{Agniva}, the authors study not just completely monotone maps --
aka Laplace transforms of positive (finite) measures $\mu$ on the dual
cone $\Gamma^*$ -- but also Stieltjes--Fantappi\'e transforms of such
measures:
\begin{equation}\label{Eagniva}
f(z) = \int_{\Gamma^*} \exp(- z \cdot y) d \mu(y), \qquad
g(z) = \int_{\Gamma^*} \frac{d \mu(y)}{z_0 + z \cdot y}.
\end{equation}
(These two operations are related by the Borel transform.)
Here, $T_\Gamma := \inte \Gamma + i\mathbb{R}^d \subset \mathbb{C}^d$ and
$T_{\Gamma \times \mathbb{R}_+} := (\inte \Gamma \times
\mathbb{R}_+) + i\mathbb{R}^{d+1} \subset \mathbb{C}^{d+1}$ are tube
domains, and $z \in T_\Gamma, (z,z_0) \in T_{\Gamma \times
\mathbb{R}_+}$.

The motivation behind the work~\cite{Agniva} was to interpolate and
approximate these integrals~\eqref{Eagniva} by finite sums. Note that in
most of~\cite{Agniva} these transforms are scalar-valued -- i.e.\ $y \in
\Gamma^*$ -- but the authors also explore the matrix-valued versions in
their Section~7.

Our next result explains how these constructions extend to
operator-valued maps. As the following result shows, it suffices to work
with only ``local'' complete monotonicity at the origin, which unifies
the two distinguished settings~\eqref{Eagniva} in~\cite{Agniva}:

\begin{utheorem}[Operator-valued complete monotonicity implies commuting]\label{Tagniva}
Suppose $\epsilon > 0$ and $f \in C^\infty(-\epsilon,\epsilon)$ has a
power series expansion $\sum_{m=0}^\infty c_m t^m$ with all $(-1)^m c_m
\geq 0$ and $c_m \neq 0$ for infinitely many $m>0$, not all of the same
parity. Also fix a Hilbert space $\mathcal{H}$, self-adjoint operators
$A_1, \dots, A_d \in \mathcal{B}(\mathcal{H})$ for $d \geq 1$, and a
closed convex acute solid cone $\Gamma \subset \mathbb{R}^d$.

Then each part implies the next.
\begin{enumerate}
\item For a small ball of radius $0 < \epsilon_0 \ll 1$, define the map
$\varphi_{\bf A} : \Gamma \cap B({\bf 0}, \epsilon_0) \to
\mathcal{B}(\mathcal{H})$ via: $x \mapsto f(x \cdot {\bf
A})$.\footnote{\textit{Notation:} Here and below, the tuple ${\bf A} :=
(A_1, \dots, A_d)$, and one denotes $x \cdot {\bf A} := \sum_{j=1}^d x_j
A_j$ for $x \in \mathbb{R}^d$.}
Also fix linearly independent $\bxi_1, \dots, \bxi_d \in \Gamma$. Then
for all $1 \leq i < j \leq d$, we have {\em operator-valued complete
monotonicity} (in 2 variables) at the origin:
\[
(-1)^m (D_{\bxi_i}^{m-1} D_{\bxi_j} \varphi_{\bf A})({\bf 0}) \geq 0
\]
for infinitely many $m>0$, not all of the same parity, such that $c_m
\neq 0$.

\item The operators $A_j$ commute pairwise.

\item If ${\bf A}$ is moreover in the {\em dual
$\mathcal{B}(\mathcal{H})$-cone} -- which is defined to be
\begin{equation}
\Gamma^*_{\mathcal{B}(\mathcal{H})} := \{ (A_1, \dots, A_d) \in
\mathcal{B}(\mathcal{H})^d \ | \ A_j = A_j^*, \ \sum_{j=1}^d x_j A_j \geq
{\bf 0}\ \forall x \in \Gamma \},
\end{equation}
then $\sigma_{joint}({\bf A}) \subseteq \Gamma^*$.
\end{enumerate}
\end{utheorem}

See Definition~\ref{Djointspec} for the joint spectrum
$\sigma_{joint}({\bf A})$ of the commuting operators $A_j$. Also note
that the above definition of the dual $\mathcal{B}(\mathcal{H})$-cone
$\Gamma^*_{\mathcal{B}(\mathcal{H})}$ is rather natural in this setting,
and if $\mathcal{H} = \mathbb{C}$ then
$\Gamma^*_{\mathcal{B}(\mathcal{H})}$ specializes to $\Gamma^*$, the dual
cone.

\begin{remark}
Note that the class of maps above includes not only the ``Laplace'' and
``Stieltjes--Fantappi\'e'' maps $\exp(-yt)$ and $(y_0 + yt)^{-1}$ for
$y_0,y>0$ as in~\eqref{Eagniva}, but also other completely monotone-type
maps: $-\log(1+yt)$ or $(1+yt)^{-\gamma}$ for $\gamma,y>0$, or even
$g_y(t) := 1 - t + t^{20} \exp(-yt)$ for $y>0$, among others. Moreover,
$g_y^{(19)}(t)$ is the power series $20! \, t + t^2 p(t)$ for some power
series (which is absolutely convergent on $\mathbb{R}$). It follows that
$g_y^{(19)}(t) > 0$ for small $0 < t \ll 1$. Hence this class of maps is
strictly larger than the (locally) completely monotone maps.
\end{remark}

\begin{remark}
Theorem~\ref{Tagniva} strengthens the main result of
\cite[Section~7]{Agniva} in numerous ways at the same time. Firstly, it
``decouples'' the implications $(1) \implies (2)$ and $(2) \implies (3)$.
Second, to deduce commutativity in part~(2), the assertion in part~(1)
does not assume (as in~\cite{Agniva})
\begin{itemize}
\item either that the $\bxi_j$ (which were called $v^{(j)}$
in~\cite{Agniva}) lie on extreme rays in $\Gamma$,
\item or that the $\bxi_j \cdot {\bf A}$ (which were called
$L(v^{(j)})$ in~\cite{Agniva}) are positive operators,
\item or that $S_{m,1}$ -- which was called $C_{k,1}$ in~\cite{Agniva} --
is positive for every $m$. (Indeed, we only need to assume it is positive
for a zero-density infinite set of $m$, not all even/odd.)
\end{itemize}
Finally, Theorem~\ref{Tagniva} is valid for operators over any Hilbert
space, and not merely matrices.
\end{remark}

The next result returns to the quantum-physics origins of the BMV
conjecture, again without taking traces:

\begin{utheorem}[Commuting via Gibbs thermal weights]\label{Tthermal}
Suppose $\mathcal{H}$ is a Hilbert space, and $H_0, H_1 \in
\mathcal{B}(\mathcal{H})$ are self-adjoint.
\begin{enumerate}
\item Then $H_0, H_1$ commute if the {\em Gibbs thermal weight operator}
$e^{-\beta(H_0 + t H_1)}$ is monotone non-increasing at $H_0$ along $H_1$
(in the Loewner order):
\[
\left. -\frac{d}{dt} e^{-\beta(H_0 + t H_1)} \right|_{t=0} \geq 0,
\]
along some unbounded sequence of inverse-temperatures $0 < \beta_k
\uparrow \infty$.
\item Conversely, if $H_0, H_1$ commute and moreover $H_1 \geq 0$,
then $e^{-\beta(H_0 + t H_1)}$ is monotone non-increasing for all $\beta
> 0$.
\end{enumerate}
\end{utheorem}

We deduce the harder part, Theorem~\ref{Tthermal}(1), from a more general
result. See Theorem~\ref{TC2}.

\begin{remark}
The BMV conjecture has well-known consequences in quantum mechanics. Now
Theorems~\ref{Thilbert} and~\ref{Tnxn} study a stronger phenomenon: the
positivity of the operator $S_{m,l}(A,B)$, which captures the complete
quantum interference obtained from $l$ interactions and $m-l$ free
propagations, across all time-ordered paths. These theorems translate
into the following quantum-physics result in a setting where (unlike
classical mechanics) the two Hamiltonians $A,B \geq 0$ do not commute. If
so, then for every order $l \geq 1$, the perturbation terms
$S_{m,l}(A,B)$ (e.g.\ in $\widetilde{Z}(s,t) := e^{sA+tB}$) stop being
(dilations of) positive unnormalized states beyond some degree $m \geq
m_l > 0$.
\end{remark}

\begin{remark}
Theorem~\ref{Tthermal} reveals another (see Remark~\ref{Rdichotomy})
dichotomy: fix any $\mathcal{H}$ and operators $H_0 = H_0^*$ and $H_1
\geq 0$. Then the operator derivative $\left. -\frac{d}{dt} e^{-\beta(H_0
+ t H_1)} \right|_{t=0}$ is either positive for every $\beta > 0$, or
not positive for every sufficiently large $\beta \gg 0$.
\end{remark}

As consequences of Theorem~\ref{Tthermal} we obtain multiple results in
the literature that relate commuting to Loewner monotonicity of maps --
the main theorems of Moln\'ar~\cite{Molnar} (actually, a weakening) and
Wu~\cite{Wu}. The latter says that a $C^*$-algebra is commutative if and
only if $\exp$ is Loewner-monotone. This is similar in spirit to a
well-known result of Ogasawara~\cite{Ogasa}, which uses the square in
place of the exponential -- and another by Pedersen~\cite{Pedersen}, who
used the $p$th power map for $p>1$. We also obtain both of these results
(for $p \geq 2$) as consequences of a different main result (see
Theorem~\ref{Tuchi}).

\begin{utheorem}[Commutativity via Loewner monotonicity]\label{Tmolnar}
Suppose $\mathcal{A}$ is a unital $C^*$-algebra.
\begin{enumerate}
\item {\rm (``Weak Moln\'ar''.)}
A (fixed) element $H_0 = H_0^* \in \mathcal{A}$ is central (i.e.,
commutes with all of $\mathcal{A}$) if and only if $e^{\beta H} \geq
e^{\beta H_0}$ for all $H = H^* \geq H_0$ in $\mathcal{A}$ and $\beta>0$
(equivalently, for one unbounded sequence $0 < \beta_k \uparrow \infty$).

\item {\rm (Wu.)}
$\mathcal{A}$ is commutative if and only if $\exp$ is Loewner-monotone,
i.e.\ $H \geq H_0$ for self-adjoint $H_0, H$ implies $e^H \geq e^{H_0}$.

\item {\rm (Ogasawara, and the Moln\'ar-variant.)}
An element $0 \leq H_0 \in \mathcal{A}$ (now not necessarily unital) is
central in $\mathcal{A}$ if and only if $H^2 \geq H_0^2$ whenever $H \geq
H_0 \geq 0$ in $\mathcal{A}$. Thus, $\mathcal{A}$ is commutative if and
only if the square-map is Loewner-monotone on positive operators.

\item {\rm (Pedersen, and the Moln\'ar-variant.)}
Fix $p \geq 2$. An element $0 \leq H_0 \in \mathcal{A}$ (again not
necessarily unital) is central in $\mathcal{A}$ if and only if $H^p \geq
H_0^p$ whenever $H \geq H_0 \geq 0$. Thus, $\mathcal{A}$ is commutative
if and only if the $p$th power map is Loewner-monotone on positive
operators.
\end{enumerate}
\end{utheorem}

See Sections~\ref{Sagniva} and~\ref{Suchi} for the proofs. Also, it is
not hard to see (and our proof will do so) that even the ``weak
Moln\'ar'' condition is much stronger than our hypothesis in
Theorem~\ref{Tthermal}; yet \textit{a posteriori}, they are both
equivalent to the commutation of $H_0, H_1$.

\subsection{Strengthening the BMV conjecture to monomial-positive
symmetric functions}

Note that all of our results so far were ``conditional'', i.e.\ involving
the positivity of Hurwitz products, which means commuting of $A,B$. Our
next application is ``unconditional'' -- i.e., for all $A,B \geq 0$ --
and goes further than the original BMV conjecture itself. In particular,
we answer the ``Strong Lieb--Seiringer'' Question~\ref{Qlieb}(2).

Recall that our original goal was to understand Lieb and Seiringer's
questions on the positivity of $S_{m,l}(A,B)$ or of $W+W^*$ (and we were
also motivated by geometric considerations, see Appendix~\ref{sec:geom}).
Specifically:
\begin{itemize}
\item We answered Question~\ref{Qlieb}(1) for each fixed $l \geq 1$: when
is $e_k(S_{m,l}(A,B)) := e_k(\sigma(S_{m,l}(A,B)))$ nonnegative for all
$1 \leq k \leq n$ and $m \geq l$. The answer is ``conditional'': $AB=BA$.

\item This is close to, but disparate from, the BMV conjecture / Stahl's
theorem -- which says that $e_1(S_{m,l}(A,B)) = {\rm Tr}(S_{m,l}(A,B))
\geq 0$ ``unconditionally'' for all $A,B \geq 0$.
\end{itemize}

Given these, a natural strengthening of the BMV conjecture would hold if
at least the following question has a positive answer:\smallskip

\begin{center}
\textit{Is it true that $e_2(S_{m,l}(A,B)) \geq 0$ for all $A,B \geq 0$
and $m \geq l \geq 1$?}
\end{center}\smallskip

\noindent It turns out this is false -- and badly so -- already for the
smallest values possible:

\begin{example}\label{Exeventual}
Even $e_2(S_{2,1}(A,B))$ can be negative for $3 \times 3$ matrices $A,B$.
More generally, let $l = 1$, $n \geq 3$, and $a \geq 8 \binom{n}{2}$, and
set
\begin{equation}\label{EspecialB}
A = \begin{pmatrix} a & 0 \\ 0 & 1 \end{pmatrix}
\oplus {\rm Id}_{n-2}, \qquad
\mathscr{B} := \begin{pmatrix} 2 & 1 \\ 1 & 1 \end{pmatrix} \oplus {\rm
Id}_{n-2}.
\end{equation}

We claim that $e_2(S_{m,1}(A,\mathscr{B}))$ is far from being always
nonnegative -- in fact, it is \textit{always negative} for $m \geq 2$ --
fully opposite to the BMV theorem for $e_1$! The proof is deferred to
Appendix~\ref{Appendix3}. \qed
\end{example}

Given this, the ``next best hope'' would be to extend the BMV conjecture
as follows:\smallskip

\begin{center}
\textit{Fix arbitrary $A > 0$. Which $e_k(S_{m,l}(A,B))$ are guaranteed
to be nonnegative,}

\textit{for all $B > 0$, $l \geq 1$, and (all, or at least) infinitely
many $m \geq l$?}
\end{center}\smallskip

This is precisely the ``Strong Lieb--Seiringer'' condition, see
Question~\ref{Qlieb}(2). We provide a complete -- and tight -- answer to
this question:

\begin{theorem}\label{Tej}
Given a Hermitian matrix $M_{n \times n}$ (with spectrum $\sigma(M)
\subset \mathbb{R}$), let $d_M := \dim \ker(M - \max(\sigma(M)) {\rm
Id}_n)$ denote the dimension of the eigenspace for the largest
eigenvalue.

Then eventual BMV-positivity holds for all $k \leq d_A$: given $n,l \geq
1$ and positive definite $A,B \in \mathbb{C}^{n \times n}$, there is a
threshold $m_l \geq l$ such that
\[
e_k(S_{m,l}(A,B)) > 0, \qquad \forall 1 \leq k \leq d_A,
\]
for all $m \geq m_l$. This is tight in two ways:
(a)~$d_A$ is sharp, and
(b)~``eventual-BMV'' $m \gg 0$ cannot be improved to $m \geq l$ as in the
BMV conjecture, for any $k \geq 2$.
\end{theorem}

In particular, this also explains why Example~\ref{Exeventual} fails:
$d_A = 1$ here, and we consider $e_2$ with $AB \neq BA$.

\begin{defn}
The tightness of (a)~$d_A$ in Theorem~\ref{Tej} thus refines the original
BMV conjecture, for each matrix $A$ individually. Given this, we term the
\textbf{BMV-degree} of a positive definite matrix $A$ to be $d_A$.
\end{defn}

\begin{remark}
There was another ``tightness'' in Theorem~\ref{Tej}: of (b)~the fact
that the threshold $m_l$ cannot be brought down to $l$ for $e_2, e_3,
\dots$. Thus, the original BMV conjecture/theorem -- for $e_1$ -- is
indeed special, because in this case we have $e_1(S_{m,l}) \geq 0$ for
all $m \geq 1$, not merely for ``all large $m$'' (i.e., universal and not
just eventual).
\end{remark}

We now extend Theorem~\ref{Tej} well beyond elementary symmetric
polynomials $e_k$, via the theory of symmetric functions. We first
describe the successive strengthenings, before stating our result.
\begin{itemize}
\item Recall that a \textit{partition} $\lambda = (\lambda_1 \geq
\lambda_2 \geq \cdots)$ is a non-increasing sequence of nonnegative
integers with finite sum, e.g.\ the $k$-tuple of ones $(1,\dots,1) =
(1)^k$ for $k \geq 1$.
By the \textit{number of (positive) parts}, or \textit{length} of
$\lambda$, we  mean $\max \{ k \geq 0 : \lambda_k > 0 \}$ -- and the
empty/zero partition has length $0$.
For each partition one has several families of
symmetric functions, such as the monomial symmetric functions
$m_\lambda(x_1, x_2, \dots) = \sum_{\tau \, \text{a permutation}}
\prod_{k \geq 1} x_{\tau(k)}^{\lambda_k}$ (where one sums over only the
distinct monomials), or Schur polynomials $s_\lambda({\bf x})$. Both of
these generalize the elementary symmetric polynomials: $e_k({\bf x}) =
s_{(1)^k}({\bf x}) = m_{(1)^k}({\bf x})$.

\item The Schur polynomials $s_\lambda$ are a $\mathbb{Z}_{\geq
0}$-linear combination of monomial functions $m_\mu$ (over all partitions
$\mu$ majorized by $\lambda$ -- there are only finitely many such $\mu$
for any $\lambda$). One can work more generally over the function field
$\mathbb{Q}(\alpha)$; now the Schur polynomials are a special case of
\textit{Jack polynomials} $J^{(\alpha)}_\lambda({\bf x})$. These were
defined by Jack~\cite{Jack}, and generalize zonal polynomials
$Z_\lambda({\bf x})$ (which have connections to
$GL_n(\mathbb{R})$-representations and multivariate statistics), Schur
polynomials $s_\lambda({\bf x})$, and monomial symmetric functions
$m_\lambda({\bf x})$. Moreover, they are eigenfunctions of the
Calogero--Moser--Sutherland operator in mathematical physics~\cite{SV}.

Every Jack polynomial $J_\lambda^{(\alpha)}({\bf x})$ is again a
$\mathbb{Q}(\alpha)$-linear combination of the monomials $m_\mu({\bf x})$
(with $\mu$ majorized by $\lambda$). The coefficient of each $m_\mu$ in
$J_\lambda^{(\alpha)}$ is a rational function in $\alpha$ with
coefficients in $\mathbb{Z}_{\geq 0}$, and we normalize the denominators
such that $m_\lambda$ has coefficient $1$ -- see~\cite{KS}, and e.g.\
\cite[Chapter~VI.10]{Macdonald} for more details.

\item Even more general are the \textit{Macdonald polynomials}
$P_\lambda^{(q,t)}({\bf x})$, where $q,t \in (0,1)$ -- see
\cite[Chapter~VI]{Macdonald}. Once again, for each fixed $q,t \in (0,1)$,
$P_\lambda^{(q,t)}({\bf x})$ is an $\mathbb{R}_{\geq 0}$-linear
combination of the same $m_\mu({\bf x})$ as above -- and the coefficient
of $m_\lambda({\bf x})$ is $1$ under a certain normalization.
\end{itemize}

Each family above yields the preceding family via specializing some
parameters -- and to get Jack polynomials, taking a limit. Thus, we
begin with a significant strengthening of Theorem~\ref{Tej} and hence of
the BMV theorem:

\begin{utheorem}[BMV for symmetric functions]\label{Tej2}
Fix integers $n,l \geq 1$ and a positive definite matrix $A \in
\mathbb{C}^{n \times n}$.
\begin{enumerate}
\item For each partition $\lambda$ with at most $d_A$ parts (resp.\
$d_{A^{-1}}$ parts) and positive definite $B_{n \times n}$, there exists
a threshold $m_l > 0$ such that eventual BMV-positivity holds for the
monomial symmetric function in the spectrum: for all $m_\pm \geq m_l$,
\begin{equation}\label{Eeventual}
m_\lambda(\sigma({S_{m_+,l}(A,B)})) > 0, \qquad \text{resp.} \qquad
m_\lambda(\sigma((-1)^l S_{-m_-,l}(A,B))) > 0.
\end{equation}
We adopt here and below the \textbf{convention}: we treat all spectra in
the arguments of symmetric functions as multisets with $n$ elements, so
that each eigenvalue occurs with its algebraic multiplicity.

\item If $AB=BA$ then the word ``eventual'' and bounds $d_{A^{\pm 1}}$
may be omitted: \eqref{Eeventual} holds for all $l \geq 1$, $m_+ \geq l$,
$m_- \geq 1$, and partitions $\lambda$ with at most $n$ parts.

\item Let $g({\bf x})$ be any nonzero symmetric polynomial which is an
$\mathbb{R}_+$-linear combination of $m_\lambda$ for partitions $\lambda$
with at most $d_A$ parts, resp.\ $d_{A^{-1}}$ parts. Then for any
positive definite $B$,
\[
g(\sigma(S_{m,l}(A,B))) > 0, \qquad \text{resp.,} \qquad g(\sigma((-1)^l
S_{-m,l}(A,B))) > 0
\]
for all sufficiently large $m \geq m_l(A,B,g) \gg 0$.

\item The BMV theorem also extends from $e_1$ to all Jack/Macdonald
polynomials $J_\lambda^{(\alpha)}({\bf x}), P_\lambda^{(q,t)}({\bf x})$
indexed by $\lambda$ with at most $d_{A^{\pm 1}}$ parts and with $\alpha
\in [0,\infty]$ or $q,t \in (0,1)$ -- and hence for all of their
specializations, including all Schur, Hall--Littlewood, zonal, complete
homogeneous, elementary symmetric, and monomial functions -- with at most
$d_{A^{\pm 1}}$-many parts. (In particular, the monomials include power
sums too.)
\end{enumerate}
\end{utheorem}

We show Theorem~\ref{Tej2} via an even stronger result -- see
Theorem~\ref{Tstronger}.

\begin{remark}
Theorem~\ref{Tej} for $j=1$ describes the original BMV conjecture being
``eventually true''. This is the very first and simplest special case of
Theorem~\ref{Tej2}(1) -- where $\lambda$ is the smallest partition
possible: $\lambda = (1)$.
Further within $j=1$, the $m \geq m_l$ sub-case of Theorem~\ref{Tej}
recovers the main result of Fleischhack and
Friedland~\cite{FF}.\footnote{See also~\cite[Section~1]{FF} for a
detailed account of early progress on the BMV conjecture, mostly for
small $m$.}
\end{remark}

The vigilant reader will have noticed that Theorem~\ref{Tej2} does not
address the two tightness statements in Theorem~\ref{Tej}. These are now
stated more precisely and hold for every $A > 0$, more strongly than one
might have assumed from Theorem~\ref{Tej}:

\begin{utheorem}[Tightness in BMV]\label{Ttight}
Theorem~\ref{Tej2}(1) -- for all $m_\lambda$ or even upon restricting to
elementary symmetric functions -- is tight in two ways:
\begin{enumerate}
\item For the $e_k \neq e_1$ which satisfy eventual BMV-positivity, i.e.\
for $e_k$ with $2 \leq k \leq d_A$, the threshold $m_l$ cannot be brought
down to $l$ (for non-scalar $A$), unlike the BMV conjecture. More
precisely: given any $m_0 \geq 2$ and positive definite non-scalar matrix
$A \in \mathbb{C}^{n \times n}$ (i.e.\ $d_A < n$), there exists a
positive definite $B_{n \times n}$ satisfying
\[
e_2(\sigma(S_{m,1}(A,B))), \dots, e_n(\sigma(S_{m,1}(A,B))) < 0,
\quad \forall 2 \leq m \leq m_0.
\]

\item For eventual BMV-positivity, the bound of $d_A$ is best possible.
More precisely: given any positive definite matrix $A \in \mathbb{C}^{n
\times n}$ with $d_A < n$, and any $l \geq 1$, there exists a positive
definite $B_{n \times n}$ such that we have ``eventual negativity'' above
the BMV-degree $d_A$:
\[
e_k(\sigma(S_{m,l}(A,B))) < 0, \quad \forall k \in (d_A, n],
\]
simultaneously for all $m>0$ large enough.
\end{enumerate}
\end{utheorem}

Thus in the second part, something stronger than ``negativity for
infinitely many $m$'' (i.e., the negation of ``eventual BMV positivity'')
can hold for each $k$ above the tight BMV-degree $d_A$: eventual
BMV-negativity for all $k > d_A$, and for every positive definite $A_{n
\times n}$.

\subsection{Application 5: Strengthening Uchiyama's result and answering
Ricker--V\"ath's 2002~question}\label{Suchi1}

After proving the above results, we discovered that while Lieb was
perhaps the first to ask the question about eigenvalues of BMV-words
being positive, one specific symmetrized word was studied by
Uchiyama~\cite{Uchiyama} a decade earlier. In our convention, Uchiyama
had shown:

\textit{If $0 \leq A,B \in \mathcal{B}(\mathcal{H})$, then $AB=BA$ if and
only if $0 \leq A^n B + B A^n$ for all $n \geq 1$.}

This is the $(d,N)=(0,1)$ special case of a subsequent generalization by
Ricker and V\"ath:

\begin{theorem}[\cite{RickerVath}]\label{Trv} Given operators $0 \leq A,B
\in \mathcal{B}(\mathcal{H})$, we have $AB=BA$ if and only if there exist
$N \in \mathbb{Z}_{>0}$ and $d \in [0,\infty)$ such that $A^n B + B A^n +
d A^n \geq 0$ for all $N | n$ with $n>0$. \end{theorem}

Given the above results and proofs, it is natural to try to prove the
above results under weaker hypotheses. Our next theorem is a vast
strengthening of the Ricker--V\"ath result, whence also of Uchiyama's
theorem. This requires the following notion:

\begin{defn} Given $\mathcal{H}$ and $A = A^* \in
\mathcal{B}(\mathcal{H})$, define the class $\mathcal{F}_{adm}(A)$ of
``admissible'' functions to be all continuous maps $f : [\min \sigma(A),
\max \sigma(A)] \to \mathbb{R}$, such that $|f|$ is: (1)~positive and
nondecreasing on $(0,\max \sigma(A)]$, and (2)~positive and nonincreasing
on $[\min \sigma(A), 0)$. \end{defn}

We now proceed. Note that Ricker--V\"ath showed in
\cite[Theorem~3]{RickerVath} that $AB=BA$ for \textit{matrices} $A,B \geq
0$ if and only if there exists $d \geq 0$ such that $A^\beta B + B
A^\beta + d A^\beta \geq 0$ for an increasing unbounded sequence of
\textit{real} exponents $0 < \beta = \beta_k \uparrow \infty$. Both
before and after this result, Ricker--V\"ath explicitly \textbf{ask} if
this result is true and can be showed in infinite-dimensional Hilbert
spaces. Our next result answers this affirmatively.

\begin{utheorem}\label{Tuchi}
Let $\mathcal{H}$ be a Hilbert space and $A \in \mathcal{B}(\mathcal{H})$
be self-adjoint. Fix $f \in \mathcal{F}_{adm}(A)$, an integer $m \geq 1$,
an $m$-tuple of operators ${\bf C} = (C_1, \dots, C_m) \in
\mathcal{B}(\mathcal{H})^m$ (these may not be self-adjoint, and may be
created out of $A, B$, and other operators), and a set of continuous
functions
\[
\mathcal{P} := \{ p_j, q_j, r \in C(\sigma(A), \mathbb{R}) \ : \ 1 \leq j
\leq m \}.
\]
Now define the function
\begin{equation}\label{Euchi}
W_{\mathcal{P},{\bf C},f,n}(A,-) : \mathcal{B}(\mathcal{H}) \to
\mathcal{B}(\mathcal{H}), \qquad B \mapsto A^n f(A) B + \sum_{j=1}^m
p_j(A) C_j q_j(A) + A^n r(A).
\end{equation}
If $A \geq 0$, we allow $n$ to be a positive {\em real} exponent, denoted
here by $\beta$.
\begin{enumerate}
\item If $B = B^*$ and
\[
W_{\mathcal{P},{\bf 0},f,n}(A,B) + 
W_{\mathcal{P},{\bf 0},f,n}(A,B)^* = A^n f(A) B + B A^n f(A) + 2 A^n r(A)
\]
is a positive operator for infinitely many $n \geq 1$, then $A^2 B = B
A^2$.

\item Suppose $B = B^*$ and $A - {\rm id}_{\mathcal{H}} \geq 0$. If
$W_{\mathcal{P},{\bf C},f,\beta}(A,B) + W_{\mathcal{P},{\bf
C},f,\beta}(A,B)^* \geq 0$ for an unbounded sequence $0 < \beta = \beta_k
\uparrow \infty$, then $AB=BA$.

\item Suppose $B = B^*$ and $A \geq 0$. If $W_{\mathcal{P}, {\bf 0}, f,
\beta}(A,B) + W_{\mathcal{P}, {\bf 0}, f, \beta}(A,B)^* \geq 0$ for an
unbounded sequence $0 < \beta = \beta_k \uparrow \infty$, then $AB=BA$.

\item Theorem~\ref{Trv} by Ricker--V\"ath (whence Uchiyama's result)
follows from the first part (with ${\bf C} = {\bf 0}$), even without
requiring $d \geq 0$ and $B \geq 0$. It also follows from the second
part (using $\beta_k = Nk$, say), not only if $d \geq 0$ and $B \geq 0$,
but in fact whenever $2B + d \, {\rm id}_{\mathcal{H}} \geq 0$.

Similarly, their \cite[Theorem~3]{RickerVath} follows from the third
part.
\end{enumerate}
\end{utheorem}

We make two observations.
\begin{enumerate}
\item In addition to the far larger family of words here than in
\cite{Uchiyama} or the above result from \cite{RickerVath},
Theorem~\ref{Tuchi} reveals that an infinite set of $n$ (possibly with
zero density) is already sufficient to ensure $A^2B=BA^2$ -- and $AB=BA$
if $A \geq 0$ (or more generally, using Theorem~\ref{Thilbert}(3)).

\item The method of proof is akin to that of Theorem~\ref{Thilbert}.
However, that proof additionally uses a key algebraic telescoping
identity $[A^m,B] = [A,S_{m,1}(A,B)]$ -- which is not needed here as
there are only finitely many terms in $W_{\mathcal{P},{\bf C},f,n}(A,B)$.
Thus, Theorem~\ref{Thilbert} goes through because of the above
``miracle'' of telescoping.
\end{enumerate}

In the spirit of Theorem~\ref{Tuchi}, using the family
$\mathcal{F}_{adm}(A)$ we now provide a significant strengthening of our
original ``Strong BMV'' Theorem~\ref{Thilbert}(1):

\begin{utheorem}[Sums of Hurwitz products characterize commutation]\label{Tadmissible}
Let $\mathcal{H}$ be a Hilbert space and $A \in \mathcal{B}(\mathcal{H})$
be self-adjoint.
Let a power series $f(x) = \sum_{i \geq 0} c_i x^i \in
\mathcal{F}_{adm}(A)$ be such that $\widetilde{f}(x) := \sum_{i \geq 0}
|c_i| x^i$ is also convergent at some $\delta > \|A\|$, whence on
$(-\delta,\delta)$.
\begin{enumerate}
\item For every self-adjoint $B \in \mathcal{B}(\mathcal{H})$ and $m \geq
0$, the series $W_m(A,B) := \sum_{i \geq 0} c_i S_{m+i,1}(A,B)$ converges
in norm -- as does $\sum_{i \geq 0} |c_i| S_{m+i,1}(A,B)$.
\item For such $B$, if $W_m(A,B) \geq 0$ for infinitely many $m \geq 1$,
then $A^2 B = B A^2$.
\end{enumerate}
\end{utheorem}

\begin{remark}
Our main Theorem~\ref{Thilbert}(1) is the simplest possible special case:
$f(x) = 1$. Moreover, every nontrivial odd/even power series with
nonnegative Maclaurin coefficients which is convergent on
$(-(\|A\|+\varepsilon),\|A\|+\varepsilon)$ for some $\varepsilon>0$, is
admissible.
\end{remark}

As a final remark, by spectral permanence the following holds, akin to
Theorem~\ref{Trecast}:

\begin{theorem}
The assertions in Theorems~\ref{Tuchi} and~\ref{Tadmissible} go through
verbatim for $A,B$ self-adjoint elements of any unital $C^*$-algebra.
\end{theorem}

\subsection{Motivation and results in geometry (by Vamsi P.\ Pingali)}

The Richberg technique for gluing and regularization of metrics on line
bundles has many applications in geometric analysis and several complex
variables. The analogue of such a gluing technique is hitherto
unavailable for Griffiths positively curved metrics on holomorphic vector
bundles of higher rank. Had it been available, one could have attempted
to make progress on the Griffiths conjecture on the existence of
Griffiths positively curved metrics on ample bundles. Unfortunately, the
Griffiths conjecture has been disproven recently~\cite{DuXie}.
Nonetheless, having a gluing technique is still important in view of
fully nonlinear systems of PDE like the deformed
Hermitian--Yang--Mills/Leung--Yau--Zaslow equation. We use the strong BMV
property (for matrices, and for $l=1$) to prove that an important
ingredient of pluripotential theory and Richberg-style methods fails in
the case of vector bundles.

\begin{utheorem}\label{Tvamsi}
There is no maximum function of two Hermitian metrics that preserves
Griffiths negativity.
\end{utheorem}

See Appendix~\ref{sec:geom} for the precise definitions and theorems.

\subsection{Digression on commutativity}

We end with a small digression, to view our results from the ``dual''
perspective. Theorem~\ref{Tnxn} -- which we sought to prove in order to
answer Lieb and Seiringer's question in~\cite{LS2} -- turns out to
``conversely'' yield a novel and striking characterization of
commutativity of two positive operators, which is close to the BMV
conjecture from quantum physics. Note that this condition\smallskip

\begin{center}
\textit{``if $A,B \geq 0$ and $l \geq 1$, then $AB = BA$ if and only if
$S_{m,l}(A,B) \geq 0$ for infinitely many $m \geq l$''}
\end{center}\smallskip

\noindent is novel even for $2 \times 2$ positive (semi)definite
matrices, and even for $l=2$, say. (Also for $l=1$ -- now one can even
drop the assumption of $A,B \geq 0$, if those infinite $m$ are even and
one is odd.)

From this perspective, we end with a separate, striking (and unexpected)
characterization of commutativity of positive definite matrices -- this
time involving Banach spaces:

\begin{theorem}\label{Tpolymath}
Given an integer $n \geq 1$ and two $n \times n$ positive definite matrices
$A,B$, let $\tangle{A,B}$ denote the subgroup they generate in $GL(n,
\mathbb{C})$. Then $AB = BA$ if and only if there is a group monomorphism
from $\tangle{A,B}$ into a lattice in a $1$ or $2$-dimensional real
Banach space $(\mathbb{R}^d, +, \| \cdot \|)$ for $d \in \{ 1, 2 \}$.

More generally, this characterization holds for $A,B$ invertible positive
(self-adjoint) operators over any Hilbert space $\mathcal{H}$.
\end{theorem}

That is, either $A,B$ are rational powers of each other, or words in them
form a lattice in a 2-dimensional real Banach space.

\begin{proof}
The ``if'' assertion is clear. The converse uses the
positive-definiteness hypothesis, and we first show it for matrices. If
$AB=BA$ then any element of $\tangle{A,B}$ is of the form $A^p B^q$ with
$p,q \in \mathbb{Z}$, hence positive definite as well. We claim that the
group $\tangle{A,B}$ is torsionfree. Indeed, since $A,B$ are
simultaneously unitarily diagonalizable, we may assume they are diagonal,
and now it is easy to see that $(A^p B^q)^m = {\rm Id}_n$ implies $m=0$
or $A^p B^q = {\rm Id}_n$.

But now we are done by the following result, shown by one of us in a
Polymath project~\cite{Polymath}, and applied here to $G = \tangle{a,b}$:
\textit{Two arbitrary elements $a,b$ in an arbitrary torsionfree group
$G$ commute, if and only if there is a group monomorphism from
$\tangle{a,b}$ into a real Banach space $(\mathbb{B}, +, \|\cdot\|)$.}

Note that the images of $a,b$ span a ($d=1$ or $2$-dimensional) lattice
in $\mathbb{R}^d$.

The same proof goes through for $0 \leq A,B \in
\mathcal{B}(\mathcal{H})$. Use the multiplication operator-version of the
Spectral theorem for commuting operators: there exists a common measure
space $(X, \Sigma, \mu)$ and unitary operator $U : \mathcal{H} \to
L^2(X,\mu)$, and \textit{positive} real-valued essentially bounded
measurable maps $\lambda_A, \lambda_B : X \to (0,\infty)$ such that $A =
U^* M_{\lambda_A} U$ and $B = U^* M_{\lambda_B} U$. Now proceed as above.
\end{proof}

\begin{remark}
Theorem~\ref{Tpolymath} immediately extends to the commuting of any two
psd matrices or positive operators $A',B' \geq 0$ by working with
$\tangle{A,B}$ where $A = A' + {\rm Id}, B = B' + {\rm Id}$. (This is
because $[A',B'] = [A,B]$, where $[A,B]$ is the commutator (Lie) bracket
$AB-BA$.) This result is a group-theoretic ``metric-characterization of
commutativity''; for a ring-theoretic version, we refer the reader to the
pioneering work by Le Page~\cite{LePage}. This showed that a unital
Banach algebra $\mathcal{A}_1$ is commutative if and only if it satisfies
the following property:
\begin{equation}\label{Elepage}
There\ exists\ \alpha > 0\ such\ that\ \|xy\| \leq \alpha \|yx\|\ \forall
x,y \in \mathcal{A}_1.
\end{equation}

Thus, $A'B' = B'A'$ if and only if the matrix (polynomial) algebra
$\mathcal{A}_1 := \mathbb{C}[A',B'] \subseteq \mathbb{C}^{n \times n}$
satisfies~\eqref{Elepage}. The same conclusion holds in the
$C^*$-subalgebra $\mathcal{A} := A' \mathcal{A}_1 + B' \mathcal{A}_1$,
because \textit{a $C^*$-algebra is also commutative if and only if it
satisfies~\eqref{Elepage}}. Li--Liu--Peralta~\cite{LLP} quickly deduced
this from an unpublished result of Kaplansky, which was communicated by
Kadison to Dixmier and appears in~\cite[Exercise~2.12.21]{Dix}. For
completeness: Kaplansky's result itself says that \textit{a $C^*$-algebra
is commutative if and only if it has no nonzero nilpotents}. (Also note
the formulations of commutativity in Theorem~\ref{Tmolnar} by Ogasawara,
Pedersen, Wu, and Moln\'ar.)
\end{remark}

\subsection*{Organization}

The rest of the paper is devoted to proving the above results. In the
first part: In Section~\ref{sec:Linalg} we prove Theorems~\ref{Tcommute1}
and~\ref{Tcommute2}, characterizing the Lieb--Seiringer and ``symmetrized
Lieb'' conditions -- as well as the strong BMV property -- via the
commuting of $A,B$. Then in Section~\ref{Shilbert} we show how to use
spectral measures (in the absence of eigenvectors) to prove the Hilbert
space counterpart, Theorem~\ref{Thilbert}(1),(3) (followed by
Theorem~\ref{Tfinite}). We also include a brief Appendix~\ref{Appendix}
with some basic notions (e.g.\ spectral measures, the Spectral theorem,
and the Weyl calculus), as well as some of the more routine proofs of
results in the paper.
Section~\ref{Shilbert2} proves Theorems~\ref{Thilbert}(2) and~\ref{Tnxn},
which characterize commutation via the positivity of higher order Hurwitz
products $S_{m,l}(A,B)$ (i.e., for each $l \geq 2$) using an orthonormal
eigenbasis for $A$.

The ``second part'' of this paper is devoted to applications of the above
results (beyond answering the Lieb--Seiringer questions in
Theorem~\ref{Tcommute1} above). Anderson's result on the Weyl functional
calculus is extended by Theorem~\ref{Tanderson} in Section~\ref{Sweyl}.
Next, Section~\ref{Sagniva} characterizes commuting of bounded operators
via Gibbs thermal weights (Theorem~\ref{Tthermal}); and extends
Bhowmik--Chatterjee--Putinar's multivariate completely monotone functions
with matrix arguments to bounded operators (Theorem~\ref{Tagniva}).
It also proves the Wu and weak Moln\'ar results in
Theorem~\ref{Tmolnar}, characterizing commutativity in any $C^*$-algebra.

Section~\ref{SeventualBMV2} extends the BMV conjecture/theorem to
monomial-positive symmetric functions in $\sigma(S_{m,l}(A,B))$
(Theorems~\ref{Tej2} and~\ref{Ttight}).
Next, Section~\ref{Suchi} proves our strengthening of the Uchiyama and
Ricker--V\"ath results (and answers the latter's question),
Theorem~\ref{Tuchi}, followed by the Ogasawara and Pedersen
characterizations of commutativity of a $C^*$-algebra, see
Theorem~\ref{Tmolnar}(3),(4). We also show the strengthening of our main
Theorem~\ref{Thilbert}(1): Theorem~\ref{Tadmissible}.
A final application is in Appendix~\ref{sec:geom} by Vamsi
P.~Pingali. This explains the motivations from geometry, and provides a
geometric application of Theorem~\ref{Tcommute2} to constructing
Hermitian metrics on holomorphic vector bundles.

\subsection*{Acknowledgements}

The authors thank Agniva Chatterjee, Purvi Gupta, Chandan Pradhan, Mihai
Putinar, and Aninda Sinha for useful discussions.
S.S.D.\ is partially supported by a Fulbright--Nehru Doctoral Fellowship,
being availed at MIT Cambridge. This work was (partly conducted and)
concluded while he was at MIT.
A.K.\ is partially supported by a Shanti Swarup Bhatnagar Award and an
ARG grant (ANRF/ARG/2025/011665/MS), and V.P.P.\ by a MATRICS grant
(ANRF/ARGM/2025/000313/MTR) -- from CSIR and ANRF (Govt.\ of India).
This work was also partly supported by the DST FIST program-2021
[TPN-700661].

A.K.\ and S.S.D.\ would also like to acknowledge support from the
International Centre for Theoretical Physics (ICTP, Trieste, Italy)
through the Associates Programme and from the Simons Foundation, whose
grant (Record ID: SFI-MPS-T-Institutes-00012057, AD) and ICTP's excellent
working conditions -- including the Luigi Stasi Seminar Room! --
supported this paper in May 2026.

\subsection*{AI Usage}

This paper originated with the proof of Theorem~\ref{thm:genmax}.
Previously we had shown Theorem~\ref{thm:Tposdef}, and in May~2026 we
used this and AI (specifically, ChatGPT5-Pro) to help prove
Theorem~\ref{thm:genmax}. AI also provided useful observations on our
proofs in Section~\ref{sec:Linalg}, and helpful corrections for our
typos, including in stating some results.

At this point, we discovered the connection to the BMV conjecture and the
Weyl calculus, and the paper started to expand to its present form. We
decided to \textit{not} show the paper to AI from May--November 2026. In
November 2026, we again ran the paper by AI (ChatGPT6 / Astra), which
(1)~suggested we cite Uchiyama and Moln\'ar's works (we ended up
subsuming their results into Theorems~\ref{Tmolnar} and~\ref{Tuchi}),
(2)~improved on the treatment of Theorem~\ref{thm:genmax}, and
(3)~gave us the special case of Theorem~\ref{Tthermal}: for matrices
(which we then proved for general operators which need not have
eigenvalues, and in fact a stronger result in Theorem~\ref{TC2}). Astra
also gave us the expected minor fixes, as well as valuable suggestions
that helped improve the presentation of -- and clarify somewhat ambiguous
phrases in -- several parts of the paper. It was a truly proficient
proofreader.


\section{Psd matrices commute if and only if their linearized BMV
polynomials are psd}\label{sec:Linalg}

In this section, we prove Theorems~\ref{Tcommute1} and~\ref{Tcommute2}.
We begin with a folklore lemma that helps better understand the ``weak
Lieb--Seiringer'' condition in Question~\ref{Qlieb}(1).

\begin{lemma}\label{Ldescartes}
For a complex Hermitian matrix $M_{n \times n}$, the elementary symmetric
polynomials in its eigenvalues $e_1(\sigma(M)), \dots, e_n(\sigma(M))$
are all nonnegative if and only if $M$ is positive semidefinite.
\end{lemma}

Applying Lemma~\ref{Ldescartes} to the characteristic polynomial of
$\sum_{i \geq 1} W_{i,m,l}(A,B) = S_{m,l}(A,B)$, it follows that the
Lieb--Seiringer question in Theorem~\ref{Tcommute1}(1) is equivalent to
the strong BMV property in Theorem~\ref{Tcommute1}(3).

\begin{proof}
The ``if'' part is obvious. Conversely, say $M$ has eigenvalues
$\lambda_1, \dots, \lambda_n \in \mathbb{R}$. Consider the characteristic
polynomial of $-M$, with roots $-\lambda_j$:
\[
\prod_{j=1}^n (t+\lambda_j) = \sum_{j=0}^n t^{n-j} e_j(\lambda_1, \dots,
\lambda_n) = \sum_{j=0}^n t^{n-j} e_j(M).
\]
But by Descartes' rule of signs~\cite{Descartes}, the number of positive
roots of this polynomial is at most the number of sign changes in the
coefficients, which is zero by hypothesis. Hence all $-\lambda_j \leq 0$.
\end{proof}

The next observation is that if $A,B$ commute then every word in $A,B$
-- i.e., every summand of every $S_{m,l}(A,B)$ -- can be written as
$A^{m-l} B^l$, and if $A,B$ are positive (semi)definite then so is this
word-matrix $A^{m-l} B^l = B^{l/2} A^{m-l} B^{l/2}$. But then so is the
polynomial $S_{m,l}(A,B)$, being a sum of such word-matrix pairs $W_i +
W_i^*$ (if $W_i$ is already a ``symmetric word'' then we use $(W_i +
W_i^*)/2$).

Thus, to complete the proof of Theorems~\ref{Tcommute1}
and~\ref{Tcommute2}, it remains to show the contrapositive of the
remaining cyclic implication: if $AB \neq BA$ and $A,B$ are psd, then we
will show that $S_{m,1}(A,B)$ can be psd for at most finitely many $m
\geq 1$. We begin by showing that if $B$ has rank one and no zero entry,
then $S_{m,1}(A,B)$ is never psd.

\begin{remark}
In the rest of this paper, when working with matrices we may
simultaneously unitarily conjugate $A,B$, and hence take $A$ to be a
diagonal matrix. Since $AB \neq BA$, now $B$ is necessarily not diagonal.
\end{remark}

\begin{theorem}\label{thm:semidefinite}
Choose any integer $n \geq 2$, any vector $u \in \mathbb{C}^n$ with no
zero coordinate, and any positive definite $n \times n$ diagonal matrix
$A$ that is not a scalar $c \, {\rm Id}_n$. Then $S_{m,1}(A,u u^*)$ has
a negative eigenvalue for every integer $m \geq 2$. (Hence, so does the
even larger family of matrices $S_{m,1}(P A P^*, Pu u^* P^*)$ for every
unitary matrix $P$.)
\end{theorem}

Here one has to necessarily exclude scalar matrices $A$, as each
$P A P^*$ (with $P$ unitary) is then a scalar matrix $c{\rm Id}_n$ with
$c \in (0,\infty)$, hence commutes with all $B$, and so the above
analysis applies to $S_{m,1}(A,B)$ for all $m \geq 1$.

\begin{proof}
Let $A$ have diagonal entries $\lambda_1, \dots, \lambda_n$; and $H =
(h_{ij})_{i,j=1}^n$ be a matrix that we will specify later. We directly
compute:
\[
(A^r H A^s)_{ij} = \lambda_i^r h_{ij} \lambda_j^s, \quad \forall r,s \in
\mathbb{Z}_{\geq 0}.
\]

From this it follows that
\begin{equation}\label{Ecalculation}
S_{m,1}(A,H)_{ij} = \sum_{l=0}^{m-1} \lambda_i^l h_{ij} \lambda_j^{m-1-l}
= \begin{cases} \displaystyle h_{ij} \, \frac{\lambda_i^m -
\lambda_j^m}{\lambda_i - \lambda_j}, \qquad & \text{if } \lambda_i \neq
\lambda_j,\\
h_{ij} \, m \lambda_j^{m-1}, & \text{otherwise}.
\end{cases}
\end{equation}

Thus, $S_{m,1}(A,H)$ is simply the entrywise/Schur product of $H$ with
the \textit{Loewner matrix of $f(x) = x^m$}. Note that such matrices
$L_f(\lambda_1, \dots, \lambda_n)$, whose $(i,j)$ entry is
$\frac{f(\lambda_i) - f(\lambda_j)}{\lambda_i - \lambda_j}$ if $\lambda_i
\neq \lambda_j$ and $f'(\lambda_i)$ otherwise, have prominently featured
in the Daletskii--Krein formula \cite{DaletskiiKrein} and in Loewner's
characterization of operator monotone functions \cite{Loewner34}.

We now fix $H = u u^*$ as above (and write $u = (u_1, \dots, u_n)^T$), so
that
\[
S_{m,1} (A, u u^*) = (u u^*) \circ L_m(\lambda_1, \dots, \lambda_n) = D_u
L_m(\lambda_1, \dots, \lambda_n) D_u^*,
\]
where $D_u$ is the diagonal matrix with $(j,j)$ entry $u_j$ and
$L_m(\cdot) := L_{x^m}(\cdot)$. We are to show this matrix has a negative
eigenvalue; as $D_u$ is nonsingular, by Sylvester's law of inertia
\cite{Sylvester} it suffices to show this for the Loewner matrix
$L_m(\lambda_1, \dots, \lambda_n)$, with at least two $\lambda_j > 0$
unequal.

We first address the point of repeated entries $\lambda$. If say
$\lambda_1 = \lambda_2$ then the (symmetric) matrix $L_m(\lambda_1,
\dots, \lambda_n)$ has two identical rows and columns. In such a case,
one can ``compress'' this matrix to one with distinct rows and columns,
e.g.:
\[
\begin{pmatrix}
1 & 1 & 2 & 3 & 2 \\
1 & 1 & 2 & 3 & 2 \\
2 & 2 & 5 & 7 & 5 \\
3 & 3 & 7 & 8 & 7 \\
2 & 2 & 5 & 7 & 5
\end{pmatrix}
\quad \longrightarrow \quad
\begin{pmatrix}
1 & 2 & 3 \\
2 & 5 & 7 \\
3 & 7 & 8
\end{pmatrix}
\]
corresponding to the indices partitioned by: $\{ 1, 2 \}, \{ 3, 5 \}, \{
4 \}$. And such a compression -- or its reverse ``inflation'' map -- only
changes the number of zero-eigenvalues, as e.g.\ the following result
shows (via using Sylvester's law of inertia again):

\begin{lemma}[{\cite[Lemma 3.9]{BGKP-inertia}}]
The inflation or compression -- as above -- of a real symmetric matrix
$A$ has exactly as many positive and negative eigenvalues (counted with
multiplicity) as $A$.
\end{lemma}

Thus, to show that $L_m(\lambda_1, \dots,
\lambda_n)$ has a negative eigenvalue, via this lemma we may assume that
$\lambda_1, \dots, \lambda_d$ are positive and pairwise distinct, say $0
< \lambda_1 < \cdots < \lambda_d$ with $d \geq 2$. But computing the
inertia of this matrix\footnote{The \textit{inertia} of a Hermitian
matrix is the ordered triple of numbers of (positive, zero, negative)
eigenvalues.} is precisely what the main result of
Bhatia--Friedland--Jain \cite{BFJ} does. While their result was for all
real powers $x^r, r > 0$, we recall here only the cases of use to us:
integer powers.

\begin{theorem}[{\cite[Theorem~1.1]{BFJ}}]
If $0 < p_1 < \cdots < p_d$ are reals, and $m \geq 1$ is an integer, then
\[
{\rm Inertia}(L_m(p_1, \dots, p_d)) = \begin{cases}
(\lceil m/2 \rceil, d-m, \lfloor m/2 \rfloor), \qquad & \text{if } 1 \leq
m \leq d,\\
{\rm Inertia}(L_d(p_1, \dots, p_d)), & \text{if } m > d.
\end{cases}
\]
\end{theorem}

In our setting, $m,d \geq 2$, so there are at least
$\min(\lfloor d/2 \rfloor, \lfloor m/2 \rfloor) \geq 1$
negative eigenvalues.
\end{proof}

The next case that we consider is if $B$ is nonsingular, starting with
the $2 \times 2$ situation.

\begin{theorem}\label{thm:Tposdef}
For all positive definite matrices $B = \begin{pmatrix} a & b \\
\overline{b} & c \end{pmatrix}$ with $b \neq 0$, and all positive
definite diagonal matrices $A = \begin{pmatrix} \lambda & 0 \\ 0 & \mu
\end{pmatrix}$ with $\lambda \neq \mu$, the following hold.
\begin{enumerate}
\item For all sufficiently large $m$, the matrix $S_{m,1}(A,B)$ has a
negative eigenvalue.
\item If moreover $\lambda \geq \mu \cdot 4ac/|b|^2$, then $S_{m,1}(A,B)$
has a negative eigenvalue for every $m \geq 2$.
\end{enumerate}
\end{theorem}

As a consequence, we obtain a family of positive definite counterexamples
for all $n \geq 2$:

\begin{corollary}
Let $n \geq 2$ and $B \in \mathbb{C}^{n \times n}$ be any (Hermitian,
but) non-diagonal positive definite matrix. Choose any two indices $1
\leq i < j \leq n$ such that $b := b_{ij} \neq 0$, and let the
corresponding principal submatrix be $\begin{pmatrix} a & b \\
\overline{b} & c \end{pmatrix}$. Also fix a positive definite diagonal
matrix $A$ such that $a_{ii} \neq a_{jj}$.
\begin{enumerate}
\item For all sufficiently large $m$, $S_{m,1}(A,B)$ has a negative
eigenvalue.

\item If moreover $\displaystyle \frac{\max \{ a_{ii}, a_{jj} \}}{\min \{
a_{ii}, a_{jj} \}} \geq 4 \frac{ac}{|b|^2}$, then $S_{m,1}(A,B)$ has a
negative eigenvalue for all $m \geq 2$.
\end{enumerate}
\end{corollary}

\begin{proof}
It suffices to show that $S_{m,1}(A,B)$ has a $2 \times 2$ principal
submatrix that is not psd (via considering principal minors). But this
follows immediately from Theorem~\ref{thm:Tposdef}.
\end{proof}

We next provide:

\begin{proof}[Proof of Theorem~\ref{thm:Tposdef}]
Assume without loss of generality that $\lambda > \mu > 0$.
Using~\eqref{Ecalculation},
\begin{equation}\label{Ecompute2x2}
S_{m,1}(A,B) = \begin{pmatrix} m \lambda^{m-1} a & \displaystyle
\frac{\lambda^m - \mu^m}{\lambda - \mu} b \\
\displaystyle \frac{\lambda^m - \mu^m}{\lambda - \mu} \overline{b} &
m \mu^{m-1} c
\end{pmatrix} \qquad \forall m \geq 1,
\end{equation}

\noindent and so letting $\tau := \lambda / \mu > 1$, we have
\[
\mu^{-2(m-1)} |b|^{-2} \det S_{m,1}(A,B) = \frac{ac}{|b|^2} \cdot m^2
\tau^{m-1} - \left( \frac{\tau^m - 1}{\tau- 1} \right)^2.
\]

To show $S_{m,1}(A,B)$ is not positive semidefinite, we consider its
principal minors. As the diagonal entries are positive for all $m$, and
the determinant is clearly negative for sufficiently large $m$, the first
assertion follows. To prove the second, we claim that
\[
\left( \frac{\tau^m - 1}{\tau- 1} \right)^2 \geq m^2 \tau^{m-1}
\cdot \frac{ac}{|b|^2}, \qquad \forall m \geq 2.
\]

Indeed, writing out the geometric series on the left, and multiplying
both sides by $\frac{\tau^{m-1}}{\tau^{2(m-1)}}$, this is restated as:
\[
\tau^{m-1} (1 + \tau^{-1} + \cdots + \tau^{-(m-1)})^2 \geq m^2
(ac/|b|^2).
\]

Since $\tau > 0$ and $ac/|b|^2 > 1$, it suffices in turn to show that
\[
\tau^{m-1} \geq m^2 (ac/|b|^2)^{m-1} \geq m^2 (ac/|b|^2), \qquad \forall
m \geq 2.
\]

We now invoke that $\tau \geq 4 \cdot ac/|b|^2$; hence the first
inequality easily follows from the elementary inequality $4^{m-1} \geq
m^2$ for all $m \geq 0$.
\end{proof}

Finally, we use Theorems~\ref{thm:semidefinite} and~\ref{thm:Tposdef} to
prove the remaining implication in Theorems~\ref{Tcommute1}
and~\ref{Tcommute2}. We need the following preliminary result.

\begin{lemma}\label{Lcompute}
Suppose $A,B$ are $n \times n$ matrices with $A$ a diagonal matrix. If
$I \subseteq \{ 1, \dots, n \}$ is a nonempty subset of indices, and
$A_{I \times I}$ denotes the corresponding principal submatrix of $A$,
then
\[
(AB)_{I \times I} = A_{I \times I} B_{I \times I},
\qquad (BA)_{I \times I} = B_{I \times I} A_{I \times I}.
\]
In particular, $S_{m,1}(A,B)_{I \times I} = S_{m,1}(A_{I \times I},
B_{I \times I})$ for all $m \geq 1$.
\end{lemma}

\begin{proof}
The first claim is a straightforward computation. The second is obtained
from the first by transposing. The third is obtained by repeatedly
applying the first two claims and adding up the results.
\end{proof}

With Lemma~\ref{Lcompute} at hand, we complete the proof.

\begin{proof}[Proof of Theorems~\ref{Tcommute1} and~\ref{Tcommute2}]
As discussed at the start of this section, it remains to assume $AB \neq
BA$, where $A,B$ are psd and $A$ is diagonal -- and to show that
$S_{m,1}(A,B)$ is psd for at most finitely many $m \geq 1$. We now indeed
show this -- i.e., that $S_{m,1}(A,B)$ is not psd for all $m \gg 1$ large
enough.

Since $AB \neq BA$, there exist indices $i \neq j$ such that the
$(i,j)$-entries are not equal, i.e., $(a_{ii} - a_{jj}) b_{ij} \neq 0$.
Since $B$ is psd, we have $b_{ii}, b_{jj} > 0$. Applying
Lemma~\ref{Lcompute} with $I = \{ i, j \}$, it suffices to show via
Sylvester's criterion that $S_{m,1}(A_{I \times I}, B_{I \times I})$ has
a negative eigenvalue for all $m \gg 1$ large enough. For this, write
\begin{equation}\label{E2x2}
A_{I \times I} = \begin{pmatrix} \lambda & 0 \\ 0 & \mu \end{pmatrix},
\qquad B_{I \times I} = \begin{pmatrix} a & b \\ \overline{b} & c
\end{pmatrix}.
\end{equation}

Here $a,b,c \neq 0$ and $\lambda \neq \mu$. Now we
have~\eqref{Ecompute2x2}, so if either $\lambda$ or $\mu$ is zero, it is
clear that $\det S_{m,1}(A_{I \times I},B_{I \times I}) < 0$ for all $m
\geq 2$, proving the desired claim. Henceforth, we thus assume $\lambda,
\mu > 0$. Now if $B_{I \times I}$ has rank one,
Theorem~\ref{thm:semidefinite} for $n=2$ immediately shows the claim (in
fact, for all $m \geq 2$). Similarly, if $B_{I \times I}$ is invertible
then Theorem~\ref{thm:Tposdef} finishes the proof.
\end{proof}

The vigilant reader will have noticed that the proof is not complete, as
we have not yet shown why ``infinitely many $m \geq 1$ such that
$S_{m,1}(A,B)$ is psd'' is a necessary hypothesis in
Theorem~\ref{Tcommute2}. This is addressed in a more precise form by our
next and final result in this section.

\begin{theorem}\label{Tthreshold}
Fix a dimension $n \geq 2$. Let $\mu > 0$ and $\lambda \geq e^2 \mu$, and
let $A = \begin{pmatrix} \lambda & 0 & 0 \\ 0 & \mu & 0 \\ 0 & 0 & {\rm
Id}_{n-2} \end{pmatrix}$. Then for all $m_0 \geq 2$, there exists a
positive definite matrix $B_{m_0}$ such that $A B_{m_0} \neq B_{m_0} A$
and $S_{m,1}(A,B_{m_0})$ has a negative eigenvalue if and only if $m >
m_0$. The same applies to the psd matrix $B_1 = \begin{pmatrix} 1 & 1 \\
1 & 1 \end{pmatrix} \oplus {\rm Id}_{n-2}$.
\end{theorem}

\begin{proof}
Let the scalar $\eta := \log \sqrt{\lambda/\mu} \geq 1$. We first define
the following function of $m$ -- but where $m$ is a continuous variable:
\[
\varphi(m) := \frac{e^{2m\eta} - 1}{e^{(m-1)\eta} (e^{2\eta}-1) m} =
\frac{1}{e^{2\eta}-1} \frac{e^{(m+1)\eta} - e^{-(m-1)\eta}}{m}, \qquad m
> 1.
\]

Clearly, $\varphi$ is differentiable and
\[
\varphi'(m) = \frac{1}{e^{2\eta}-1} \frac{m \eta (e^{2m\eta} + 1) -
e^{2m\eta} + 1}{e^{(m-1)\eta} m^2}.
\]

We claim that $\varphi'$ is positive on $(1,\infty)$ since $\eta \geq 1$.
Indeed, the denominators are positive, and we now lower-bound the
numerator:
\[
m \eta (e^{2m\eta} + 1) - e^{2m\eta} + 1 \geq (m-1) e^{2m\eta} + m\eta +
1 > 1.
\]

Thus, the sequence $\varphi(1), \varphi(2), \dots$ is strictly
increasing. Hence the following sequence is strictly decreasing for
integers $m \geq 1$,
\begin{align}\label{Evarphi}
\begin{aligned}
\frac{1}{\varphi(m)^2} = \frac{(e^{2\eta})^{m-1} (e^{2\eta}-1)^2
m^2}{((e^{2\eta})^m-1)^2}
= &\ \frac{m^2 (\lambda/\mu)^{m-1}
((\lambda/\mu)-1)^2}{((\lambda/\mu)^m-1)^2}\\
= &\ \frac{m \lambda^{m-1} \cdot m \mu^{m-1} \cdot
(\lambda-\mu)^2}{(\lambda^m-\mu^m)^2}.
\end{aligned}
\end{align}
Moreover, this sequence has initial term $1/\varphi(1)^2 = 1$.

We now define the desired sequence of matrices, via
\[
B_{m_0} := \begin{pmatrix}
\varphi(m_0)^2 & 1 & 0 \\ 1 & 1 & 0 \\ 0 & 0 & {\rm Id}_{n-2} \end{pmatrix}.
\]
As $\varphi(1)^2 = 1 < \varphi(2)^2 < \cdots$, the matrix $B_1$ is
positive semidefinite while $B_2, B_3, \dots$ are positive definite.
Moreover, $A B_{m_0} \neq B_{m_0} A$ for all $m_0 \geq 1$.

Now using~\eqref{Ecompute2x2}, we compute:
\[
S_{m,1}(A,B_{m_0}) = 
\begin{pmatrix} m \lambda^{m-1} \varphi(m_0)^2 &
\frac{\lambda^m - \mu^m}{\lambda - \mu} & 0 \\
\frac{\lambda^m - \mu^m}{\lambda - \mu} & m \mu^{m-1} & 0 \\
0 & 0 & m {\rm Id}_{n-2} \end{pmatrix}
= \begin{pmatrix} a_m \varphi(m_0)^2 & b_m & 0 \\ b_m & c_m & 0 \\
0 & 0 & m {\rm Id}_{n-2} \end{pmatrix},
\]
say. Since the top two diagonal entries are positive and this is a block
diagonal matrix, it has a negative eigenvalue if and only if the leading
principal $2 \times 2$ minor is negative, i.e.:
\[
a_m c_m \varphi(m_0)^2 < b_m^2.
\]

But $\displaystyle \frac{a_m c_m}{b_m^2} \varphi(m_0)^2 =
\frac{1}{\varphi(m)^2} \varphi(m_0)^2$ by~\eqref{Evarphi}, and from above
this is a strictly decreasing positive sequence that equals $1$ at
$m=m_0$. Hence $S_{m,1}(A,B_{m_0})$ has a negative eigenvalue if and only
if $m>m_0$, as claimed.
\end{proof}

\section{Strong BMV property over Hilbert spaces: General $A$}\label{Shilbert}

The goal of this section is to prove Theorem~\ref{Thilbert}(1),(3) and
related results (part~(2) will be shown in the next section). We first
note that if one assumes $A$ has a complete orthonormal eigenbasis in
$\mathcal{H}$, then one can essentially push through the matrix-theoretic
arguments used above, in order to show Theorem~\ref{Tcommute2}. However,
for a general bounded positive operator there may not exist eigenvectors,
and so one has to resort to spectral measures instead. We now explain
this proof, perhaps in somewhat greater detail than the expert reader may
require.

For the non-expert reader, we define a spectral measure $E$ on
$\mathcal{H}$ and the related measures $E_{g,h}$ for vectors $g,h \in
\mathcal{H}$ -- as well as state the Spectral theorem -- in
Appendix~\ref{Appendix}. We also recall the continuous functional
calculus there; and over here, a result from the more general
\textit{Borel functional calculus} (see \cite[Section~5.6]{Theo}, e.g.\
the proof of Theorem~5.6.3):

\begin{theorem}\label{Tborel}
Let $\mathcal{H}$ be a complex Hilbert space and $A = A^* \in
\mathcal{B}(\mathcal{H})$ an operator. Also fix a Borel subset $F \in
\mathcal{B}(\sigma(A))$ with spectral projection $E(F)$, and a bounded
Borel function $f : \sigma(A) \to \mathbb{C}$. Then,
\[
\tangle{g, f(A) E(F) h} = \int_{\sigma(A)} f(\lambda) {\bf 1}_F(\lambda)
\, d E_{g,h}(\lambda) = \int_{F} f(\lambda) \, d E_{g,h}(\lambda), \qquad
\forall g,h \in \mathcal{H}.
\]
\end{theorem}

We record some consequences that are useful below.

\begin{corollary}\label{Cborel}
$\mathcal{H}, A, F, f$ as in Theorem~\ref{Tborel}. Fix a vector $h \in
{\rm Ran}(E(F))$.
\begin{enumerate}
\item If $f \in C(\sigma(A))$ then $\|f(A)h\| \leq \|f\|_{\infty,F}
\|h\|$, where $\|f\|_{\infty,F}$ is the supnorm on $F$.

\item Suppose $g, h \in {\rm Ran}(E(F))$ and $B \in
\mathcal{B}(\mathcal{H})$. Also fix $\alpha,\alpha' \in \mathbb{Z}_{\geq
0}$ and $m \in \mathbb{Z}_{>0}$. Then,
\[
|\tangle{A^\alpha g, S_{m,1}(A,B) A^{\alpha'} h}| \leq
\|B\| \, \|g\| \, \|h\| \cdot m (\sup |F|)^{\alpha+\alpha'+m-1},
\]
where $|F| := \{ |\lambda| : \lambda \in F \}$.
\end{enumerate}
\end{corollary}

\begin{proof}
For the first part, compute using Theorem~\ref{Tborel} and that $h =
E(F)h$ and $\sigma(A) \subseteq [-\|A\|, \|A\|]$:
\begin{align*}
\| f(A) h \|^2 = &\ \tangle{f(A) h, f(A) h} = \tangle{h, f(A)^* f(A)
E(F)h} = \int_{\sigma(A)} |f(\lambda)|^2 {\bf 1}_F(\lambda) \, d
E_{h,h}(\lambda)\\
\leq &\ \int_F \|f\|_{\infty,F}^2 \, d E_{h,h}(\lambda) =
\|f\|_{\infty,F}^2 \tangle{h, E(F)h} = \|f\|_{\infty,F}^2 \|h\|^2,
\end{align*}
as desired. For the second part, compute using Cauchy--Schwarz and the
first part:
\begin{align*}
|\tangle{A^\alpha g, S_{m,1}(A,B) A^{\alpha'} h}| \leq &\ \sum_{j=0}^{m-1} |
\tangle{ A^\alpha g, A^j B A^{m-1-j+\alpha'} h}|
= \sum_{j=0}^{m-1} | \tangle{ A^{\alpha+j} g, B A^{m-1-j+\alpha'} h}|\\
\leq &\ \sum_{j=0}^{m-1} \|A^{\alpha+j}g\| \|B A^{m-1-j+\alpha'} h\|\\
\leq & \sum_{j=0}^{m-1} (\sup |F|)^{\alpha+j} \|g\| \cdot \|B\| \cdot
(\sup |F|)^{m-1-j+\alpha'} \|h\|,
\end{align*}
and we are done.
\end{proof}

\subsection{Proof of Theorem~\ref{Thilbert}}

This proof will also require understanding $S_{-m,l}(A,B)$ for $m>0$ --
specifically, which words $W_i$ occur in it.

\begin{prop}\label{Psml}
Given integers $m \geq 1$ and $l \geq 0$, and operators $A > 0$ and $B =
B^*$ in $\mathcal{B}(\mathcal{H})$,
\[
(-1)^l S_{-m,l}(A,B) = A^{-1/2} S_{m-1+l, l}(A^{-1}, B') A^{-1/2}, \qquad
\text{where} \qquad B' = A^{-1/2} B A^{-1/2}.
\]
\end{prop}

\begin{proof}
If $l=0$ then both sides equal $A^{-m}$, so we suppose henceforth that
$l \geq 1$. By \cite[Equation~(7)]{LS2}, we have
\[
\frac{1}{l!} \frac{d^l}{dt^l} (A+tB)^{-m} = (-1)^l \sum_{\substack{j_1,
\dots, j_{l+1} \geq 1,\\ j_1 + \cdots + j_{l+1} = m+l}} (A+tB)^{-j_1} B
\cdots B (A+tB)^{-j_{l+1}}.\footnote{This was stated without proof
in~\cite{LS2}; here is a sketch for the curious reader. Differentiate
both sides of $(A+tB) (A+tB)^{-1} = {\rm id}_\mathcal{H}$ to solve for
$\frac{d}{dt} (A+tB)^{-1}$; induct on $m \geq 1$ to obtain $\frac{d}{dt}
(A+tB)^{-m}$; then induct on $l \geq 1$.}
\]

Now let $B'$ as in the statement, and let $j'_i := j_i-1$ for all $i$.
Since $A+tB$ remains invertible for small $|t| \geq 0$, letting $t \to 0$
yields:
\[
(-1)^l S_{-m,l}(A,B) = \sum_{\substack{j'_1, \dots, j'_{l+1} \geq 0,\\
j'_1 + \cdots + j'_{l+1} = m-1}} A^{-1/2} \cdot A^{-j'_1} B' A^{-j'_2} B'
\cdots B' A^{-j'_{l+1}} \cdot A^{-1/2},
\]
and this is precisely $ A^{-1/2} S_{m-1+l, l}(A^{-1}, B') A^{-1/2}$.
\end{proof}

Following these preliminaries, we come to:

\begin{proof}[Proof of Theorem~\ref{Thilbert}(1),(3)]
We begin by dispensing with part~(3). Suppose $\sigma(A) \cap \sigma(-A)
\subseteq \{ 0 \}$. Consider the map $f : \sigma(A) \to [0,\infty)$,
sending $\lambda \mapsto \lambda^2$. This is clearly continuous, has
domain the compact set $\sigma(A)$, and is a bijection onto $\sigma(A^2)$
by Spectral mapping (see Theorem~\ref{Tfnalcalc}(4)). Hence $f$ is a
homeomorphism, and so $f^{-1} : \sigma(A^2) \twoheadrightarrow \sigma(A)$
is also a continuous bijection. Thus, it can be approximated uniformly on
$\sigma(A^2)$ by a sequence of polynomials $p_k(x)$. Since $A^2$ commutes
with $B$, so does each $p_k(A^2)$. Taking limits, $B$ commutes with
\[
\lim_{k \to \infty} \int_{\sigma(A^2)} p_k(\lambda') \,
dE_{A^2}(\lambda') = \int_{\sigma(A^2)} f^{-1}(\lambda') \,
dE_{A^2}(\lambda') = \int_{\sigma(A)} \lambda \, dE_A(\lambda) = A,
\]
where the last equality uses ``change of variables'' or push-forwards of
measures $dE_{vw}(\cdot)$.\medskip

We now come to Theorem~\ref{Thilbert}(1). Note that the hypotheses must
hold for (i)~infinitely many $m = m_k>0$ or (ii)~infinitely many $m =
m_k<0$. Suppose we can show that $[A^2,B]=0$ in case~(i). Then it holds
in case~(ii) as follows: if $-S_{-m,1}(A,B) \geq 0$ for infinitely many
$m>0$, then by Proposition~\ref{Psml} and case~(i), $S_{m,1}(A^{-1},B')
\geq 0$ infinitely often. Hence $A^{-2} B' = B' A^{-2}$, and multiplying
on the left and right by $A^{5/2}$ yields: $BA^2=A^2B$.

It remains to show $[A^2,B]=0$ in case~(i). We prove the contrapositive,
splitting the argument into steps for ease of exposition.\medskip

\noindent \textbf{Step 1:}
Henceforth, we use without further reference the standard notation:
\[
|A| := \sqrt{A^2} = \sqrt{A^*A} \geq 0.
\]
Note that the spectral projections of $|A|$ and $A$ are related via:
$E_{|A|}(F) = E_A(F \cup -F)$, for all $F \subseteq [0,\infty)$.
Moreover, if $B$ commutes with $|A|$ then $[A^2,B]=0$.

To show Theorem~\ref{Thilbert}(1) in case~(i), we suppose for
contradiction that $[A^2,B] \neq 0$ -- whence $[|A|,B] \neq 0$. We will
use here (and below) the following technical lemma.

\begin{lemma}\label{Lstasi}
Suppose $\mathcal{H}$ is a Hilbert space and $A' \geq 0, B =
B^*$ are in $\mathcal{B}(\mathcal{H})$ with
\[
\sigma(A') \, \subseteq \, [\lambda_0, \|A'\|] \, \subseteq \, [0,
\|A'\|] \qquad \text{for some } \lambda_0 \geq 0.
\]
For a Borel subset $F \subset \mathbb{R}$, define $F^{\sigma'} := F \cap
\sigma(A')$. If $A'B \neq BA'$, then there exist scalars
\begin{equation}\label{Eclaim2}
\lambda_0 \leq c < a \leq b \leq \|A'\| \quad \text{such that} \quad
E(F^{\sigma'}) \cdot B \cdot E([\lambda_0, c]^{\sigma'}) \neq 0, \ \
\frac{b+c}{2} < a,
\end{equation}
where $F$ is a (closed/open) interval with endpoints $a, b$.
\end{lemma}

We apply this with $\lambda_0 = 0$ below, with a unique exception while
proving Theorem~\ref{Tuchi}(2).

\begin{remark}
If e.g.\ $A'$ was a Hermitian matrix, then $A'B \neq BA'$ would mean $B$
does not preserve some $\mu$-eigenspace $V_\mu$ of $A'$, meaning that
some off-diagonal entry in $P_{V_\lambda} B V_\mu$ is nonzero. In the
context of Lemma~\ref{Lstasi}, $\mu \in [\lambda_0,c]$ and $\lambda \in
[a, b]$.
\end{remark}

\begin{proof}[Proof of Lemma~\ref{Lstasi}]
Say $[A',B] \neq 0$. We first claim:
\[
B \cdot E([\lambda_0,c]^{\sigma'}) \neq E([\lambda_0,c]^{\sigma'}) \cdot
B \qquad \text{ for some } c \leq \|A'\|,
\]
which would in turn imply $c < \|A'\|$. This is a standard ``good sets''
argument: consider the collection of Borel subsets of the spectrum,
\[
\Sigma := \{ F \in \mathcal{B}(\sigma(A')): B \cdot E(F) = E(F) \cdot B
\}.
\]

Now if the above claim is false, then $\Sigma$ is a $\sigma$-algebra
containing all ``intervals'' $[\lambda_0,c]^{\sigma'}$, and since the
family $\{ [\lambda_0,c]^{\sigma'} : c \geq 0 \}$ generates
$\mathcal{B}(\sigma(A'))$, $\Sigma = \mathcal{B}(\sigma(A'))$. But then
$A'B = BA'$ by the Spectral theorem~\ref{Tspectral}, which is
false.\medskip

Thus, let $F = [\lambda_0,c]$, and set $P := E(F^{\sigma'})$. Then $BP
\neq PB$, and we claim that $P^\perp B P \neq 0$, where $P^\perp = {\rm
id}_{\mathcal{H}} - P = E((F^c)^{\sigma'}) = E((c, \|A'\|]^{\sigma'})$.
Indeed, if $P^\perp B P = 0$, then ${\rm Ran}(P)$ is invariant for $B$,
hence reducing for $B = B^*$, and so $PB = BP$, which is
false.\footnote{In basic terms, if $P^\perp B P = 0$, adding $PBP$
yields: $BP = PBP$. Now taking adjoints yields $PB = (BP)^* = (PBP)^* =
PBP = BP$, which is false.} Hence,
\[
E((c, \|A'\|]^{\sigma'}) \cdot B \cdot E([\lambda_0,c]^{\sigma'}) \neq 0.
\]

This implies there exists an integer $k \gg 0$ such that
\[
c + \frac{1}{k} < \|A'\| \quad \text{and} \quad E([c +
\textstyle\frac{1}{k}, \|A'\|]^{\sigma'}) \cdot B \cdot
E([\lambda_0,c]^{\sigma'}) \neq 0.
\]

Partition $[c+\frac{1}{k}, \|A'\|]$ into intervals of length $<
\frac{1}{k}$. At least one of these, say $F$ with endpoints $a < b$,
satisfies
\[
E(F^{\sigma'}) \cdot B \cdot E([\lambda_0,c]^{\sigma'}) \neq 0.
\]
Moreover, $b + c  < (a + \frac{1}{k}) + c \leq a + a$, and the proof is
complete.
\end{proof}

\noindent \textbf{Step 2:}
For this proof, given a subset $F \subseteq [-\|A\|, \|A\|]$, define
$F^\sigma := F \cap \sigma(A)$. Applying Lemma~\ref{Lstasi} with $A' =
|A|$ and $\lambda_0 = 0$, and noting that
\[
E_{A'}(F^{\sigma'}) = E_A(F^\sigma) + E_A((-F)^\sigma) \qquad \text{and}
\qquad E_{A'}([0,c]^{\sigma'}) = E_A([-c,c]^\sigma),
\]
there exists $\varepsilon \in \{ \pm 1 \}$ such that
\[
E_A((\varepsilon F)^\sigma) \cdot B \cdot E_A([-c,c]^\sigma) \neq 0.
\]

Hence there exists $v_\mu \in {\rm Ran}(E_A([-c,c]^\sigma))$ such that
$v := E_A((\varepsilon F)^\sigma) B v_\mu$ is nonzero.\footnote{This
choice of notation $v_\mu, v$ is in order to be consistent with the next
section.} We will obtain a contradiction to this fact, by showing that
$\| v \|^2 \leq 0$.\medskip

\noindent \textbf{Step 3:}
We use the algebraic identity $[A^m, B] = [A, S_{m,1}(A,B)]$, valid in
any associative ring. Taking inner products,
\begin{align*}
\tangle{v, A^m B v_\mu} = &\ \tangle{v, B A^m v_\mu} + \tangle{v, [A^m,
B] v_\mu} = \tangle{v, B A^m v_\mu} + \tangle{v, [A, S_{m,1}(A,B)]
v_\mu}\\
= &\ \tangle{v, B A^m v_\mu} + \tangle{Av, S_{m,1}(A,B) v_\mu} -
\tangle{v, S_{m,1}(A,B) A v_\mu}.
\end{align*}

Therefore,
\begin{equation}\label{Ecommutator}
|\tangle{v, A^m B v_\mu}| \leq |\tangle{v, B A^m v_\mu}| + |\tangle{Av,
S_{m,1}(A,B) v_\mu}| + |\tangle{v, S_{m,1}(A,B) A v_\mu}|.
\end{equation}

The first term on the right in~\eqref{Ecommutator} is bounded via
Cauchy--Schwarz and Corollary~\ref{Cborel}(1):
\begin{equation}\label{Ebound2}
|\tangle{v, B A^m v_\mu}| \leq \|v\| \cdot \|B\| \|A^m v_\mu\|
\leq \|v\| \, \|B\| \, c^m \, \|v_\mu\|.
\end{equation}

We next lower-bound the left-hand term in~\eqref{Ecommutator}:
\begin{align*}
\tangle{v, A^m B v_\mu} = &\ \tangle{ E_A((\varepsilon F)^\sigma)^2 v,
A^m B v_\mu} = \tangle{v, E_A((\varepsilon F)^\sigma)^2 A^m B v_\mu}\\
= &\ \tangle{ v, A^m E_A((\varepsilon F)^\sigma)^2 B v_\mu}
= \tangle{ v, A^m E_A((\varepsilon F)^\sigma) v}
= \int_{(\varepsilon F)^\sigma} \lambda^m \, d E_{v,v}(\lambda),
\end{align*}
where the final equality is by Theorem~\ref{Tborel}. Thus,
\begin{equation}\label{Ebound1}
|\tangle{v, A^m B v_\mu}| = \left| \int_{(\varepsilon F)^\sigma}
\lambda^m \, d E_{v,v}(\lambda) \right| \geq a^m \int_{\sigma(A)} {\bf
1}_{(\varepsilon F)^\sigma}(\lambda) \, d E_{v,v}(\lambda) = a^m
\tangle{v, E((\varepsilon F)^\sigma) v} = a^m \|v\|^2.
\end{equation}

\noindent \textbf{Step 4:}
The rest of the right-hand side of~\eqref{Ecommutator} is where we use
the hypothesis in Theorem~\ref{Thilbert}(1) that $S_{m_k,1}(A,B)$ is
positive for all $k$. In this case $(g,h) \mapsto \tangle{g,
S_{m_k,1}(A,B)h}$ is a semi-inner product on $\mathcal{H}$, and so
Cauchy--Schwarz holds for it. Using this and suitably applying
Corollary~\ref{Cborel}(2) (for $F \leadsto (\varepsilon F)^\sigma,
[-c,c]^\sigma$), we obtain:
\begin{align*}
|\tangle{Av, S_{m_k,1}(A,B)v_\mu}| \leq &\
\tangle{Av, S_{m_k,1}(A,B) Av}^{1/2} \tangle{v_\mu, S_{m_k,1}(A,B)
v_\mu}^{1/2}\\
\leq &\ (\|B\| \, \|v\|^2 m_k b^{m_k+1})^{1/2} 
(\|B\| \, \|v_\mu\|^2 m_k c^{m_k-1})^{1/2}.
\end{align*}

A similar calculation upper-bounds the final term:
\[
|\tangle{v, S_{m_k,1}(A,B)Av_\mu}| \leq 
(\|B\| \, \|v\|^2 m_k b^{m_k-1})^{1/2} 
(\|B\| \, \|v_\mu\|^2 m_k c^{m_k+1})^{1/2}.
\]

Combining these bounds yields:
\begin{equation}\label{Ebound3}
| \tangle{v, [A, S_{m_k,1}(A,B)] v_\mu} | \leq \|B\| \, \|v\| \, \| v_\mu
\| \cdot m_k (bc)^{(m_k-1)/2} (b+c).
\end{equation}

\noindent \textbf{Step 5:}
We now finish the proof. Suppose $S_{m_k,1}(A,B)$ is positive for all
$k$. Combining the bounds in~\eqref{Ebound2}, \eqref{Ebound1},
and~\eqref{Ebound3} with the inequality~\eqref{Ecommutator},
\begin{align*}
a^{m_k} \|v\|^2 \leq c^{m_k} \|B\| \, \|v\| \, \|v_\mu\| + \|B\| \,
\|v\| \, \| v_\mu \| \cdot m_k (bc)^{(m_k-1)/2} (b+c).
\end{align*}

But $a > (b+c)/2 > \sqrt{bc} \geq c \geq 0$ by Lemma~\ref{Lstasi}. Hence
dividing by $a^{m_k}$ and taking $k \to \infty$,
\[
\|v\|^2 \leq \|B\| \, \|v\| \, \| v_\mu \| \lim_{k \to \infty} \left(
(c/a)^{m_k} + \frac{b+c}{a} m_k (bc / a^2)^{(m_k-1)/2} \right) = 0.
\]
This produces the desired contradiction to the final paragraph in
Step~2.
\end{proof}

\subsection{Finitely many Hurwitz products don't work}

We next settle the question of \textit{minimality}.
Namely, one can combine the above results to obtain:

\begin{theorem}\label{Tcommute4}
The following are equivalent for positive (semi)definite matrices $A,B$
-- or more generally for $0 \leq A,B \in \mathcal{B}(\mathcal{H})$ with
$A$ invertible:
\begin{enumerate}
\item For all polynomials $p(A,B)$ that are $\mathbb{R}_+$-linear
combinations of monomials/words in $A,B \geq 0$, the
symmetrized matrix $p(A,B) + p(A,B)^*$ is psd.

\item For infinitely many integers $m$, the linear BMV-coefficient
${\rm sgn}(m) \cdot S_{m,1}(A,B)$ is psd.

\item For any fixed $l \geq 1$ and infinitely many integers $m$, the
BMV-coefficient ${\rm sgn}(m)^l \cdot S_{m,l}(A,B)$ is psd. (If $l \geq
2$ then we assume the hypotheses of Theorem~\ref{Tnxn}.)
\item $AB=BA$.
\end{enumerate}
\end{theorem}

\noindent (After all, clearly $(4) \implies (1) \implies (2),(3)$, and we
show above and below respectively that $(2),(3)$ each imply $(4)$.) In
particular, any family of symmetrized test polynomials that is in between
the uncountable family~(1) and the countable sets in~(2) or~(3), suffices
to guarantee $AB=BA$.

Thus, it is natural to ask if the (un)countability can be replaced by
finiteness. The next proof resolves this,
negatively.\footnote{Without the symmetrization hypothesis, indeed this
uncountable family can be replaced by just one word: $AB$, since $AB \geq
0 \implies AB = (AB)^* = BA$.}

\begin{proof}[Proof of Theorem~\ref{Tfinite}]
We will display two $2 \times 2$ real psd matrices $A,B$ with this
property; the result for a general (real or complex) Hilbert space is
then obtained by letting $\mathcal{H}_0$ be a 2-dimensional subspace --
with a fixed orthonormal basis -- on which $A,B$ act; and extending $A,B$
to each act on all of $\mathcal{H}$ by fixing $\mathcal{H}_0^\perp$
pointwise.

Fix positive real scalars $\lambda \neq \mu$, and define
\[
A = \begin{pmatrix} \lambda & 0 \\ 0 & \mu \end{pmatrix}, \qquad
B_\varepsilon := \begin{pmatrix} 1 & \varepsilon \\ \varepsilon & 1
\end{pmatrix}, \ \varepsilon \in [0,1).
\]

Then $A, B_\varepsilon$ are positive definite, and as each word
$W(A,B_\varepsilon)$ in each $\sgn(m)^l \cdot S_{m,l}(A,B_\varepsilon)$ has
finite length, one checks via Proposition~\ref{Psml} that
\[
\lim_{\varepsilon \to 0^+} W(A,B_\varepsilon) + W(A,B_\varepsilon)^* =
W(A, {\rm Id}_2) + W(A, {\rm Id}_2)^* = 2 A^{m-l}, \qquad \forall l \in
\mathbb{Z}_{\geq 0}, m \in \mathbb{Z},
\]
and this is positive definite. Hence there exists a small positive
$\varepsilon_0 = \varepsilon(\{ m_k, l_k : 1 \leq k \leq n \})$ below
which all word-symmetrizations $W(A, B_\varepsilon) + W(A,
B_\varepsilon)^*$ are positive definite within each Hurwitz product
$\sgn(m_k)^l \cdot S_{m_k, l_k}(A, B_\varepsilon)$ (for all $0 <
\varepsilon < \varepsilon_0$). But since $\varepsilon$ and $\lambda \neq
\mu$ are positive, $A^2 B_\varepsilon \neq B_\varepsilon A^2$.
\end{proof}

\subsection{Going from commutation to Hurwitz products}

Here we prove a ``converse'' result to Theorem~\ref{Thilbert}. That is,
the following ``goes the other way'', and provides sufficient conditions
to deduce the positivity of Hurwitz products from commutation.

\begin{theorem}\label{Tnotpos3}
Fix a Hilbert space $\mathcal{H}$ and self-adjoint $A,B \in
\mathcal{B}(\mathcal{H})$.
\begin{enumerate}
\item Suppose $A^2B = BA^2$, and one of the following also holds:
\begin{enumerate}
\item $S_{m_\pm,1}(A,B) \geq 0$ for some odd $m_- > 0$ and even $m_+ >
0$.
\item $\pm B \geq 0$, and $S_{m_+,1}(A,B) \geq 0$ for some even $m_+ >
0$.
\end{enumerate}
Then $B$ commutes with $A$.

\item Say $AB = BA$, $l \geq 1$, and $S_{m_\pm,l}(A,B) \geq 0$ for one
even $m_+ \geq l$ and one odd $m_- \geq l$. Then $S_{m,l}(A,B) \geq 0$
for all $m > l$.
\end{enumerate}
\end{theorem}

\begin{corollary}
Theorem~\ref{Tnotpos3} holds verbatim in an arbitrary $C^*$-algebra.
\end{corollary}

The proof of Theorem~\ref{Tnotpos3} becomes easier if one works with
operator-matrices. Thus, we first introduce this notation; the technology
in the next remark is used not just here, but in the next section to
prove Theorems~\ref{Tnxn} (for $l=1$) and~\ref{Thilbert}(2) as well.

\begin{remark}[Conformally partitioned operator-matrices]\label{Rconformally}
Suppose $W_i$ $(i=1,2,3)$ are closed subspaces of $\mathcal{H}$ with
respective projection operators $P_i$, satisfying
\begin{equation}
W_1 \oplus W_2 \oplus W_3 = \mathcal{H}, \qquad \text{so} \qquad
P_1 + P_2 + P_3 = {\rm id}_{\mathcal{H}}.
\end{equation}

Then given an operator $T \in \mathcal{B}(\mathcal{H})$, one can write it
as a conformally partitioned matrix (with operator entries):
\begin{equation}
T = \sum_{i=1}^3 P_i \cdot T \cdot \sum_{j=1}^3 P_j \longleftrightarrow
(T_{ij})_{i,j=1}^3, \qquad T_{ij} = (P_{W_i} \circ T)|_{W_j} = P_{W_i}
\circ T \circ P_{W_j}.
\end{equation}
This immediately shows that $(T_{ij})^* = (T^*)_{ji}$ for $1 \leq i,j
\leq 3$. It also shows that one can compose operators like matrix
multiplication:
\begin{equation}\label{Emult}
(T \circ B)_{ij} = P_i T \, {\rm id}_{\mathcal{H}} B P_j = P_i T (P_1 +
P_2 + P_3) B P_j = \sum_{p=1}^3 (P_i T P_p)(P_p B P_j) = \sum_{p=1}^3
T_{ip} B_{pj}.
\end{equation}
\end{remark}

\begin{proof}[Proof of Theorem~\ref{Tnotpos3}]
We begin with some preparation. Define the spectral projections
\begin{equation}\label{Espectralproj}
P_1 := E((0,\|A\|]), \qquad
P_2 := E([-\|A\|,0)), \qquad
P_3 := E(\{0\}),
\end{equation}
and denote their images by $W_1, W_2 = \mathcal{H}_-$ (as in
Theorem~\ref{Tnxn}(3)), and $W_3 = \ker(A)$, respectively. Then
Remark~\ref{Rconformally} applies here.
Writing $A$ as a $3 \times 3$ block matrix, and since each $P_i$ commutes
with $A$ by the Spectral theorem, we obtain $P_i A P_j = A P_i P_j = 0$
for $i \neq j$, and so $A$ is block diagonal. Write $|A|$ and $B$ in the
same way, say
\begin{equation}\label{EAB3x3}
A = \begin{pmatrix} A_+ & 0 & 0 \\ 0 & A_- & 0 \\ 0 & 0 & 0
\end{pmatrix}, \qquad
|A| = \begin{pmatrix} A_+ & 0 & 0 \\ 0 & -A_- & 0 \\ 0 & 0 & 0
\end{pmatrix}, \qquad
B = \begin{pmatrix} a & b & c \\ b^* & d & e \\ c^* & e^* & f
\end{pmatrix}.
\end{equation}

Note that $\pm A_\pm$ are positive (hence self-adjoint) and injective
operators. Now if $B$ commutes with $A^2$, then $[B,|A|] = 0$ as in the
proof of Corollary~\ref{Ccstar}. Equating $|A|B = B|A|$ and comparing
entries yields:
\[
aA_+ = A_+a, \quad dA_- = A_-d, \quad A_+ b = -b A_-, \quad A_+ c = 0,
\quad A_- e = 0.
\]
But $A_\pm$ are injective, so $c,e$ are zero, whence $B$ is block
diagonal on $(W_1 \oplus W_2) \oplus W_3$.

\begin{remark}\label{Rdense}
This remark will be used repeatedly in the sequel. The operators $A_\pm$
are not necessarily invertible, e.g.\ if $0$ is a limit point in
$\sigma(A)$. However, they are injective and have dense ranges in $W_1,
W_2$, respectively. Hence the same holds for $A_+^\alpha$ for $\alpha >
0$. If now $A_+^\alpha X A_+^\alpha \geq 0$ for some $X = X^*$, then $X
\geq 0$ by this denseness.
\end{remark}

To prove $AB=BA$, it suffices to show $b=0$, so we may work with the
leading principal $2 \times 2$ submatrices of $A,B$. To that end, let
$\mathcal{H}' := W_1 \oplus W_2$, and for any $3 \times 3$
block-partitioned operator-matrix $M \in \mathcal{B}(\mathcal{H})$ as
above, denote by $M'$ its leading principal $2 \times 2$
operator-submatrix:
\begin{equation}
M' := (P_1 + P_2) M (P_1 + P_2).
\end{equation}

This implies for instance that if $M \geq 0$ then $M' \geq 0$.
Moreover, since $A,B$ are block diagonal from above,
\[
0 = 0' = [|A|,B]' = [|A'|,B'], 
\]
and $[A,B] = 0$ if and only if $[A',B'] = 0$.

Write $A = |A| \cdot U = U \cdot |A|$, where $U := {\rm id}_{W_1} \oplus
-{\rm id}_{W_2} \oplus {\bf 0}_{W_3}$. Then $U' = \begin{pmatrix} {\rm
id}_{W_1} & 0 \\ 0 & - {\rm id}_{W_2} \end{pmatrix}$, a self-adjoint and
unitary operator.
\begin{itemize}
\item[(1)(a)] Since $B$ commutes with $|A|$, we have $B'$ commutes with
$|A'|$, and so
\begin{align*}
0 = 0' \leq S_{m_\pm,1}(A,B)' = &\ S_{m_\pm,1}(U' |A'|, B') =
|A'|^{m_\pm-1} S_{m_\pm,1}(U',B')\\
= &\ |A'|^{\frac{m_\pm-1}{2}} S_{m_\pm,1}(U',B')
|A'|^{\frac{m_\pm-1}{2}}.
\end{align*}

Suitably applying Remark~\ref{Rdense} yields: $S_{m_\pm,1}(U',B') \geq
0$. As $m_+$ is even and $m_-$ is odd, an explicit computation reveals:
\begin{equation}\label{ESmplus}
S_{m_+,1}(U',B') = \begin{pmatrix} m_+a & 0 \\ 0 & -m_+d \end{pmatrix},
\qquad
S_{m_-,1}(U',B') = \begin{pmatrix} m_-a & b \\ b^* & m_-d \end{pmatrix}.
\end{equation}

Since these are both positive, pre- and post- composing with $P_2$ shows
$\pm d \geq 0$, so $d=0$. Now suppose for contradiction that $b \neq 0$.
Then there exists $v_\mu \in W_2$ such that $v := b v_\mu \in W_1
\setminus \{ 0 \}$. Compute for any $\varepsilon > 0$, using from above
that $d=0$:
\begin{equation}\label{Evarv}
0 \leq \varepsilon^{-1} \tangle{ \varepsilon v \oplus -v_\mu,
S_{m_-,1}(U',B') (\varepsilon v \oplus -v_\mu)} = \varepsilon m_-
\tangle{v, a v} - 2 \|v\|^2.
\end{equation}
But this is false, since the final expression is negative for small $0 <
\varepsilon \ll 1$. Hence $b=0$.

\item[(1)(b)]
First say $B \geq 0$. Choose $m_+ > 0$ even such that $S_{m_+,1}(A,B)
\geq 0$. By the same computations for $m_+$ as above (ending
in~\eqref{ESmplus}),
\[
0 \leq S_{m_+,1}(U',B') = \begin{pmatrix} m_+a & 0 \\ 0 & -m_+d
\end{pmatrix},
\]
whence $-d \geq 0$. But $d = P_2 B P_2 \geq 0$ as well, so $d = 0$. Again
working as in~\eqref{Evarv} (but with $S_{m_-,1}(U',B')$ replaced by $B'
\geq 0$), it follows that $b=0$.

If instead $-B \geq 0$, then $0 \leq S_{m_+,1}(A,B) = S_{m_+,1}(-A,-B)$,
and $-B \geq 0$ commutes with $(-A)^2$. Hence the previous working gives
$0 = [-A,-B] = [A,B]$.

\item[(2)] Recall~\eqref{Espectralproj}, \eqref{EAB3x3}, and the notation
of Remark~\ref{Rconformally}, decomposing $\mathcal{H} = W_1 \oplus W_2
\oplus W_3$. Now $AB=BA$, so $B$ commutes with $P_1, P_2, P_3$ and hence
is block-diagonal too:
\[
A = \begin{pmatrix} A_+ & 0 & 0 \\ 0 & A_- & 0 \\ 0 & 0 & 0
\end{pmatrix}, \qquad
B = \begin{pmatrix} B_+ & 0 & 0 \\ 0 & B_- & 0 \\ 0 & 0 & B_0
\end{pmatrix}.
\]

Transversally to $\sigma(A)$, decompose $\mathcal{H} = (\ker B) \oplus
(\ker B)^\perp$. Now define
\[
W_j^0 := W_j \cap \ker B, \qquad W_j^\perp := W_j \cap (\ker B)^\perp,
\qquad j = 1,2,3.
\]
Then with respect to the decomposition into six closed subspaces
\[
\mathcal{H} =
(W_1^0 \oplus W_1^\perp) \oplus
(W_2^0 \oplus W_2^\perp) \oplus
(W_3^0 \oplus W_3^\perp),
\]
we write
\begin{equation}\label{E6block}
A = {\rm diag}(A_+^0, A_+^\perp, A_-^0, A_-^\perp, 0, 0), \qquad
B = {\rm diag}(0, B_+^\perp, 0, B_-^\perp, 0, B_0^\perp).
\end{equation}
Note that $\pm A_\pm^0, \pm A_\pm^\perp$ are positive, and $B_\pm^\perp,
B_0^\perp$ are self-adjoint.

Now for each $1 \leq p \leq 6$, the $p$th block diagonal components of
$A$ and $B$ commute, since $AB=BA$. Thus, for $m \geq l$ we have:
\[
\binom{m}{l}^{-1} S_{m,l}(A,B) = {\rm diag} (0, (A_+^\perp)^{m-l}
(B_+^\perp)^l, 0, (A_-^\perp)^{m-l} (B_-^\perp)^l, 0, \delta_{m,l}
(B_0^\perp)^l).
\]

Since this matrix is positive for $m = m_\pm \geq l$ (even and odd), the
four operators
\begin{equation}\label{Eabpmperp}
(A_\pm^\perp)^{m_+-l} (B_\pm^\perp)^l, \quad
(A_\pm^\perp)^{m_--l} (B_\pm^\perp)^l
\end{equation}
are positive.
Now it suffices to prove that $(A_\pm^\perp)^{m-l} (B_\pm^\perp)^l \geq
0$ for all $m > l$. Indeed, by the Spectral theorem it suffices to show
that for every joint spectral point $(x,y) \in
\sigma_{joint}(A_+^\perp,B_+^\perp)$, the scalar $x^{m-l} y^l \geq 0$ for
all $m > l$ (the proof is similar for $A_-^\perp, B_-^\perp$). Suppose
this is not true for some $x,y,m$. Then $x^{m-l}y^l < 0$, so $x \neq 0$.
If $m$ has the same parity as $m_-$, then we get
\[
0 > (x^2)^{(m_--m)/2} \cdot (x^{m-l}y^l) = x^{m_--l} y^l \geq 0,
\]
and this contradiction proves the claim. Similarly if $m,m_+$ have the
same parity.\qedhere
\end{itemize}
\end{proof}

\subsection{Tightness of the hypotheses for $l=1$}\label{S1tight}

We conclude this section by exploring results related to the $l=1$ cases
of Theorems~\ref{Thilbert} and~\ref{Tnotpos3}. A first -- and immediate
-- consequence is:

\begin{corollary}\label{Cnotpos}
Fix a Hilbert space $\mathcal{H}$ and self-adjoint $A,B \in
\mathcal{B}(\mathcal{H})$.
\begin{enumerate}
\item If $S_{m,1}(A,B) \geq 0$ for infinitely many $m \geq 1$, then
$A^2B=BA^2$.

\item If not all of these $m$ are of the same parity, then $AB=BA$.

\item Similarly, if not all of these $m$ are odd, and $B \geq 0$, then
$AB=BA$.
\end{enumerate}
\end{corollary}

Thus, it is natural to ask if the assumptions in these additional
conclusions are tight. We show here that they are:

\begin{prop}
Given any Hilbert space $\mathcal{H}$ of dimension in $[2,\infty]$,
there exist self-adjoint $A,B_\pm \in \mathcal{B}(\mathcal{H})$ such that
\begin{enumerate}
\item $[A, B_\pm] \neq 0$;
\item $B_-$ is positive but $A,B_+$ are not;
\item $S_{m,1}(A,B_+) \geq 0$ for all even $m>0$; and
\item $S_{m,1}(A,B_-) \geq 0$ for all odd $m>0$.
\end{enumerate}
\end{prop}

\begin{proof}
By ``padding by zeros'' or ``extending by identity'' to an
ortho-complement, it suffices to work with $\dim \mathcal{H} = 2$. Now
define
\[
A := \begin{pmatrix} 1 & 0 \\ 0 & -1 \end{pmatrix},
\qquad
B_- := \begin{pmatrix} 1 & 1 \\ 1 & 1 \end{pmatrix},
\qquad
B_+ := \begin{pmatrix} 1 & 1 \\ 1 & -1 \end{pmatrix}.
\]

Then it is easily verified that properties~(1), (2) hold. Moreover,
\[
S_{m,1}(A,B_+) = \begin{pmatrix} m & (1 - (-1)^m)/2 \\ (1 - (-1)^m)/2 &
m(-1)^m \end{pmatrix},
\]
and for $m$ even this is $m {\rm Id}_2 \geq 0$. Similarly,
\[
S_{m,1}(A,B_-) = \begin{pmatrix} m & (1 - (-1)^m)/2 \\ (1 - (-1)^m)/2 &
m(-1)^{m-1} \end{pmatrix},
\]
and for $m$ odd this too is psd.
\end{proof}

\section{Strong BMV in Hilbert spaces via higher order
terms}\label{Shilbert2}

Having shown the above results for $l=1$ first over matrices and then
over Hilbert spaces, we next turn to the higher $l$ results:
Theorems~\ref{Thilbert}(2) and~\ref{Tnxn}. While writing this section, we
first showed Theorem~\ref{Thilbert}(2) for matrices -- even before
working out the Hilbert space results in Section~\ref{Shilbert}. It turns
out, even this result for matrices is already significantly more
challenging than Theorems~\ref{Tcommute1} and~\ref{Tcommute2}, for two
reasons:
\begin{enumerate}
\item Firstly, even for $2 \times 2$ matrices, the $l=1$ case above was
worked out using explicit computations. This approach works for $l=2$ as
well, as we discovered. However, it quickly becomes much harder to
compute $S_{m,l}(A,B)$ explicitly for larger $l$ -- even for $2 \times 2$
matrices. This is because computing higher powers of $B$ is more involved
when $B$ is not necessarily diagonal (we are not aware of a closed-form
expression).

The workaround is to not consider exact formulas for the entries of
$S_{m,l}(A,B)$, but only their ``fastest growing terms'' / asymptotics in
$m$. This requires tools that may be termed ``analytic combinatorics''
and are akin to techniques used in~\cite{CK,FS} (as mentioned in the
introduction).

Moreover, if $l \geq2$ then the asymptotics for the $(2,2)$-entry
in~\eqref{Ecompute2x2} differ from the $l=1$ case; see
Section~\ref{Sasymptotics}.

\item There is a further challenge -- that of going from $2 \times 2$ to
$n \times n$. Note that the $l=1$ case for $2 \times 2$ matrices helped
show the same result for $n \times n$ matrices $A,B$ for any $n$, using
Lemma~\ref{Lcompute}. However, the last line of this lemma is no longer
true if we replace $S_{m,1}(A,B)$ with $S_{m,l}(A,B)$ for any $l \geq 2$.

The workaround for this challenge was to not work with the entries of
$S_{m,l}(A,B)$ -- as we did for $l=1$ -- but with its \textit{blocks} (as
in Section~\ref{Shilbert}, which we worked out after
Theorem~\ref{Thilbert}(2) for matrices). We did this first to show
Theorems~\ref{Thilbert}(2) and~\ref{Tnxn} for $l=2$; and then wrote up a
similar but more intricate proof for all $l \geq 3$. Both cases used
analytic combinatorics.
\end{enumerate}

We now come to the proof of Theorems~\ref{Thilbert}(2) and~\ref{Tnxn} in
their full generality: over Hilbert spaces, and for all $l \geq 2$. (The
$l=1$ case was proved in Section~\ref{Shilbert}.) Our proof for $n \times
n$ matrices can be tweaked to work over Hilbert spaces (and the reader
interested in matrices not operators, can read this section verbatim,
assuming $\mathcal{H} = \mathbb{C}^n$). The key is to use ``blocks'' of
$B$, i.e.\ the technology of Remark~\ref{Rconformally}.

\subsection{Proof of Theorem~\ref{Tnxn}: two implications}\label{S2impl}

We now begin the proof of Theorem~\ref{Tnxn}. Clearly, $(1) \implies
(2)$. Here we show that $(3) \implies (1)$, starting with the equivalence
in~(3). Indeed, if $B P_{\mathcal{H}_-} = 0$ then $P_{\mathcal{H}_-} B
P_{\mathcal{H}_-} = 0$. Conversely, if $P_{\mathcal{H}_-} B
P_{\mathcal{H}_-} = 0$ then we use $\sqrt{B}$ (since $B \geq 0$) to
compute:
\[
\| \sqrt{B} P_{\mathcal{H}_-} v \|^2 = \tangle{v, P_{\mathcal{H}_-} B
P_{\mathcal{H}_-} v} = 0, \qquad \forall v \in \mathcal{H}.
\]
This implies $\sqrt{B} P_{\mathcal{H}_-} = 0$ and so $B P_{\mathcal{H}_-}
= 0$.

Next, we prove $(3) \implies (1)$ if $A>0$. Indeed, then $B$ commutes
with $A$, hence with $\sqrt{A}$ as in the proof of
Corollary~\ref{Ccstar}. Hence $S_{m,l}(A,B) = \binom{m}{l} A^{(m-l)/2}
B^l A^{(m-l)/2} \geq 0$ for $m \geq 0$. In turn, this implies (since
$A^{-1}, B' = A^{-1/2} B A^{-1/2}$ commute) that $(-1)^l S_{-m,l}(A,B)
\geq 0$ for all $m>0$, by Proposition~\ref{Psml}.

The ``true'' proof of $(3) \implies (1)$ is when $A$ is not positive
definite, but we do have $B P_{\mathcal{H}_-} = 0$. We are to show
$S_{m,l}(A,B) \geq 0$ for $m \geq l \geq 1$.
For this, adopt the notation of Remark~\ref{Rconformally}, for the
subspaces and projections at the start of the proof of
Theorem~\ref{Tnotpos3}. Writing $A, |A|, B$ as in~\eqref{EAB3x3}, the
equation $B P_{\mathcal{H}_-} = 0 \leq B$ implies $b,d,e$ are all zero,
and so
\[
A = \begin{pmatrix} A_+ & 0 & 0 \\ 0 & A_- & 0 \\ 0 & 0 & 0 \end{pmatrix},
\qquad
B = \begin{pmatrix} a & 0 & c \\ 0 & 0 & 0 \\ c^* & 0 & f \end{pmatrix}.
\]
It also follows that $[A_+,a]=0$. Moreover, as $B$ preserves $\ker(A)$,
for all $v \in \ker(A) = W_3$ we have
$c v = P_1 B P_3 v = P_1 Bv = 0$,
since $Bv \in {\rm Ran}(P_3)$ and $P_1 P_3 = 0$. Therefore $c = 0$, so
\[
S_{m,l}(A,B) = \binom{m}{l} {\rm diag}(A_+^{m-l} a^l,0, \delta_{m,l}f^l)
\geq 0.
\]
This concludes the proof. \qed

\subsection{Proof of Theorem~\ref{Tnxn}$(2) \implies (3)$ and
Theorem~\ref{Thilbert}(2) -- for negative $m$ or for $l = 1$}

We first assume the $m \geq 0$ case for both of these implications --
which are all that remain to be shown -- and prove their counterparts
when $(-1)^l S_{-m,l}(A,B) \geq 0$ for infinitely many $m > 0$.
The arguments are both similar: from the hypotheses and
Proposition~\ref{Psml}, we deduce by pre- and post- multiplying $(-1)^l
S_{-m,l}(A,B)$ by $A^{1/2}$, that $S_{m+l-1,l}(A^{-1}, B') \geq 0$, with
$B' = A^{-1/2} B A^{-1/2} \geq 0$ if $B \geq 0$. Hence the $m \geq 0$
case of either result gives us the commutation of $B$ with $A^2$ or with
$A$. Finally, if $A>0$ then $\mathcal{H}_- = 0$, and so $B
P_{\mathcal{H}_-} = 0$.

We next prove the $l = 1, m \geq 0$ case of Theorem~\ref{Tnxn}
$(2)\implies (3)$. Given~(3), Theorem~\ref{Thilbert} implies $A^2B=BA^2$,
and Theorem~\ref{Tnotpos3}(1)(b) in turn gives $AB=BA$. Moreover, the
argument in the proof of Corollary~\ref{Ccstar} gives $[|A|,B]=0$ as
well. Now write $A, |A|, B$ as in~\eqref{EAB3x3}, and
use~\eqref{Ecompute2x2} to compute for infinitely many even $m>1$:
\[
P_2 S_{m,1}(A,B) P_2 = S_{m,1}(A_-,d) = S_{m,1}(-A_-,-d) \geq 0.
\]

As seen above, $-A_- = P_2 |A| P_2$ is positive, so $0 = [-d,-A_-] =
[d,A_-]$ by Theorem~\ref{Thilbert}(1),(3). Thus, the above equation
becomes:
\[
(-A_-)^{(m-1)/2} \cdot (-m d) \cdot (-A_-)^{(m-1)/2} \geq 0.
\]
By suitably applying Remark~\ref{Rdense}, we get $-md \geq 0 \leq P_2 B
P_2 = d$. Thus $0 = d = P_2 B P_2$. \qed

\subsection{Proof of Theorem~\ref{Tnxn}$(2) \implies (3)$ and
Theorem~\ref{Thilbert}(2) -- for $m \geq 0, l = 2$}\label{Squad}

Having dispensed with the non-principal cases, we now begin the two
proofs -- via common notation and arguments. Note that the conclusions we
desire for $l \geq 2$ (with $B \geq 0$ given) are: $A^2B=BA^2$ and
$AB=BA$ -- which implies $A^2B=BA^2$.

Thus, we work via contradiction for both proofs -- for now allowing any
$l \geq 2$, and starting with $A^2B \neq BA^2$. This means the positive
operator $A^2$ cannot be a scalar, so it has at least two eigenvalues. We
now divide the proof into steps, for ease of exposition. The first two steps are for all $l \geq 2$, and then we specialize to $l=2$ in this
subsection.\medskip

\noindent \textbf{Step 1:}
\textit{(Common preparation for all $l \geq 2$.)}

Since we have the positivity $S_{m_k,l}(A,B) \geq 0$ for an infinite
increasing sequence
\[
{\bf m} := (m_1 < m_2 < \cdots) \subseteq [l, \infty) \cap \mathbb{Z},
\]
we can assume all $m_k$ have the same parity, and the following limit
exists:
\begin{equation}\label{Ebse}
\bse_l({\bf m}) := \lim_{k \to \infty} (-1)^{m_k-l}, \qquad l \geq 1.
\end{equation}

Since $\sigma(A) \setminus \{ 0 \}$ is discrete, so is $\sigma(A^2)
\setminus \{ 0 \} = \sigma(A)^2 \setminus \{ 0 \}$ (by Spectral mapping). 
Arrange these eigenvalues in decreasing order:
\[
\sigma(A^2) \setminus \{ 0 \} = \{ \lambda_1 > \lambda_2 > \cdots \}
\subset (0, \|A\|^2],
\]

Next, define for each real scalar $\lambda$ the (eigen)spaces
\begin{equation}\label{EAsquare}
V'_\lambda := \ker(A^2 - \lambda \, {\rm id}_{\mathcal{H}}), \qquad
V_\lambda := \ker(A - \lambda \, {\rm id}_{\mathcal{H}}).
\end{equation}
Then $V'_0 = V_0$,
\[
V'_{\lambda} = V_{\sqrt{\lambda}} \oplus V_{-\sqrt{\lambda}} \ \forall
\lambda>0, \qquad \text{and} \qquad \mathcal{H} = V'_0 \oplus
\bigoplus_{n \geq 1} V'_{\lambda_n}.
\]

Following these preliminaries, we proceed. Suppose $B(V'_{\lambda_1})
\subseteq V'_{\lambda_1}$. Then $B = B^*$ also preserves
$(V'_{\lambda_1})^\perp$, as does $A$. Moreover, on $V'_{\lambda_1}$, $B$
commutes with $A^2|_{V'_{\lambda_1}} \equiv \lambda_1 \, {\rm
id}_{\mathcal{H}}$, and $S_{m,l}(A,B)$ splits into a block-diagonal
operator-matrix, on $V'_{\lambda_1} \oplus (V'_{\lambda_1})^\perp$. Thus,
we may reduce to working in $(V'_{\lambda_1})^\perp$.

Continuing similarly, if the above holds for all $\lambda_n > 0$ then we
reduce to working on $V'_0$, on which $A^2 \equiv 0$ and so $[A^2,B] =
0$. Thus, there exists $\lambda'_n > 0$ such that $B$ does not preserve
$V'_{\lambda_n}$. We take the smallest such $n$, act by $A,B$ on
$(V'_{\lambda_1} \oplus \cdots \oplus V'_{\lambda_{n-1}})^\perp$, and
without loss of generality let $n=1$.

This yields $0 \leq \eta < \lambda_1$ in $\sigma(A^2)$ such that
$(P_{V'_\eta} \circ B)(V'_{\lambda_1}) \neq 0$. Let
\begin{equation}
\lambda := \sqrt{\lambda_1}.
\end{equation}

Recalling that $B = B^*$, the nonvanishing of $P_{V'_\eta} \circ B \circ
P_{V'_{\lambda_1}}$ and hence $P_{V'_{\lambda_1}} \circ B \circ
P_{V'_\eta}$ yields the following ``$2 \times 2$ submatrix'' being
nonzero (see~\eqref{EAsquare}):
\begin{equation}\label{E2x2block}
\begin{pmatrix}
P_{V_{\sqrt{\eta}} } \circ B \circ P_{V_\lambda} &
P_{V_{\sqrt{\eta}} } \circ B \circ P_{V_{-\lambda}} \\
P_{V_{-\sqrt{\eta}} } \circ B \circ P_{V_\lambda} &
P_{V_{-\sqrt{\eta}} } \circ B \circ P_{V_{-\lambda}}
\end{pmatrix}
= P_{V'_{\eta}} \circ B \circ P_{V'_\lambda} \neq 0.
\end{equation}
(If $\eta = 0$ then we just consider the first row.) Hence at least one
row is nonzero, and we set that eigenvalue of $A$ to be $\mu \in \{ \pm
\sqrt{\eta}  \}$. In particular, $|\mu| < \lambda$.

We now define four closed subspaces:
\begin{align}\label{Eabove1}
\begin{aligned}
W_1 := &\ V_{\lambda}, \qquad
W_2 := V_{-\lambda}, \qquad
W_3 := V_{\mu},\\
W_4 := &\ (W_1 \oplus W_2 \oplus W_3)^\perp = \bigoplus_{\nu \in
\sigma(A) \setminus \{ \mu, \pm \lambda \}} V_\nu,
\end{aligned}
\end{align}
and $P_j := P_{W_j}$ for all $j$. Then $\mathcal{H} = \oplus_{j=1}^4
W_j$, and akin to Remark~\ref{Rconformally}, we will write any $T \in
\mathcal{B}(\mathcal{H})$ as a conformally partitioned $4 \times 4$
operator-matrix. In particular, write
\begin{equation}\label{Eabove2}
A = {\rm diag}( \lambda P_1,\ -\lambda P_2,\ \mu P_3,\ D := P_4 A P_4 ),
\qquad B = (B_{ij} := P_i B P_j)_{i,j=1}^4
\end{equation}
and note from above that
\[
|\mu|, \|D\| \leq \sqrt{\lambda_2} < \lambda = \sqrt{\lambda_1}, \qquad
(B_{13}, B_{23}) \neq (0,0),
\]
the latter since $[A^2,B] \neq 0$.\medskip

\noindent \textbf{Step 2:} \textit{(Three technical lemmas.)}

We now adopt the notation of Remark~\ref{Rconformally} and proceed. We
require several technical lemmas later on, so we present them here and
then present the remaining proof of Theorem~\ref{Tnxn} ``uninterrupted''.
The first lemma is folklore even for Hilbert spaces; it involves the
Aitken factorization~\cite[\S III.29]{Aitken} of a block $2 \times 2$
``matrix'' using the Schur complement~\cite{Schur}:

\begin{lemma}\label{Laitken}
Suppose $\mathcal{H}_1, \mathcal{H}_2$ are complex Hilbert spaces, and
$\mathcal{H} = \mathcal{H}_1 \oplus \mathcal{H}_2$. Write $P_j :=
P_{\mathcal{H}_j}$ for $j=1,2$. Now suppose $C \in
\mathcal{B}(\mathcal{H})$ is written as a block $2 \times 2$
operator-matrix $\begin{pmatrix} Q & R \\ S & T \end{pmatrix}$, with $Q =
P_1 C P_1$ invertible. Then we have the factorization:
\[
\begin{pmatrix} P_1 & 0 \\ -S Q^{-1} & P_2 \end{pmatrix}
\begin{pmatrix} Q & R \\ S & T \end{pmatrix}
\begin{pmatrix} P_1 & -Q^{-1} R \\ 0 & P_2 \end{pmatrix}
= \begin{pmatrix} Q & 0 \\ 0 & T - S Q^{-1} R \end{pmatrix}.
\]

In particular, if $C = C^*$ (so $S = R^*$) and if $Q$ is moreover
positive on $\mathcal{H}_1$, then $C$ is positive on $\mathcal{H}$ if and
only if the Schur complement $T - R^* Q^{-1} R$ is so.
\end{lemma}

This is a direct computation. The next lemma evaluates the limit of a
family of positive operators.

\begin{lemma}\label{Lspectralcompute}
Suppose $\mathcal{H}$ is a Hilbert space and $B \in
\mathcal{B}(\mathcal{H})$ is self-adjoint and positive. Fix positive real
scalars
\[
\alpha, \beta, a, b, c > 0, \qquad \text{with} \qquad \alpha \geq \beta,
\ b > ac,
\]
and define the family of functions
\[
f_\varepsilon(B) := c B^{\alpha-\beta} - b B^\alpha (a B^\beta +
\varepsilon \, {\rm id}_{\mathcal{H}})^{-1}, \qquad \varepsilon > 0.
\]

Also fix a vector $v \in \mathcal{H}$ such that if $\alpha=\beta$ then
$E_B(\{0\}) v = 0$, where $E_B(\{0\})$ is the spectral projection.
Then $\tangle{v, f_\varepsilon(B)v} \to \frac{ac-b}{a} \tangle{v,
B^{\alpha-\beta}v}$ as $\varepsilon \to 0^+$.
\end{lemma}

\begin{proof}
If $v=0$, both sides of the limit equal zero from the start. So we assume
henceforth that $v \neq 0$. Now define the function
\[
f_\varepsilon(\lambda) := c\lambda^{\alpha-\beta} -
\frac{b\lambda^\alpha}{a\lambda^\beta + \varepsilon}, \qquad \lambda \in
[0, \|B\|], \ \varepsilon>0.
\]
(Here if $\lambda = 0$ and $\alpha = \beta$, we define $0^0 := 1$.)
It is easily checked for $\lambda,\varepsilon$ in these ranges that
\begin{equation}\label{EDCT}
\frac{ca-b}{a} \|B\|^{\alpha-\beta} \leq \frac{ca-b}{a}
\lambda^{\alpha-\beta} \leq f_\varepsilon(\lambda) \leq c
\lambda^{\alpha-\beta} \leq c \|B\|^{\alpha-\beta}.
\end{equation}

We now proceed. By the Spectral theorem~\ref{Tspectral},
\[
\lim_{\varepsilon \to 0^+} \tangle{v, f_\varepsilon(B)v} =
\lim_{\varepsilon \to 0^+} \int_{\sigma(B)} f_\varepsilon(\lambda) \, d
E_{v,v}(\lambda) =
\int_{\sigma(B)} \lim_{\varepsilon \to 0^+} f_\varepsilon(\lambda) \, d
E_{v,v}(\lambda),
\]
where the second equality follows from the Dominated Convergence Theorem,
using~\eqref{EDCT}.

We consider two cases. First suppose $\alpha>\beta$. Then for both
$\lambda=0$ and $\lambda>0$, one checks: $\lim_{\varepsilon \to 0^+}
f_\varepsilon(\lambda) = \frac{ac-b}{a} \lambda^{\alpha-\beta}$. Hence
the Spectral theorem yields:
\[
\lim_{\varepsilon \to 0^+} \tangle{v, f_\varepsilon(B)v} =
\int_{\sigma(B)} \frac{ac-b}{a} \lambda^{\alpha-\beta} \, d E_{v,v}(\lambda)
= \frac{ac-b}{a} \tangle{v, B^{\alpha-\beta}v}.
\]

The other case is $\alpha=\beta$. Now
$\lim_{\varepsilon \to 0^+} f_\varepsilon(\lambda)$ equals
$\frac{ac-b}{a} \, \lambda^{\alpha-\beta}$ if $\lambda>0$,
and $c$ if $\lambda=0$. Hence we are again done by the Spectral theorem:
\begin{align*}
\lim_{\varepsilon \to 0^+} \tangle{v, f_\varepsilon(B)v} = &\
\int_{\sigma(B)} \frac{ac-b}{a} \lambda^{\alpha-\beta} \, d
E_{v,v}(\lambda) + \frac{b}{a} \int_{\{0\}} d E_{v,v}(\lambda)\\
= &\ \frac{ac-b}{a} \tangle{v, B^{\alpha-\beta}v} + \frac{b}{a}
\tangle{v,E_B(\{0\})v}. \qedhere
\end{align*}
\end{proof}

Our final technical lemma involves a block-partitioned positive
operator-matrix:

\begin{lemma}\label{LB12}
Suppose $B = B^* \in \mathcal{B}(\mathcal{H})$. Write $\mathcal{H} = W_1
\oplus W_2$ with $W_1, W_2$ closed subspaces, and let $B =
(B_{ij})_{i,j=1}^2$ conformally. If $B$ is a positive operator, then
$E_{B_{11}}(\{0\}) B_{12} = 0$.
\end{lemma}

Informally, if $B$ was a matrix and hence could be diagonalized via a
unitary change-of-basis, then $E_{B_{11}}(\{0\})$ represents the
projection onto the $0$-eigenspace. Now Lemma~\ref{LB12} says that in the
new basis, the corresponding rows of $B_{12}$ (and hence of any vector $v
= B_{12} v_\mu$) are zero.

\begin{proof}
Let $P_j$ denote the projection onto $W_j$ for $j=1,2$, and abbreviate $E
:= E_{B_{11}}(\{0\})$. Then the positivity of $B$ implies that of
\[
\begin{pmatrix} E & 0 \\ 0 & P_2 \end{pmatrix} \cdot B \cdot
\begin{pmatrix} E & 0 \\ 0 & P_2 \end{pmatrix}
= \begin{pmatrix} E B_{11} E & E B_{12} \\ B_{21} E & B_{22} \end{pmatrix}.
\]

We claim that $E B_{11} E = 0$. Indeed, since $E,B_{11}$ commute by the
Spectral theorem, if any $w \in \mathcal{H}$ is written as $w_1 \oplus
w_2$, then $E B_{11} E w = E (B_{11} P_1)w = E B_{11} w_1$. Now write
$w_1 = w_0 \oplus w_\perp$, where $w_0 \in \ker(B_{11})$ and $w_\perp
\perp \ker(B_{11})$. Then,
\[
E B_{11} E w = E B_{11} w_0 + B_{11} E w_\perp = 0 + 0,
\]
proving the claim.

Finally, suppose for contradiction that $E B_{12} \neq 0$. Then there
exists $v_\mu \in W_2$ such that $v := E B_{12} v_\mu \neq 0$. Now
compute for $\varepsilon > 0$ (using that $E = E^2$):
\[
0 \leq \varepsilon^{-1} \left\langle \begin{pmatrix} v \\ -\varepsilon
v_\mu \end{pmatrix}, \begin{pmatrix} 0 & E B_{12} \\ B_{21} E & B_{22}
\end{pmatrix} \begin{pmatrix} v \\ -\varepsilon v_\mu \end{pmatrix}
\right\rangle = -2 \|v\|^2 + \varepsilon \tangle{v_\mu, B_{22} v_\mu}.
\]
But this quantity is negative for small positive values of $\varepsilon$,
which gives the desired contradiction.
\end{proof}

\noindent \textbf{Step 3:}
Now we begin the proofs for $l=2$. These repeatedly use two computations
-- the first is the two-variable summation
\begin{equation}\label{Eapgp}
\sum_{\substack{\alpha, \beta, \gamma \in \mathbb{Z}_{\geq 0},\\
\alpha+\beta+\gamma = m-2}} \mu^{\alpha+\gamma} \lambda^\beta =
\sum_{\beta=0}^{m-2} (m-1-\beta) \lambda^\beta \mu^{m-2-\beta} =
\frac{\lambda^m - m \mu^{m-1}\lambda + (m-1) \mu^m}{(\lambda-\mu)^2}.
\end{equation}

The second is the analogous three-variable summation:
\begin{equation}\label{Echs}
\sum_{\substack{\alpha, \beta, \gamma \in \mathbb{Z}_{\geq 0},\\
\alpha+\beta+\gamma = m-2}} x^\alpha y^\beta z^\gamma =
\frac{x^m}{(x-y)(x-z)} + \frac{y^m}{(y-x)(y-z)} + \frac{z^m}{(z-x)(z-y)},
\end{equation}
for any $x,y,z$ in a commutative ring in which $x-y, y-z, z-x$ are
invertible. This latter follows from the theory of symmetric
functions~\cite{Macdonald} -- specifically, the equivalence of 
Littlewood's ``SSYT-definition'' of the \textit{complete homogeneous
symmetric function}
\[
h_{m-2}(x,y,z) := \sum_{\substack{\alpha, \beta, \gamma \in
\mathbb{Z}_{\geq 0},\\ \alpha+\beta+\gamma = m-2}} x^\alpha y^\beta
z^\gamma
\]
with Cauchy's bialternant definition of $h_{m-2}$ as a Schur polynomial
for $m \geq 2$, aka the Weyl character formula in type $A$ representation
theory:
\begin{align*}
h_{m-2}(x,y,z) = &\ s_{(m-2,0,0)}(x,y,z) = \frac{\det {\bf
x}^{(m,1,0)}}{\det {\bf x}^{(2,1,0)}}\\
= &\ \frac{1}{(x-y)(y-z)(x-z)} \det \begin{pmatrix} x^m & y^m & z^m \\ x
& y & z \\ 1 & 1 & 1 \end{pmatrix}.
\end{align*}

\noindent \textbf{Step 4:}
Returning to the setting at the end of Step~1, we explicitly compute the
Hurwitz product
\[
S_{m,2}(A,B) = \sum_{\substack{\alpha, \beta, \gamma \in \mathbb{Z}_{\geq
0},\\ \alpha+\beta+\gamma = m-2}} A^\alpha B A^\beta B A^\gamma
\]
block by block. We only compute the $(i,j)$ blocks for $i,j \leq 3$ since
in the sequel, we are concerned with the leading principal $3 \times 3$
submatrix.

Set $(\nu_1, \nu_2, \nu_3) := (\lambda,-\lambda, \mu)$. Then for $1 \leq
i,j \leq 3$, we have
\begin{align}\label{Easymptotics}
&\ S_{m,2}(A,B)_{ij}\\
= &\ \sum_{\substack{\alpha, \beta, \gamma \in
\mathbb{Z}_{\geq 0},\\ \alpha+\beta+\gamma = m-2}}
(B_{i1} \nu_i^\alpha \nu_1^\beta \nu_j^\gamma B_{1j} +
 B_{i2} \nu_i^\alpha \nu_2^\beta \nu_j^\gamma B_{2j} +
 B_{i3} \nu_i^\alpha \nu_3^\beta \nu_j^\gamma B_{3j} +
 B_{i4} \nu_i^\alpha D^\beta \nu_j^\gamma B_{4j}).\notag
\end{align}

Unfortunately (for everyone!), we now explicitly evaluate this for the
six cases $1 \leq i \leq j \leq 3$ (since the $i>j$ cases follow from the
self-adjointness of $S_{m,2}(A,B)$) -- and follow it up by computing the
asymptotics too.
In what follows, we use~\eqref{Eapgp} and~\eqref{Echs} (as stated before
Step~1) without mention -- and in the asymptotics, we will assume $m$
belongs to an infinite sequence of \textbf{all even or all odd $m_k$ --
see the notation~\eqref{Ebse}.}
\begin{enumerate}
\item $(i,j)=(1,1)$:
The first summand is $B_{11} \lambda^{m-2} B_{11}$, which adds up to
$\binom{m}{2} \lambda^{m-2} B_{11}^2$.

The second summand is $\lambda^{\alpha+\gamma} (-\lambda)^\beta B_{12}
B_{21}$, which adds up to
\[
\frac{(-1)^m + 2m-1}{4} \lambda^{m-2} B_{12} B_{21}.
\]

The third summand is $\lambda^{\alpha+\gamma} \mu^\beta B_{13} B_{31}$,
which adds up to
\[
\frac{\mu^m - m \lambda^{m-1} \mu + (m-1) \lambda^m}{(\lambda-\mu)^2}
B_{13} B_{31}.
\]

The fourth summand is $B_{14} \cdot \lambda^{\alpha+\gamma} D^\beta \cdot
B_{41}$. This case is more interesting: since $\|D\| \leq
\sqrt{\lambda_2} < \lambda$, we have $\lambda \, {\rm id}_{\mathcal{H}} -
D$ is invertible. Now~\eqref{Eapgp} applies in the commutative subring
$\mathbb{C}[[D]] \subseteq \mathcal{B}(W_4)$ given by convergent power
series in $D$. Hence the fourth summand adds up to
\[
B_{14} \cdot (\lambda \, {\rm id}_{\mathcal{H}} - D)^{-2} (D^m - m
\lambda^{m-1} D + (m-1)\lambda^m {\rm id}_{W_4}) \cdot B_{41}.
\]

Now for the asymptotics. The fastest growing quantity comes from the
first summand, and is $m^2 \lambda^{m-2}$. Since $\mu$ and $\|D\|$ are at
most $\sqrt{\lambda_2} < \lambda_1$, we conclude that
\begin{equation}\label{Egrowth11l2}
\lim_{m \to \infty} \frac{1}{m^2 \lambda^{m-2}} S_{m,2}(A,B)_{11} =
\frac{1}{2} B_{11}^2.
\end{equation}

\item $(i,j)=(1,2)$: The four summands add up, respectively, to:
\begin{align*}
&\ \frac{(-1)^m + 2m-1}{4} \lambda^{m-2} B_{11} B_{12}, \qquad
\frac{1 + (-1)^m(2m-1)}{4} \lambda^{m-2} B_{12} B_{22},\\
&\ \left( \frac{\lambda^m}{2\lambda(\lambda-\mu)} +
\frac{(-\lambda)^m}{2\lambda(\lambda+\mu)} +
\frac{\mu^m}{(\mu-\lambda)(\mu+\lambda)} \right) B_{13} B_{32},\\
&\ B_{14} \left( \frac{\lambda^m}{2\lambda} (\lambda - D)^{-1} +
\frac{(-\lambda)^m}{2\lambda} (\lambda + D)^{-1} + D^m (D - \lambda)^{-1}
(D + \lambda)^{-1} \right) B_{42},
\end{align*}
where $D+c$ stands for $D + c \, {\rm id}_{W_4}$. Since $\|D\| \leq
\sqrt{\lambda_2} < \lambda$, we get the asymptotics
\begin{equation}
\lim_{m \to \infty} \frac{1}{m^2 \lambda^{m-2}} S_{m,2}(A,B)_{12} = 0.
\end{equation}

\item $(i,j)=(1,3)$:
The first three summands add up, respectively, to:
\begin{align*}
&\ \frac{\mu^m - m \lambda^{m-1} \mu + (m-1) \lambda^m}{(\lambda -
\mu)^2} B_{11} B_{13}, \\
&\ \left( \frac{\lambda^m}{2\lambda(\lambda-\mu)} +
\frac{(-\lambda)^m}{2\lambda(\lambda+\mu)} +
\frac{\mu^m}{(\mu-\lambda)(\mu+\lambda)} \right) B_{12} B_{23},\\
&\ \frac{\lambda^m - m \mu^{m-1} \lambda + (m-1) \mu^m}{(\lambda -
\mu)^2} B_{13} B_{33}.
\end{align*}

We also upper-bound the fourth summand's sum in norm, using that
$\max(\|D\|, \mu) \leq \sqrt{\lambda_2} < \lambda$:
\begin{align*}
0 \leq &\ \left\| \sum_{\substack{\alpha, \beta, \gamma \in
\mathbb{Z}_{\geq 0},\\ \alpha+\beta+\gamma = m-2}} B_{14} \cdot
\lambda^\alpha D^\beta \mu^\gamma \cdot B_{43} \right\| \leq \|B_{14}\|
\, \|B_{43}\| \sum_{\substack{\alpha, \beta, \gamma \in \mathbb{Z}_{\geq
0},\\
\alpha+\beta+\gamma = m-2}} \lambda^\alpha \| \mu^\gamma D^\beta \|\\
&\ \leq \|B_{14}\| \, \|B_{43}\| \sum_{\substack{\alpha, \beta, \gamma
\in \mathbb{Z}_{\geq 0},\\ \alpha+\beta+\gamma = m-2}} \lambda^\alpha
(\sqrt{\lambda_2})^{\beta+\gamma}\\
= &\ \frac{\lambda^m - m (\sqrt{\lambda_2})^{m-1} \lambda + (m-1)
(\sqrt{\lambda_2})^m}{(\lambda - \sqrt{\lambda_2})^2} \|B_{14}\| \,
\|B_{43}\|.
\end{align*}

Now verify that this upper bound, as well as the second and third sums
prior to it, when divided by $m \lambda^{m-2}$, go to zero. This yields:
\begin{equation}
\lim_{m \to \infty} \frac{1}{m \lambda^{m-2}} S_{m,2}(A,B)_{13} =
\frac{\lambda}{\lambda-\mu} B_{11} B_{13}.
\end{equation}

\item $(i,j)=(2,2)$:
The four summands add up, respectively, to:
\begin{align*}
&\ \frac{1 + (-1)^m(2m-1)}{4} \lambda^{m-2} B_{21} B_{12}, \qquad
\binom{m}{2} (-\lambda)^{m-2} B_{22}^2,\\
&\ \frac{\mu^m - m \mu (-\lambda)^{m-1} + (m-1)
(-\lambda)^m}{(\mu+\lambda)^2} B_{23} B_{32},\\
&\ B_{24} \cdot (D^m - m (-\lambda)^{m-1} D + (m-1) (-\lambda)^m {\rm
id}_{W_4}) (D + \lambda)^{-2} \cdot B_{42}.
\end{align*}

Akin to $i=j=1$, the asymptotics we compute -- for
infinitely many odd or even $m_k$ -- are (see~\eqref{Ebse}):
\begin{equation}
\lim_{k \to \infty} \frac{1}{m_k^2 \lambda^{m_k-2}} S_{m_k,2}(A,B)_{22} =
\frac{\bse_2({\bf m})}{2} B_{22}^2.
\end{equation}

\item $(i,j)=(2,3)$:
The first three summands add up, respectively, to:
\begin{align*}
&\ \left( \frac{(-\lambda)^m}{2\lambda(\lambda+\mu)} +
\frac{\lambda^m}{2\lambda(\lambda-\mu)} +
\frac{\mu^m}{(\mu-\lambda)(\mu+\lambda)} \right) B_{21} B_{13},\\
&\ \frac{\mu^m - m \mu (-\lambda)^{m-1} + (m-1)
(-\lambda)^m}{(\mu+\lambda)^2} B_{22} B_{23}, \qquad
\frac{(-\lambda)^m + m \lambda \mu^{m-1} + (m-1) \mu^m}{(\mu+\lambda)^2}
B_{23} B_{33}.
\end{align*}

As for $(i,j)=(1,3)$, the fourth summand's sum can be bounded above in
norm by
\[
\frac{\lambda^m - m (\sqrt{\lambda_2})^{m-1} \lambda + (m-1)
(\sqrt{\lambda_2})^m}{(\lambda - \sqrt{\lambda_2})^2} \|B_{24}\| \,
\|B_{43}\|.
\]

From these, we obtain the asymptotics
\begin{equation}
\lim_{k \to \infty} \frac{1}{m_k \lambda^{m_k-2}} S_{m_k,2}(A,B)_{23} =
\bse_2({\bf m}) \frac{\lambda}{\lambda+\mu} B_{22} B_{23}.
\end{equation}

\item $(i,j)=(3,3)$:
The first three summands add up, respectively, to:
\[
\frac{\lambda^m - m \mu^{m-1} \lambda + (m-1) \mu^m}{(\lambda -
\mu)^2} B_{31} B_{13}, \qquad
\frac{(-\lambda)^m + m \mu^{m-1} \lambda + (m-1) \mu^m}{(\lambda +
\mu)^2} B_{32} B_{23},
\]
and $\binom{m}{2} \mu^{m-2} B_{33}^2$. The fourth can be bounded above in
norm (as for $(i,j)=(2,3)$) by $\binom{m}{2} (\sqrt{\lambda_2})^{m-2} \,
\|B_{34}\| \, \|B_{43}\|$.
This yields the asymptotics
\begin{equation}
\lim_{k \to \infty} \frac{1}{\lambda^{m_k-2}} S_{m_k,2}(A,B)_{33} =
\frac{\lambda^2}{(\lambda-\mu)^2} B_{31} B_{13} + \bse_2({\bf m})
\frac{\lambda^2}{(\lambda+\mu)^2} B_{32} B_{23}.
\end{equation}
\end{enumerate}

\noindent \textbf{Step 5:}
Taking stock of the above asymptotics, suppose $A^2 B \neq B A^2$ and
$S_{m,2}(A,B) \geq 0$ for infinitely many $m = m_k$ -- so we may assume
they are all odd or all even. Then so is
\begin{equation}\label{Eput0}
\frac{1}{\lambda^{m_k-2}} T S_{m_k,2}(A,B) T^*, \qquad \text{where}
\qquad T = \begin{pmatrix}
m_k^{-1} P_1 & 0 & 0 & 0 \\ 0 & m_k^{-1} P_2 & 0 & 0 \\ 0 & 0 & P_3 & 0
\end{pmatrix}.
\end{equation}

Taking limits and using the asymptotics in~\eqref{Egrowth11l2} and the
next five cases with all $m_k \geq 2$ odd or all even),
\begin{align}\label{El2C}
0 \leq &\ \lim_{k \to \infty} \frac{1}{\lambda^{m_k-2}} T S_{m_k,2}(A,B)
T^*\notag\\
&\ = \begin{pmatrix}
B_{11}^2/2 & 0 & \frac{\lambda}{\lambda-\mu} B_{11} B_{13}\\
0 & \bse_2({\bf m}) B_{22}^2/2 & \bse_2({\bf m})
\frac{\lambda}{\lambda+\mu} B_{22} B_{23}\\
\frac{\lambda}{\lambda-\mu} B_{31} B_{11} &
\bse_2({\bf m}) \frac{\lambda}{\lambda+\mu} B_{32} B_{22} &
\frac{\lambda^2}{(\lambda-\mu)^2} B_{31} B_{13} + \bse_2({\bf m})
\frac{\lambda^2}{(\lambda+\mu)^2} B_{32} B_{23}
\end{pmatrix}.
\end{align}

Denote this final operator-matrix by $C$.\medskip

\noindent \textbf{Step 6:}
In this step, we prove Theorem~\ref{Thilbert}(2) for infinitely many
$m_k$ of the same parity as $l$ (in this case, even since $l=2$). This
means that $\bse_2({\bf m}) = 1$. (The case with all $m_k-l$ odd reveals
a stronger conclusion: the remaining unproved implication in
Theorem~\ref{Tnxn}.)

We now invoke the assumption that $B \geq 0$ -- note this is needed,
because if $B_{11} = B_{22} = 0$, then $C$ is a block-diagonal
operator-matrix with positive entries regardless of what $B_{13}, B_{23}$
are, and there is no contradiction.

We will obtain a contradiction (since we assumed that $[A^2,B] \neq 0$)
when $B_{13} \neq 0 \leq B$; the proof is similar if $B_{23} \neq 0 =
B_{13}$. Informally, the matrix $C$ is a sum of two block $2\times 2$
Hermitian matrices (each padded by zeros), and both are non-positive. So
if e.g.\ $B_{11}^2, B_{22}^2$ are invertible, whence positive definite,
then the Schur complement $C'$ of the top $2 \times 2$ block of $C$ is
positive:
\begin{align*}
C' = \frac{\lambda^2}{(\lambda-\mu)^2} B_{31} B_{13} +
\frac{\lambda^2}{(\lambda+\mu)^2} B_{32} B_{23}
- &\ \frac{\lambda}{\lambda-\mu} B_{31} B_{11} (B_{11}^2/2)^{-1}
\frac{\lambda}{\lambda-\mu} B_{11} B_{13}\\
- &\ \frac{\lambda}{\lambda+\mu} B_{32} B_{22} (B_{22}^2/2)^{-1}
\frac{\lambda}{\lambda+\mu} B_{22} B_{23}\\
= &\ - \left( \frac{\lambda^2}{(\lambda-\mu)^2} B_{31} B_{13} +
\frac{\lambda^2}{(\lambda+\mu)^2} B_{32} B_{23} \right).
\end{align*}

Now if $B_{13} \neq 0$ then there exists $v_\mu$ such that $B_{13} v_\mu
\neq 0$. Hence
\[
0 \leq \tangle{v_\mu, C' v_\mu} = 
- \frac{\lambda^2}{(\lambda-\mu)^2} \| B_{13} v_\mu \|^2 -
\frac{\lambda^2}{(\lambda+\mu)^2} \| B_{23} v_\mu \|^2 < 0,
\]
and this contradiction finishes the proof. Similarly if $B_{23} \neq 0$.

We make this argument rigorous even if $B_{11}^2, B_{22}^2 \geq 0$ are
not invertible. Since $C \geq 0$, we have $C_\varepsilon := C +
\varepsilon \, {\rm id}_{W_4^\perp} > 0$ for all $\varepsilon > 0$. Hence
the top $2 \times 2$ block of $C_\varepsilon$ is positive definite,
whence by Lemma~\ref{Laitken}, the Schur complement -- call it
$C'_\varepsilon$ -- is positive. That is:
\begin{align}\label{Eschurblock}
0 \leq \lambda^{-2} C'_\varepsilon = \frac{\varepsilon}{\lambda^2} {\rm
id}_{W_3} + &\ \frac{1}{(\lambda-\mu)^2} B_{31} \left( {\rm id}_{W_1} -
B_{11} (\textstyle{\frac{B_{11}^2}{2}} + \varepsilon \, {\rm
id}_{W_1})^{-1} B_{11} \right) B_{13}\\
+ &\ \frac{1}{(\lambda+\mu)^2} B_{32} \left( {\rm id}_{W_2} - B_{22}
(\textstyle{\frac{B_{22}^2}{2}} + \varepsilon \, {\rm id}_{W_2})^{-1}
B_{22} \right) B_{23}\notag
\end{align}

Now assume $B_{13} \neq 0$; the analysis is similar if $B_{23} \neq 0 =
B_{13}$. Choose $v_\mu \in W_3 = V_\mu$ such that $v := B_{13} v_\mu \neq
0$. Then Lemma~\ref{LB12} (applied to $B_{11} \geq 0$) enables us to
apply Lemma~\ref{Lspectralcompute} on $W_1$, with parameters as
in~\eqref{Eapplylemma} for $l=2$. This yields that $\tangle{v,
f_\varepsilon(B_{11})v} \to -\|v\|^2< 0$ as $\varepsilon \to 0^+$. Now
there are two cases:
\begin{enumerate}
\item If $B_{23}v_\mu = 0$ then this yields the desired contradiction:
applying $\tangle{v_\mu, - \cdot v_\mu}$ to both sides
of~\eqref{Eschurblock} yields:
\[
0 \leq \frac{\varepsilon}{\lambda^2} \|v_\mu\|^2 -
\frac{1}{(\lambda-\mu)^2} \|v\|^2 \qquad \text{as} \qquad \varepsilon \to
0^+.
\]

\item Else $B_{23} v_\mu \neq 0$. Again use Lemma~\ref{LB12} (for $B_{22}
\geq 0$) to apply Lemma~\ref{Lspectralcompute} with
$(\alpha, \beta, a, b, c)$ as before with $l=2$, and
$(\mathcal{H}, B, v) = (W_2, B_{22}, B_{23}v_\mu)$.
This implies: $\tangle{B_{23}v_\mu, f_\varepsilon(B_{22}) B_{23}v_\mu}
\to -\|B_{23}v_\mu\|^2< 0$ as $\varepsilon \to 0^+$. Now applying
$\tangle{v_\mu, - \cdot v_\mu}$ to both sides of~\eqref{Eschurblock}
yields:
\[
0 \leq \frac{\varepsilon}{\lambda^2} \|v_\mu\|^2 -
\frac{1}{(\lambda-\mu)^2} \|v\|^2 - \frac{1}{(\lambda+\mu)^2}
\|B_{23}v_\mu\|^2 \qquad \text{as} \qquad \varepsilon \to 0^+.
\]
\end{enumerate}
In either case, we obtain a contradiction, and the proof is complete:
$A^2B = BA^2$.\medskip

\noindent \textbf{Step 7:}
The final case is if all $m_k - l = m_k - 2$ are odd, so that
$\bse_2({\bf m}) = -1$. Then using the asymptotics for the matrix $C \geq
$ in~\eqref{El2C}, we have
\begin{equation}\label{Econtra}
0 \leq l! \, P_2 C P_2 = -B_{22}^l, \qquad 0 \leq P_2 B P_2 = B_{22}.
\end{equation}
Together, these imply $B_{22} = 0$. But akin to Section~\ref{S2impl},
this implies $B P_2 = 0 = P_2 B$, so $B_{32} = B_{23} = 0$, and the
matrix $C$ has only ``corner'' entries nonzero. As we assumed (for
contradiction) $[A^2,B] \neq 0$, we must have $B_{13} \neq 0$. But then
the working in the previous step (see case~(1)) yields the desired
contradiction. Hence $A^2B=BA^2$.

Returning to the notation in Step~1, $[A^2,B] = 0$ implies that $B$
preserves every $V'_{\lambda^2}$ (and in particular, $V'_0 = V_0 =
\ker(A)$ too). Thus, we may now restrict the action of $B$ to each
$V'_{\lambda^2} = V_\lambda \oplus V_{-\lambda}$ for $0 < \lambda$. On
this, the same computations used to obtain~\eqref{El2C}
show~\eqref{Econtra} again. Hence $0 = B_{22} = B P_2 = B P_{-\lambda} =
0$ -- for every $0 < \lambda \in \sigma(A)$ such that $-\lambda \in
\sigma(A)$. This implies $B P_{\mathcal{H}_-} = 0$. Moreover, we get that
$B_{21} = 0$ by Lemma~\ref{LB12}, so $B$ preserves not just $V_\lambda
\oplus V_{-\lambda}$, but in fact each summand too. Hence $AB=BA$. \qed

\subsection{Tightness of the hypotheses for $l=2$}\label{S2tight}

Here we work in parallel to Section~\ref{S1tight}. We begin with:

\begin{prop}\label{P2tight}
Suppose $\mathcal{H}$ is a Hilbert space and $A,B \in
\mathcal{B}(\mathcal{H})$ are self-adjoint. If $A^2 B = B A^2$, $B \geq
0$, and $S_{m_-,2}(A,B) \geq 0$ for one odd $m_- \geq 2$, then $AB=BA$.
\end{prop}

\begin{proof}
Write $A,B$ in the block-form~\eqref{EAB3x3} using the spectral
projections~\eqref{Espectralproj}. As discussed there, $B$ commutes with
$A^2$ and hence with $|A|$, and so $c=e=0$ in~\eqref{EAB3x3}. Moreover,
writing $A = |A|U = U |A|$ with $U = {\rm diag}({\rm id}_{W_1}, -{\rm
id}_{W_2}, {\bf 0}_{W_3})$, $|A|$ commutes with $U$ too. Thus,
\[
0 \leq S_{m_-,2}(A,B) = |A|^{(m_- -2)/2} S_{m_-,2}(U,B) |A|^{(m_- -2)/2}.
\]

Hence the $(2,2)$-block of this operator-matrix is positive, and it equals
\[
|A_-|^{(m_- - 2)/2} \cdot S_{m_-,2}(U,B)_{22} \cdot |A_-|^{(m_- - 2)/2}.
\]

By suitably applying Remark~\ref{Rdense}, $S_{m_-,2}(U,B)_{22} \geq 0$.
As $m_-$ is odd, an explicit computation using~\eqref{Easymptotics} (or
its $3 \times 3$ analogue) reveals:
\[
S_{m_-,2}(U,B)_{22} = \sum_{\substack{\alpha, \beta, \gamma \in
\mathbb{Z}_{\geq 0},\\ \alpha+\beta+\gamma = m_- -2}}
(-1)^{\alpha+\gamma} (b^* b + (-1)^\beta d^2)
= -\binom{m_-}{2} d^2 - \frac{m_- - 1}{2} b^* b.
\]

Because this is positive as well as negative, we obtain $b^* b = 0$,
whence $b=0$. As $c=e=0$ as well in~\eqref{EAB3x3}, we get $AB=BA$.
\end{proof}

We now present a result as in the title of this subsection (for $l=2$,
parallel say to Section~\ref{S1tight} for $l=1$). Suppose $l=2$ and $A,B
\in \mathcal{B}(\mathcal{H})$ are self-adjoint, with $B \geq 0$ and
$\sigma(A) \setminus \{ 0 \}$ discrete. We have seen multiple sufficient
conditions above to ensure $AB = BA$: $S_{m,2}(A,B) \geq 0$ for
infinitely many $m \geq 2$, and either
\begin{itemize}
\item one odd $m>2$, or
\item $\pm A \geq 0$, or the strictly weaker assumption: $\sigma(A) \cap
-\sigma(A) \subseteq \{ 0 \}$.
\end{itemize}

We end this section by proving that these sufficient conditions are tight:

\begin{prop}
Suppose $\mathcal{H}$ is a Hilbert space, and choose any $A = A^* \in
\mathcal{B}(\mathcal{H})$ such that $\sigma(A) \setminus \{ 0 \}$ is
discrete and $\sigma(A) \cap -\sigma(A) \nsubseteq \{ 0 \}$. Then there
exists $0 \leq B \in \mathcal{B}(\mathcal{H})$ such that
(a)~$AB \neq BA$ and
(b)~$S_{m,2}(A,B)$ is positive for all even $m \geq 2$,
(c)~whence $S_{m,2}(A,B)$ is not positive for any odd $m \geq 2$ by
Proposition~\ref{P2tight}.
\end{prop}

\begin{proof}
Choose any $0 < \lambda \in \sigma(A) \cap -\sigma(A)$, with orthonormal
eigenvectors $v_\pm \neq 0$ such that $A v_\pm = \pm \lambda v_\pm$. Let
$\mathcal{H}_0 := {\rm span}(v_+, v_-)^\perp$ and henceforth write
operators as block $3 \times 3$ matrices, corresponding to $(v_+, v_-,
\mathcal{H}_0)$. With respect to this decomposition of $\mathcal{H}$, we
have
\[
A := \begin{pmatrix} \lambda & 0 \\ 0 & -\lambda \end{pmatrix} \oplus
D_{\mathcal{H}_0}, \qquad \text{and now fix} \qquad B := \begin{pmatrix}
1 & 1 \\ 1 & 1 \end{pmatrix} \oplus {\bf 0}_{\mathcal{H}_0} \geq 0.
\]

Clearly, $AB \neq BA$ because they both fix $\mathcal{H}_0^\perp$ and do
not commute on it. Moreover,
\[
S_{m,2}(A,B) = S_{m,2}(A',B') \oplus {\bf 0}_{\mathcal{H}_0}, \qquad
\text{where} \qquad A' := \begin{pmatrix} \lambda & 0 \\ 0 & -\lambda
\end{pmatrix}, \ B' := {\bf 1}_{2 \times 2}.
\]
(Note also by direct inspection that $A^2$ commutes with $B$.)
Now compute explicitly using~\eqref{Easymptotics} (or its $2 \times 2$
analogue):
\begin{align*}
S_{m,2}(A',B') = &\ \sum_{\substack{\alpha, \beta, \gamma \in
\mathbb{Z}_{\geq 0},\\ \alpha+\beta+\gamma = m-2}}
(A')^\alpha B' (A')^\beta B' (A')^\gamma\\
= &\ \sum_{\substack{\alpha, \beta, \gamma \in
\mathbb{Z}_{\geq 0},\\ \alpha+\beta+\gamma = m-2}}
\begin{pmatrix}
\lambda^{m-2} + \lambda^{\alpha+\gamma} (-\lambda)^\beta &
\lambda^{\alpha+\beta} (-\lambda)^\gamma + \lambda^\alpha
(-\lambda)^{\beta+\gamma} \\
\lambda^{\beta+\gamma} (-\lambda)^\alpha + \lambda^\gamma
(-\lambda)^{\alpha+\beta} &
\lambda^\beta (-\lambda)^{\alpha+\gamma} + (-\lambda)^{m-2} 
\end{pmatrix}.
\end{align*}

Letting $m$ be even, apply~\eqref{Eapgp} to obtain:
\[
= \frac{\lambda^{m-2}}{4} \begin{pmatrix}
4 \binom{m}{2} + 2m & 4m \\ 4m & 4 \binom{m}{2} + 2m
\end{pmatrix} = \frac{m \lambda^{m-2}}{2} \begin{pmatrix} m & 2 \\ 2 & m
\end{pmatrix},
\]
and this is positive for even $m \geq 2$.
\end{proof}

\subsection{Proof of Theorem~\ref{Tnxn}$(2) \implies (3)$ and
Theorem~\ref{Thilbert}(2) -- for $m \geq 0, l \geq 3$}

Finally, we resolve the case $l \geq 3$. The setting is once again as in
Step~1 of the preceding Section~\ref{Squad}, with $B =
(B_{ij})_{i,j=1}^4$ and $(B_{13}, B_{23}) \neq (0,0)$, so that $A^2B \neq
BA^2$.

\begin{remark}
Here is a \textbf{notational change} from the $l=2$ case: henceforth in
this section, we work with $\lambda_1 = \lambda$ and $\lambda_2$ in
$\sigma(A)$, with $|\lambda_1| > |\lambda_2|$. (Recall that $\lambda_1,
\lambda_2$ were previously in $\sigma(A^2)$.)
\end{remark}

We will use without mention, the following lemma.

\begin{lemma}\label{Lclever}
Suppose $B_{13} : V_\mu \to V_{\lambda_1}$ is nonzero. Choose $v_\mu \in
V_\mu$ such that $v := B_{13} v_\mu$ is nonzero. Then (noting that $B
\geq 0$,) $\tangle{v, B_{11}^\beta v} > 0$ for all scalars $\beta \in
(0,\infty)$.
\end{lemma}

\begin{proof}
Since the $4 \times 4$ block-partitioned operator-matrix $B$ is positive,
so is
\[
\begin{pmatrix} P_1 & 0 & 0 & 0 \\ 0 & 0 & v_\mu^* & 0 \end{pmatrix}
\cdot B \cdot \begin{pmatrix} P_1 & 0 \\ 0 & 0 \\ 0 & v_\mu \\ 0 & 0
\end{pmatrix} = \begin{pmatrix} B_{11} & v \\ v^* & \alpha \end{pmatrix},
\quad \text{where} \quad \alpha = \tangle{v_\mu, B_{33} v_\mu} \in [0,
\infty).
\]

We first claim that $B_{11} v \neq 0$. Indeed, if not then for
$\varepsilon > 0$ we compute:
\[
0 \leq \varepsilon^{-1} \cdot \begin{pmatrix} v^* & -\varepsilon
\end{pmatrix} \cdot \begin{pmatrix} B_{11} & v \\ v^* & \alpha
\end{pmatrix} \cdot \begin{pmatrix} v \\ -\varepsilon \end{pmatrix} = -2
\|v\|^2 + \varepsilon \alpha,
\]
and this is negative for small $\varepsilon > 0$, which is impossible.
This contradiction shows the claim.

Now $B_{11}^\beta$ is a positive operator by Spectral mapping (see
Theorem~\ref{Tfnalcalc}(4) and Corollary~\ref{Cmapping}), so $\tangle{v,
B_{11}^\beta v} \geq 0$. Suppose for contradiction that it is zero. Then
$C v = 0$, where $C = B_{11}^{\beta/2}$ (defined via the continuous
functional calculus). Now choose a sequence $p_k(x)$ of polynomials that
converges uniformly to $x^{2/\beta}$ on the compact set $[0,\|C\|] =
[0,\|B_{11}\|^{\beta/2}]$. Define $q_k(x) := p_k(x) - p_k(0)$; then
$q_k(x) \to x^{2/\beta}$ uniformly, and $x | q_k(x)$ for all $k$.

Since $Cv = 0$, it follows that $q_k(C)v = 0$ for all $k$. Taking $k \to
\infty$, $0 = C^{2/\beta} v = B_{11}v$, and this contradicts the first
claim above.
\end{proof}

Finally, we begin the proof for $l \geq 3$. We will assume
$S_{m_k,l}(A,B)$ is positive for infinitely many integers $l \leq m_1 <
m_2 < \cdots$. This proof too is split into steps.\medskip

\noindent \textbf{Step 1:}
For fixed $l \geq 1$, define the generating function and the resolvent
\begin{equation}\label{Eresolvent}
F_l(z) := \sum_{m \geq l} z^m S_{m,l}(A,B), \qquad R(z) := ({\rm Id} -
zA)^{-1}.
\end{equation}

Note that the resolvent is defined for $\| zA \| < 1$, and in this case
the partial sums $\sum_{p=0}^N z^p A^p$ converge to $R(z)$ in the
operator norm topology. Moreover, by the above notation,
\[
R(z) = \begin{pmatrix}
\frac{1}{1-z \lambda_1} P_1 & 0 & 0 & 0 \\
0 & \frac{1}{1+z \lambda_1} P_2 & 0 & 0 \\
0 & 0 & \frac{1}{1-z \mu} P_3 & 0 \\
0 & 0 & 0 & (P_4 - zD)^{-1}
\end{pmatrix}
\]
and this is in $\mathcal{B}(\mathcal{H})$ for $|z| < 1/\|A\| =
1/\lambda_1$. Now we rewrite $F_l$:

\begin{prop}\label{Pgenfn}
Let $A,B$ be bounded linear operators on a Banach space $\mathbb{B}$.
Then for all $z$ in the disc $D(0, 1/\| A \|)$, we have an equality of
operators $F_l(z) \equiv G_l(z)$, where $F_l$ is as
in~\eqref{Eresolvent}, and
\begin{equation}\label{EGl}
G_l(z) = z^l \cdot R(z) B R(z) \cdots B R(z) = z^l R(z) (B R(z))^l.
\end{equation}
\end{prop}

\noindent (This means that for all $l, z$ as above, the operators
$F_l(z), G_l(z)$ exist; and for all vectors $v \in \mathbb{B}$ and
functionals $\varphi \in \mathbb{B}^*$, we have $\varphi( F_l(z) v) =
\varphi( G_l(z) v)$.) The proof is standard, and is deferred to
Appendix~\ref{Appendix3}.\medskip

\noindent \textbf{Step 2:}
We now examine the operator-valued map $G_l$ more closely. The next few
paragraphs are standard; see e.g.\ \cite[Chapter~5]{Theo}, and
specifically, Lemma~5.1.12 therein -- whose special case shows that
``weak holomorphic = strong holomorphic = norm-holomorphic''. That is,
operator-valued holomorphic functions mean the same thing in the norm
topology on $\mathcal{B}(\mathcal{H})$ and in the SOT/WOT (strong/weak
operator topologies). Moreover, all contour integrals are Bochner
integrals, and Cauchy's formula and the residue theorem will again be
applied in the Banach space $\mathcal{B}(\mathcal{H})$.

First note that $R(z)$ can be written in ``block diagonal'' form with
respect to the eigenspace decomposition (for $A^2$) $\mathcal{H} = V'_0
\oplus \bigoplus_{n \geq 1} V'_{\lambda_n}$ as ${\rm diag}(P_{V_0};
\frac{1}{1-z \lambda_n} P_{V_{\lambda_n}}, \frac{1}{1+z \lambda_n}
P_{V_{-\lambda_n}})$. Thus, it is (weakly/strongly) holomorphic on
$\mathbb{C} \setminus \{ \pm 1/\eta : 0 \neq \eta \in \sigma(A) \}$.

We will henceforth work inside the open disc $D(0, 1/\lambda_2)$ (if
$\lambda_2 = \mu = 0$ then this is all of $\mathbb{C}$). On this domain,
define

\begin{equation}\label{EHlGl}
S_\pm(z) := (1 \pm z \lambda_1) R(z), \qquad
H_{\pm,l}^{(1)}(z) := S_\pm(z) (B S_\pm(z))^l, \qquad
H_{\pm,l}(z) := z^l H_{\pm,l}^{(1)}(z).
\end{equation}

All three of these are meromorphic on $D(0, 1/\lambda_2)$ -- whence so
are $G_l^{(1)}(z) := R(z) (BR(z))^l$ and $G_l(z) = z^l G_l^{(1)}(z)$ --
with:
\begin{itemize}
\item $S_\pm(z)$ having simple poles at $\pm 1/\lambda_1$,
\item $H_{\pm,l}^{(1)}(z)$ therefore having poles of order at most $l+1$
at $\pm 1/\lambda_1$, and
\item $G_l^{(1)}(z), G_l(z)$ having poles of order at most $l+1$ at both
$\pm 1/\lambda_1$.
\end{itemize}

In particular, $G_l$ is holomorphic (in the SOT) on $D(0, 1/\lambda_1)$,
and agrees there with $F_l$, by Proposition~\ref{Pgenfn}. Now choose $R
\in (1/\lambda_1, 1/\lambda_2)$; if $\lambda_2=0$ (in which case
$\mu=0$), we choose $R > 1/\lambda_1$. We will henceforth work with
$G_l$, and on the circular contours $C_R := \partial D(0,R)$
and $C_{(1/\lambda_1)^-}$, which we define as $C_r$ for any $0 < r <
1/\lambda_1$. As is standard, we denote the latter by $C_{1/\lambda}$.

As is also standard, for $m \geq 0$ we denote the coefficient of $z^m$ in
a (operator-matrix valued) power series $A(z) := \sum_{m \geq 0} z^m A_m$
by
\[
[z^m]A(z) = A_m.
\]

We are interested in $[z^m]F_l(z) = [z^m]G_l(z)$ (viewing these as
holomorphic maps in $D(0, 1/\|A\|)$), and in the growth of this sequence
as $m \to \infty$. But now note by Cauchy's integral and residue formulas
(see again \cite[Lemma~5.1.12]{Theo}(3)) that
\begin{align}\label{Econtour}
[z^m]G_l(z) = &\ \frac{1}{2 \pi i} \int_{C_{1/\lambda_1}}
\frac{G_l(z)}{z^{m+1}}\, dz,\notag\\
\Res \left[ \frac{G_l(z)}{z^{m+1}}, \, z = \frac{1}{\lambda_1} \right] +
\Res \left[ \frac{G_l(z)}{z^{m+1}}, \, z = \frac{-1}{\lambda_1} \right] = &\
\frac{1}{2\pi i} \int_{C_R} \frac{G_l(z)}{z^{m+1}}\, dz -
\frac{1}{2\pi i} \int_{C_{1/\lambda_1}} \frac{G_l(z)}{z^{m+1}}\, dz,
\end{align}
where $\Res[ \cdots ]$ means that we take residues of the operator-matrix
map $G_l(z)/z^{m+1}$ ``entrywise'' -- i.e., in the sense of WOT.\medskip

\noindent \textbf{Step 3:}
We now upper-bound the first term in~\eqref{Econtour} on the right --
henceforth called the ``$C_R$-term'' in~\eqref{Econtour}. Recall that
$G_l(z) = z^l G^{(1)}_l(z)$, with $G_l^{(1)}(z)$ (weakly/strongly)
holomorphic in a small annular neighborhood of $C_R$ by Step~2. In
particular, $G^{(1)}_l$ is continuous on $C_R$, so any
``entry/evaluation'' of the $C_R$-term in~\eqref{Econtour} against
arbitrary $g,h \in \mathcal{H}$ is bounded by:
\begin{align*}
\left| \left\langle g, \frac{1}{2 \pi i} \int_{C_R}
\frac{G_l(z)}{z^{m+1}}\, dz \cdot h\right\rangle \right| := &\ \left|
\frac{1}{2 \pi i} \int_{C_R} \frac{\langle g, G_l(z) h\rangle}{z^{m+1}}\,
dz \right| \leq \frac{1}{2 \pi} \int_{C_R} \frac{|\tangle{g,
G^{(1)}_l(z)h}|}{|z|^{m+1-l}}\, |dz|\\
&\ \leq \frac{\| \tangle{g, G^{(1)}_l(z)h} \|_{\infty,C_R}}{R^{m-l}}
< \| \tangle{g, G^{(1)}_l(z)h} \|_{\infty,C_R} \cdot \lambda_1^{m-l}.
\end{align*}
The final inequality is because $R > 1/\lambda_1$.\medskip

\noindent \textbf{Step 4:}
Henceforth in this proof, $H_{\mp l}^{(p)}(z)$ will denote the $p$th
derivative of $H_{\mp,l}(z)$; in particular, $H_{\mp,l}^{(1)}(z)$ is not
the function studied above.
We next compute the residue term in~\eqref{Econtour}. Using~\eqref{EGl},
we alternately write
$\displaystyle G_l(z) = \frac{H_{\mp,l}(z)}{(1 \mp z\lambda_1)^{l+1}}$.
Thus,
\begin{align}\label{EHl}
H_{\mp,l}(z) = &\ z^l S_\mp(z) (B S_\mp(z))^l, \ \text{where e.g.}\\
S_-(z) := &\ (1-z\lambda_1)R(z) = \begin{pmatrix} P_1 & 0 & 0 & 0 \\
0 & \frac{1-z\lambda_1}{1+z\lambda_1} P_2 & 0 & 0 \\
0 & 0 & \frac{1-z\lambda_1}{1-z\mu} P_3 & 0 \\
0 & 0 & 0 & (1-z\lambda_1) (P_4 - zD)^{-1} \end{pmatrix}\notag
\end{align}
is an operator-matrix valued map, and similarly for $S_+(z)$. The residue
at the pole $1/\lambda_1$ for the function $G_l(z) / z^{m+1}$ is:
\[
\Res \left[ \frac{G_l(z)}{z^{m+1}}, \, z = \frac{1}{\lambda_1}
\right] = \frac{1}{l!} \lim_{z \to 1/\lambda_1} \frac{d^l}{dz^l} \left[
(z - 1/\lambda_1)^{l+1} \frac{G_l(z)}{z^{m+1}} \right]\\
= \frac{1}{l! \, (-\lambda_1)^{l+1}} \lim_{z \to 1/\lambda_1}
\frac{d^l}{dz^l} \frac{H_{-,l}(z)}{z^{m+1}}.
\]

By the Leibniz rule and since $\displaystyle \frac{d^p}{dz^p}
(z^{-(m+1)}) = \frac{(-1)^p}{z^{m+p+1}} \frac{(m+p)!}{m!}$, we obtain:
\begin{align}\label{Eleibniz}
\Res \left[ \frac{G_l(z)}{z^{m+1}}, \, z = \frac{1}{\lambda_1} \right] =
&\ \frac{1}{l! \, (-\lambda_1)^{l+1}} \sum_{p=0}^l \lim_{z \to
1/\lambda_1} \binom{l}{p} H_{-,l}^{(l-p)}(z) \frac{(-1)^p}{z^{m+p+1}}
\frac{(m+p)!}{m!}\notag\\
= &\ \sum_{p=0}^l \frac{(-1)^{l+1-p}}{(l-p)!}
H_{-,l}^{(l-p)}(1/\lambda_1) \cdot \binom{m+p}{p} \lambda_1^{m - (l-p)}.
\end{align}

One similarly computes the residue of $G_l(z)/z^{m+1}$ at $z =
-1/\lambda_1$, and we have:
\begin{equation}\label{Eleibniz2}
\Res \left[ \frac{G_l(z)}{z^{m+1}}, \, z = \frac{\pm 1}{\lambda_1}
\right] = \sum_{p=0}^l \frac{(\pm 1)^{m-(l-p)} (-1)^{l+1-p}}{(l-p)!}
H_{\mp,l}^{(l-p)}(\pm 1/\lambda_1) \cdot \binom{m+p}{p} \lambda_1^{m -
(l-p)}.
\end{equation}

\noindent \textbf{Step 5:}
The goal is to compute the $3 \times 3$ asymptotics of $(P_1 + P_2 + P_3)
S_{m,l}(A,B) (P_1 + P_2 + P_3)$, as in the proof of the $l=2$ case
earlier in this section. In this step we report the $(i,j)$-block
asymptotics for $1 \leq i,j \leq 2$. Since in~\eqref{Eleibniz2} the
fastest growth-rate is $m^l \lambda_1^m$ (for $p=l$), we have
\[
\lim_{k \to \infty} \frac{1}{m_k^l \lambda_1^{m_k-l}}
\Res \left[ \frac{G_l(z)}{z^{m_k+1}}, \, z = \frac{\pm 1}{\lambda_1}
\right] = -\frac{\lambda_1^l \lim_{k \to \infty} (\pm 1)^{m_k}}{l!}
H_{\mp,l}(\pm 1/\lambda_1),
\]
along any sequence $l \leq m_k \uparrow \infty$.

But $S_-(1/\lambda_1) = P_1$ and $S_+(-1/\lambda_1) = P_2$, which helps
us compute $H_{\mp,l}$ and its derivatives at $\pm 1/\lambda_1$. First
define the $j$th block-column of $B$ and the first two rows of ${\rm
id}_{\mathcal{H}}$ via
\begin{equation}\label{Edefine2}
B^{(j)} := \begin{pmatrix} B_{j1} & B_{j2} & B_{j3} & B_{j4}
\end{pmatrix}^*, \quad E_1 := \begin{pmatrix} P_1 & 0 & 0 & 0
\end{pmatrix}, \quad E_2 := \begin{pmatrix} 0 & P_2 & 0 & 0
\end{pmatrix},
\end{equation}
for $1 \leq j \leq 4$. Then using~\eqref{EHlGl}, we obtain:
\[
H_{+,l}(-1/\lambda_1) = \frac{1}{(-\lambda_1)^l} E_2^* E_2 (B E_2^*
E_2)^l = \frac{1}{(-\lambda_1)^l} E_2^* (E_2 B E_2^*)^l E_2 =
\frac{1}{(-\lambda_1)^l} P_2 B_{22}^l P_2,
\]
and similarly,
$H_{-,l}(1/\lambda_1) = \frac{1}{\lambda_1^l} P_1 B_{11}^l P_1$.

Put both these into~\eqref{Eleibniz2}, and then into~\eqref{Econtour}.
Since the $C_R$-term grows as $\lambda_1^{m-l}$ and here we divide by
$m^l \lambda_1^{m-l}$, it follows that
\begin{equation}\label{Eput1}
\lim_{k \to \infty} \frac{1}{m_k^l \lambda_1^{m_k-l}}[z^{m_k}]F_l(z) = 0
+ \frac{1}{l!} P_1 B_{11}^l P_1 + \lim_{k \to \infty}
\frac{(-1)^{m_k-l}}{l!} P_2 B_{22}^l P_2.
\end{equation}

Now using the hypotheses of infinitely many $m_k$ having the same
parity (see~\eqref{Ebse}),
\begin{equation}\label{Egrowth11}
\lim_{k \to \infty} \frac{1}{m_k^l
\lambda_1^{m_k-l}}[z^{m_k}]F_l(z) = \frac{1}{l!} \begin{pmatrix} B_{11}^l
& 0 & 0 & 0 \\ 0 & \bse_l({\bf m}) B_{22}^l & 0 & 0 \\ 0 & 0 & 0 & 0 \\ 0
& 0 & 0 & 0
\end{pmatrix}.
\end{equation}

\noindent \textbf{Step 6:}
We next compute the asymptotics of the ``off-diagonal $(3,1)$ and $(3,2)$
blocks''. Since the $p=l$ fastest growing term does not contribute (see
the end of the previous step), we examine the $p=l-1$ term. For this,
denote
\[
f(z) := z^l {\rm Id}_{\mathcal{H}}, \qquad g_\pm(z) = S_\pm(z), \qquad
h_\pm(z) = h_{\pm,l}(z) := (B S_\pm(z))^l.
\]

Thus $H_{\pm,l} = f g_\pm h_\pm$, and one computes:
\begin{align*}
H_{-,l}'(1/\lambda_1) = &\
\begin{pmatrix} P_1 & 0 & 0 & 0 \\ 0 & 0 & 0 & 0 \\ 0 & 0 & 0 & 0 \\ 0 &
0 & 0 & 0 \end{pmatrix}
(l \lambda_1^{-(l-1)} \cdot h_-(1/\lambda_1) + \lambda_1^{-l}
h'_-(1/\lambda_1))\\
&\ + \lambda_1^{-l} \begin{pmatrix} 0 & 0 & 0 & 0 \\
0 & \frac{\lambda_1^2}{-2\lambda_1} P_2 & 0 & 0\\
0 & 0 & \frac{\lambda_1^2}{\mu-\lambda_1} P_3 & 0 & 0\\
0 & 0 & 0 & S'(1/\lambda_1)_{44}
\end{pmatrix} h_-(1/\lambda_1),
\end{align*}
where the second $4 \times 4$ block-matrix displayed is
$g'_-(1/\lambda_1)$. Clearly, the $(3,j)$-blocks in the first expression
on the right are zero for $j=1,2,3$ -- as are $h_-(1/\lambda_1)_{(3,2)}$
and $h_-(1/\lambda_1)_{(3,3)}$ by~\eqref{EHminus}. Therefore,
\begin{equation}\label{E312minus}
H'_{-,l}(1/\lambda_1)_{\{ 3 \} \times \{ 1,2,3 \}} = \left(
\lambda_1^{-l} \frac{\lambda_1^2}{\mu-\lambda_1} B_{31} B_{11}^{l-1}
\quad 0 \quad 0 \right).
\end{equation}

\noindent \textbf{Step 7:}
Similarly, one computes:
\begin{equation}\label{E312plus}
H'_{+,l}(-1/\lambda_1)_{\{ 3 \} \times \{ 1,2,3 \}} = \left( 0 \quad
(-\lambda_1)^{-l} \frac{\lambda_1^2}{\mu+\lambda_1} B_{32} B_{22}^{l-1}
\quad 0 \right).
\end{equation}

Using this and~\eqref{E312minus}, and isolating the $p=l-1$ term
in~\eqref{Eleibniz2}, we get:
\begin{align}\label{Eput2}
\lim_{k \to \infty} \frac{1}{m_k^{l-1}
\lambda_1^{m_k-l}}[z^{m_k}]F_l(z)_{31} = &\ 0 -
\frac{\lambda_1^{l-1}}{(l-1)!} H'_{-,l}(1/\lambda_1)_{31} - 0
= \frac{\lambda_1}{\lambda_1 - \mu}
\frac{B_{31} B_{11}^{l-1}}{(l-1)!},\notag\\
\lim_{k \to \infty} \frac{1}{m_k^{l-1}
\lambda_1^{m_k-l}}[z^{m_k}]F_l(z)_{32} = &\ 0 - 0 - 
\frac{\lambda_1^{l-1} \lim_{k \to \infty} (-1)^{m_k-1}}{(l-1)!}
H'_{+,l}(-1/\lambda_1)_{32}\\
&\ = \frac{\lambda_1}{\lambda_1+\mu} \bse_l({\bf m}) \frac{B_{32}
B_{22}^{l-1}}{(l-1)!}.\notag
\end{align}

Here we note that these asymptotics kill the $(3,1)$- and $(3,2)$- blocks
of the $C_R$-term (as in Step~5), since $l \geq 2$. Now assuming that
(infinitely many) $m_k$ have the same parity, we obtain the asymptotics
$\bse_l({\bf m}) \frac{\lambda_1}{\lambda_1 + \mu} \frac{B_{32}
B_{22}^{l-1}}{(l-1)!}$ in the $(3,2)$-block.\medskip

\noindent \textbf{Step 8:}
Finally, we compute the asymptotics of the $(3,3)$-block of
$S_{m,l}(A,B)$. Note by the preceding steps that the $p=l$ and $p=l-1$
terms do not contribute to the $(3,3)$-block, so we now consider the
$p=l-2$ term. First write
\begin{equation}\label{E2ndder}
H_{\pm,l}'' = (fg_\pm h_\pm)'' = (fg_\pm)'' h_\pm + 2 (fg_\pm)' h_\pm' +
fg_\pm \cdot h_\pm''.
\end{equation}

We need to evaluate the $(3,3)$-block of this at $\mp 1/\lambda_1$. Note
using~\eqref{Edefine2} that
\begin{equation}\label{EHminus}
h_-(1/\lambda_1) = B^{(1)} B_{11}^{l-1} E_1, \qquad
h_+(-1/\lambda_1) = B^{(2)} B_{22}^{l-1} E_2,
\end{equation}
and both of these have third block-column zero. Similarly, each of
$g_\pm(\mp 1/\lambda_1)$ has third row zero. Hence by~\eqref{E2ndder},
the $(3,3)$ block of $H_{\pm,l}''(\mp 1/\lambda_1)$ equals
\[
H_{\pm,l}''(\mp1/\lambda_1)_{33} = \left(2 f(\mp 1/\lambda_1) \cdot
g_\pm'(\mp 1/\lambda_1) \cdot h_\pm'( \mp 1/\lambda_1) \right)_{33} =
2 (\mp \lambda_1)^{-l} \cdot \frac{\lambda_1^2}{\mu \pm \lambda_1} \cdot
h_\pm'(\mp 1/\lambda_1)_{33}.
\]

But $h_\pm'(\mp 1/\lambda_1)_{33}$ is a Hurwitz product:
\[
h'_\pm(\mp 1/\lambda_1) = \left. [(B \cdot S_\pm(z))^l]' \right|_{z=\mp
1/\lambda} = S_{l,1}(BS_\pm(\mp 1/\lambda_1), BS'_\pm(\mp 1/\lambda_1)),
\]
and $B S_\pm(\mp 1/\lambda_1)$ has $(i,3)$-block zero for all $i$. Hence
if it is the last letter in a word then that word has $(3,3)$-block zero.
Thus, the only word in this Hurwitz product that contributes, gives:
\[
h'_\pm(\mp 1/\lambda_1)_{33} = ((BS_\pm(\mp 1/\lambda_1))^{l-1}
BS'_\pm(\mp 1/\lambda_1))_{33},
\]
and an explicit computation now yields:
\[
h'_-(1/\lambda_1)_{33} = B_{31} B_{11}^{l-2} B_{13}
\frac{\lambda_1^2}{\mu-\lambda_1}, \qquad
h'_+(-1/\lambda_1)_{33} = B_{32} B_{22}^{l-2} B_{23}
\frac{\lambda_1^2}{\mu+\lambda_1}.
\]
From above, this yields:
\[
H_{-,l}''(1/\lambda_1)_{33} = 2 \frac{\lambda_1^{4-l}}{(\mu -
\lambda_1)^2} B_{31} B_{11}^{l-2} B_{13}, \qquad
H_{+,l}''(-1/\lambda_1)_{33} = 2 \frac{(-1)^l \lambda_1^{4-l}}{(\mu +
\lambda_1)^2} B_{32} B_{22}^{l-2} B_{23}.
\]

Plugging both of these into~\eqref{Eleibniz2} and then
into~\eqref{Econtour}, and using that the $C_R$-term grows slower than
$m^{l-2} \lambda_1^{m-l}$ (since $l \geq 3$ now), we get:
\begin{align}\label{Eput3}
&\ \lim_{k \to \infty} \frac{1}{m_k^{l-2}
\lambda_1^{m_k-l}}[z^{m_k}]F_l(z)_{33}\notag\\
= &\ 0 - \frac{-\lambda_1^{l-2}}{2! (l-2)!} H''_{-,l}(1/\lambda_1)_{33} -
\frac{\lambda_1^{l-2} \lim_{k \to \infty} (-1)^{m_k-1}}{2! (l-2)!}
H''_{+,l}(-1/\lambda_1)_{33}\\
= &\ \left( \frac{\lambda_1}{\mu-\lambda_1} \right)^2 \frac{B_{31}
B_{11}^{l-2} B_{13}}{(l-2)!} + \left(
\frac{\lambda_1}{\mu+\lambda_1} \right)^2 \frac{B_{32} B_{22}^{l-2}
B_{23}}{(l-2)!} \bse_l({\bf m}).\notag
\end{align}

Finally, recall the $3 \times 4$ operator-matrix $T$ as in~\eqref{Eput0},
and use the asymptotics in~\eqref{Eput1}, \eqref{Eput2},
and~\eqref{Eput3} to compute:
\begin{align}\label{Efinal}
0 \leq &\ \lim_{k \to \infty} 
\frac{1}{m_k^{l-2} \lambda_1^{m_k-l}} T S_{m_k,l}(A,B) T^*\\
&\ = \begin{pmatrix}
B_{11}^l/l! & 0 & \frac{\lambda_1}{\lambda_1-\mu}
\frac{B_{11}^{l-1}}{(l-1)!} B_{13}\\
0 & \bse_l({\bf m}) B_{22}^l/l! & \bse_l({\bf m})
\frac{\lambda_1}{\lambda_1+\mu} \frac{B_{22}^{l-1}}{(l-1)!} B_{23}\\
\frac{\lambda_1}{\lambda_1-\mu} B_{31} \frac{B_{11}^{l-1}}{(l-1)!} &
\bse_l({\bf m}) \frac{\lambda_1}{\lambda_1+\mu} B_{32}
\frac{B_{22}^{l-1}}{(l-1)!} & \frac{\lambda_1^2}{(\lambda_1-\mu)^2}
B_{31} \frac{B_{11}^{l-2}}{(l-2)!} B_{13} + \bse_l({\bf m})
\frac{\lambda_1^2}{(\lambda_1+\mu)^2} B_{32} \frac{B_{22}^{l-2}}{(l-2)!}
B_{23}
\end{pmatrix}.\notag
\end{align}

Define $C \geq 0$ to be this final matrix.\medskip

\noindent \textbf{Step 9:}
We now conclude the proof of Theorem~\ref{Thilbert}(2) for $l \geq 3$,
arguing as in the $l=2$ case in Section~\ref{Squad}.
Suppose $S_{m_k,l}(A,B) \geq 0$ for infinitely many $m_k$ of the same
parity as $l$, but $[A^2,B] \neq 0$. Then the above analysis holds for
some eigenvalues $\lambda = \lambda_1 > |\mu| \geq 0$ of $A$, and we
obtain~\eqref{Efinal} with $\bse_l({\bf m}) = 1$.

Given $\varepsilon>0$, the ``perturbation'' $C_\varepsilon := C +
\varepsilon \, {\rm id}_{W_4^\perp} > 0$, whence so is its top $2 \times
2$ block, and hence its Schur complement too. Denoting this last by
$C'_\varepsilon$,
\begin{align}\label{Eschurblock2}
0 \leq \lambda_1^{-2} C'_\varepsilon = &\ \frac{\varepsilon}{\lambda_1^2}
{\rm id}_{W_3} + \frac{1}{(\lambda_1-\mu)^2} B_{31} \left(
\frac{B_{11}^{l-2}}{(l-2)!} - \frac{B_{11}^{l-1}}{(l-1)!}
(\textstyle{\frac{B_{11}^l}{l!}} + \varepsilon \, {\rm id}_{W_1})^{-1}
\displaystyle \frac{B_{11}^{l-1}}{(l-1)!} \right) B_{13}\notag\\
&\ + \frac{1}{(\lambda_1+\mu)^2} B_{32} \left(
\frac{B_{22}^{l-2}}{(l-2)!} - \frac{B_{22}^{l-1}}{(l-1)!}
(\textstyle{\frac{B_{22}^l}{l!}} + \varepsilon \, {\rm id}_{W_2})^{-1}
\displaystyle \frac{B_{22}^{l-1}}{(l-1)!} \right) B_{23}.
\end{align}

Now assume for contradiction that $B_{13} \neq 0$ (the proof is similar
if $B_{23} \neq 0$). Choose $v_\mu \in W_3 = V_\mu$ such that $v :=
B_{13} v_\mu \neq 0$, and apply Lemma~\ref{Lspectralcompute} on $W_1$,
with
\begin{equation}\label{Eapplylemma}
(\alpha, \beta, a, b, c, B, v) \quad \leadsto \quad (2l-2, l, \textstyle
\frac{1}{l!}, \frac{1}{(l-1)!^2}, \frac{1}{(l-2)!}, B_{11}, v).
\end{equation}.

Noting that $l \geq 3$, this yields:
$\tangle{v, f_\varepsilon(B_{11})v} \to -\tangle{v, B_{11}^{l-2}v} /
(l-1)! < 0$ as $\varepsilon \to 0^+$. Now akin to the $l=2$ proof, we
work in two cases:
\begin{enumerate}
\item If $B_{23} v_\mu = 0$, then the desired contradiction follows:
\[
0 \leq \lim_{\varepsilon \to 0^+} \tangle{v_\mu, \lambda_1^{-2}
C'_\varepsilon v_\mu} = \lim_{\varepsilon \to 0^+}
\frac{\varepsilon}{\lambda_1^2} \|v_\mu\|^2 - \frac{1}{(\lambda_1-\mu)^2}
\frac{\tangle{v, B_{11}^{l-2}v}}{(l-1)!},
\]
and this is negative.

\item Else $B_{23} v_\mu \neq 0$. Use Lemma~\ref{LB12} (for $B_{22} \geq
0$) to apply Lemma~\ref{Lspectralcompute} with the same scalar parameters
as in~\eqref{Eapplylemma}, and $(\mathcal{H}, B, v) \leadsto (W_2,
B_{22}, B_{23}v_\mu)$, to get (via a suitable application of
Lemma~\ref{Lclever}):
\[
\tangle{B_{23}v_\mu, f_\varepsilon(B_{22}) B_{23}v_\mu} \to -
\frac{\tangle{B_{23}v_\mu, B_{22}^{l-2} B_{23} v_\mu}}{(l-1)!} < 0,
\quad \text{as} \quad \varepsilon \to 0^+.
\]
Now applying $\tangle{v_\mu, - \cdot v_\mu}$ to both sides
of~\eqref{Eschurblock2} yields:
\[
0 \leq \lim_{\varepsilon \to 0^+} \tangle{v_\mu, \lambda_1^{-2}
C'_\varepsilon v_\mu} = \lim_{\varepsilon \to 0^+}
\frac{\varepsilon}{\lambda_1^2} \|v_\mu\|^2 - \frac{1}{(\lambda_1-\mu)^2}
\frac{\tangle{v, B_{11}^{l-2}v}}{(l-1)!},
- \frac{1}{(\lambda_1+\mu)^2} 
\frac{\tangle{B_{23}v_\mu, B_{22}^{l-2} B_{23} v_\mu}}{(l-1)!},
\]
which is also impossible.
\end{enumerate}

This completes the proof of Theorem~\ref{Thilbert}(2) for $l \geq 3$:
$A^2B = BA^2$.\medskip

\noindent \textbf{Step 10:}
Finally, the proof of Theorem~\ref{Tnxn}$(2) \implies (3)$ -- in this
case $l \geq 3$ -- is verbatim as in the final Step~7 for the $l=2$ case,
in Section~\ref{Squad}. \qed

\subsection{Comparison of asymptotics}\label{Sasymptotics}

As the reader will have noticed, the approach to characterizing the
strong BMV property for $l=1$ was explicit, while that for $l \geq 2$
relied on growth asymptotics as $m \to \infty$. It is natural to ask if
the latter can also be applied to work when $l=1$. Here we explain why
this cannot be done ``verbatim''.

Suppose $A = \begin{pmatrix} \lambda & 0 \\ 0 & \mu \end{pmatrix}$ and $B
= \begin{pmatrix} a & b \\ \overline{b} & c \end{pmatrix}$ are psd
matrices, and $m \geq l \geq 1$ are integers. Then the explicit
computation~\eqref{Ecompute2x2} yields
\[
S_{m,1}(A,B) = \begin{pmatrix} m \lambda^{m-1} a & \displaystyle
\frac{\lambda^m - \mu^m}{\lambda - \mu} b \\
\displaystyle \frac{\lambda^m - \mu^m}{\lambda - \mu} \overline{b} &
m \mu^{m-1} c
\end{pmatrix} \qquad \forall m \geq 1,
\]
which has its asymptotics as $m \to \infty$. For $l \geq 2$ there is no
explicit formula for $S_{m,l}(A,B)$, but we computed the asymptotics
earlier in this section. These asymptotics agree with the $l=1$ case for
three of the four entries, but not for the fourth -- and we have:

\begin{prop}\label{Pasymptotics}
Suppose $A = \begin{pmatrix} \lambda & 0 \\ 0 & \mu \end{pmatrix}$ and $B
= \begin{pmatrix} a & b \\ \overline{b} & c \end{pmatrix}$ are psd
matrices with $\lambda > \mu > 0$, and $m \geq l \geq 1$ are integers.
Define
\[
C_{m,l} := \frac{l!}{m^l \lambda^{m-l}} S_{m,l}(A,B) \qquad \text{and
say} \qquad
C_{m,l} = \begin{pmatrix} a_{m,l} & b_{m,l} \\ \overline{b_{m,l}} &
c_{m,l} \end{pmatrix}.
\]
\begin{enumerate}
\item Then we have the following asymptotics, at once for $l \geq 2$ and
$l=1$:
\begin{align*}
\lim_{m \to \infty} a_{m,l} = &\ a^l, \qquad 
\lim_{m \to \infty} m b_{m,l} = \frac{\lambda}{\lambda - \mu}
l b a^{l-1}, \\
\lim_{m \to \infty} m^2 c_{m,l} = &\ \begin{cases}
\left( \frac{\lambda}{\lambda - \mu} \right)^2 |b|^2 a^{l-2} l(l-1),
\qquad & \text{if } l \geq 2,\\
0, & \text{if } l = 1.
\end{cases}
\end{align*}

\item In particular, the trace and determinant of $C_{m,l}$ satisfy --
for all $l \geq 1 $ -- the asymptotics
\[
\lim_{m \to \infty} \left(\frac{\mu}{\lambda}\right)^{m-l} m^2 {\rm
Tr}(C_{m,l}) = 0, \qquad
\lim_{m \to \infty} m^2 \det C_{m,l} = - \left( \frac{\lambda}{\lambda -
\mu} \right)^2 |b|^2 a^{2l-2} l.
\]
\end{enumerate}
\end{prop}

The matrix $C_{m,l}$ and the second part will be used below, in proving a
sharp degree-bound for which elementary symmetric polynomials satisfy the
BMV conjecture/theorem (see Theorem~\ref{Ttight}).

\begin{proof}
The first part follows from~\eqref{Ecompute2x2} for $l=1$, and
from~\eqref{Egrowth11}, \eqref{Eput2}, and~\eqref{Eput3} for $l \geq 3$
(after suitably specializing) and also from above for $l=2$. Using this
first part, the second part is easily obtained, since $\mu/\lambda \in
(0,1)$.
\end{proof}


\addtocontents{toc}{\protect\vspace{1ex}}

\section{Application 2: Anderson's result on Weyl quantization}\label{Sweyl}

We now turn to applications.
We have already shown above, the first application of our main
Theorems~\ref{Thilbert} and~\ref{Tnxn}: answers to two of the
Lieb--Seiringer questions (see Question~\ref{Qlieb}), in
Theorem~\ref{Tcommute1}. Here we first show the aforementioned
strengthening of Anderson's result: Theorem~\ref{Tanderson}. This uses
some preliminaries:

\begin{prop}\label{Panderson}
Fix a real or complex Hilbert space $\mathcal{H}$, and self-adjoint
operators $A,B \in \mathcal{B}(\mathcal{H})$.
\begin{enumerate}
\item Given scalars $a,b \in \mathbb{R}$, we have:
\[
\mathcal{W}_{A,B}((x+a)^{m-l} (y+b)^l) = \mathcal{W}_{A + a \,  {\rm
id}_{\mathcal{H}}, B + b \, {\rm id}_{\mathcal{H}} }(x^{m-l} y^l), \qquad
\forall m \geq l \geq 0.
\]

\item For all smooth functions $f,g$ of one (real) variable, we have:
\[
\mathcal{W}_{A,B}(f(x) g(y)) = \mathcal{W}_{B,A}(g(x) f(y)).
\]

\item Given self-adjoint $A_3, \dots, A_d \in \mathcal{B}(\mathcal{H})$,
the Weyl calculi are compatible:
\[
\mathcal{W}_{A,B}(x^{m-l} y^l) =
\mathcal{W}_{A,B,A_3,\dots,A_d}(x_1^{m-l} x_2^l), \qquad \forall m \geq l
\geq 0.
\]
\end{enumerate}
\end{prop}

\begin{proof}
All parts follow from \cite[Theorem~2.9]{And1} using the definitions and
uniqueness of the Weyl calculus, since rapidly decaying functions in one
variable are dense in $C^\infty(\mathbb{R})$.
\end{proof}

Using this, the main result in this section follows quickly:

\begin{proof}[Proof of Theorem~\ref{Tanderson}]
It suffices to prove that $A := A_1$ and $B := A_2$ commute. Fix scalars
$a \geq \| A \|$ and $b \geq \| B \|$; then $\mathcal{W}_{A,B, A_3,
\dots, A_d}$ is supported in $[-\|A\|,\|A\|] \times [-\|B\|, \|B\|]
\times \times_{j=3}^d [-\|A_j\|, \|A_j\|]$ by~\eqref{Enumrange}, and by
hypothesis sends $(x + a)^m (y+b)$ to a positive operator for each $m
\geq 1$. Thus by Proposition~\ref{Panderson}
\begin{align*}
0 \leq \mathcal{W}_{A,B,A_3,\dots,A_d}((x + a)^m (y+b)) = &\
\mathcal{W}_{A,B}((x + a)^m (y+b)) = \mathcal{W}_{A + a\, {\rm
id}_{\mathcal{H}}, B + b\, {\rm id}_{\mathcal{H}} }(x^m y)\\
= &\ \frac{1}{m+1} S_{m+1,1}( A + a\, {\rm id}_{\mathcal{H}}, B + b\,
{\rm id}_{\mathcal{H}})
\end{align*}
for all $m \geq 1$.
As the operators $A+a\, {\rm id}_{\mathcal{H}}, B+b\, {\rm
id}_{\mathcal{H}} \geq 0$, Theorem~\ref{Thilbert} implies they commute.
That is, $0 = [A+a\, {\rm id}_{\mathcal{H}}, B+b\, {\rm
id}_{\mathcal{H}}] = [A,B]$,
as claimed.
\end{proof}

We end with an application of the Weyl calculus: an alternate proof that
the ``higher-order case'' of strong BMV -- Theorem~\ref{Thilbert}(2) --
can be extended to avoid assuming $A \geq 0$, albeit under the
``standard'', more restrictive hypotheses than in
Theorem~\ref{Thilbert}(2):

\begin{corollary}
Let $\mathcal{H}, B \geq 0, A$ be as in Theorem~\ref{Thilbert}(2), with
$\sigma(A)$ discrete (hence finite). Also fix any integer $l \geq 1$. If
$S_{m,l}(A,B) \geq 0$ for all $m \geq l$, then $AB=BA$.
\end{corollary}

\begin{proof}
We claim that $S_{m,l}(A + \|A\| \, {\rm id}_{\mathcal{H}}, B) \geq 0$
for every $m \geq l$. This immediately implies $0 = [A + \|A\| \, {\rm
id}_{\mathcal{H}}, B] = [A,B]$, by Theorem~\ref{Thilbert}.

To show the claim, denote $a := \|A\|$ and compute using
Proposition~\ref{Panderson}:
\begin{align*}
S_{m,l}(A + a\, {\rm id}_{\mathcal{H}}, B) = &\ \binom{m}{l}
\mathcal{W}_{A + a\, {\rm id}_{\mathcal{H}},B}(x^{m-l}y^l) = \binom{m}{l}
\mathcal{W}_{A,B}((x+a)^{m-l}y^l)\\
= &\ \sum_{s=0}^{m-l} \binom{m}{l} \binom{m-l}{s} a^{m-l-s}
\mathcal{W}_{A,B}(x^s y^l)\\
= &\ \sum_{s=0}^{m-l} \binom{m}{l} \binom{m-l}{s} a^{m-l-s} \cdot
\binom{l+s}{s}^{-1} S_{l+s,l}(A,B).
\end{align*}
As each summand is positive by hypothesis, we are done.
\end{proof}

\section{Application 3: Commuting via operator complete monotonicity,
via Gibbs thermal weights, and Moln\'ar--Wu}\label{Sagniva}

In this section we show the results connected to operator-valued complete
monotonicity: Theorems~\ref{Tagniva} and~\ref{Tthermal}, as well as the
Moln\'ar and Wu results in Theorem~\ref{Tmolnar}.

\begin{proof}[Proof of Theorem~\ref{Tagniva}]
For completeness, we start not from assertion~(1), but from the even more
structured statement:

\textit{There exists $\varepsilon > 0$ such that $\varphi_{\bf A}$ is
smooth on $(-\varepsilon, \varepsilon)^d \cap (\inte \Gamma)$
and~\eqref{EmatrixCM} holds.}

\noindent This statement is the ``full'' complete monotonicity as
explored by the authors in~\cite[Section~7]{Agniva}, and it in turn
implies:

\textit{$\varphi_{\bf A}$ is ``completely monotone at $0$'', two directions
at a time. More precisely:}
\[
(-1)^m (D_{\bxi}^{m-l} D_{\boldsymbol \eta}^l \varphi_{\bf
A})({\bf 0}) \geq 0, \qquad \forall m \geq l \geq 0, \ 
\bxi, {\boldsymbol \eta} \in \Gamma.
\]

In turn, this implies assertion~(1). Next, a straightforward but careful
computation shows that
\[
D_{\bxi}^{m-1} D_{\boldsymbol \eta}(\varphi_{\bf A})({\bf 0}) = c_m \cdot
(m-1)! \, S_{m,1}(\bxi \cdot {\bf A}, {\boldsymbol\eta} \cdot {\bf A}),
\qquad \forall m \geq 1, \ \bxi, \boldsymbol\eta \in \Gamma.
\]

Thus, from the hypotheses in~(1) we get that there are infinitely many
$m>0$ for which $c_m \neq 0$, not all of the same parity, and such that
$S_{m,1}(\bxi_i \cdot {\bf A}, \bxi_j \cdot {\bf A}) \geq 0$ (for all $i
< j$). Theorems~\ref{Thilbert}(1) and~\ref{Tnotpos3}(1)(a) now give that
the operators $\bxi_i \cdot {\bf A}$ (for $1 \leq i < j \leq d$) commute
pairwise. But since the $\bxi_i$ are linearly independent, taking linear
combinations first in one argument $\bxi_i \cdot {\bf A}$ and then in the
other, yields $A_i A_j = A_j A_i$ for all $1 \leq i < j \leq d$. Hence
$(1) \implies (2)$.

To show $(2) \implies (3)$, suppose $\bla \in \sigma_{joint}({\bf A})$
and $p \in \Gamma$. By definition, there exist unit vectors $x_n \in
\mathcal{H}$ such that $\| (A_j - \lambda_j {\rm id}_{\mathcal{H}})x_n \|
\to 0$ for all $1 \leq j \leq d$ (as $n \to \infty$). Therefore,
\[
\| (p \cdot {\bf A} - (p \cdot \bla){\rm id}_{\mathcal{H}})x_n \| \leq
\sum_{j=1}^d |p_j| \cdot \| (A_j - \lambda_j {\rm id}_{\mathcal{H}})x_n
\| \to 0.
\]

But this means $p \cdot \bla \in \sigma_{joint}(p \cdot {\bf A}) =
\sigma(p \cdot {\bf A})$. As ${\bf A} \in
\Gamma^*_{\mathcal{B}(\mathcal{H})}$, this implies $p \cdot \bla \geq 0$
-- for all $p \in \Gamma$, whence $\bla \in \Gamma^*$.
\end{proof}

We next prove Theorem~\ref{Tthermal}, which provides a physics flavor to
commuting of operators.

\begin{proof}[Proof of Theorem~\ref{Tthermal}]
One implication is clear: if $H_0, H_1$ commute and $H_1 \geq 0$, then
$e^{-\beta(H_0 + t H_1)} = e^{-\beta H_0} e^{-\beta t H_1}$, whence
\[
\left. -\frac{d}{dt} e^{-\beta(H_0 + t H_1)} \right|_{t=0} = \beta H_1
e^{-\beta H_0} \geq 0, \qquad \forall \beta > 0.
\]

For the other implication, we provide two proofs, in our usual two
flavors: first the simpler computation when $H_0$ has a complete
orthonormal eigenbasis for $\mathcal{H}$, and then a quite different
argument for arbitrary self-adjoint $H_0$. As usual, we split this into
steps.\medskip

\noindent \textbf{Step 1:}
Suppose $[H_0, H_1] \neq 0$, and $H_0$ has an eigenbasis as above. Then
as in above proofs -- see e.g.~\eqref{Ecalculation} -- there exist
eigenvalues $\lambda > \mu$ in $\sigma(H_0)$ such that $P_\lambda H_1
P_\mu \neq 0$. Now treat $H_0$ as a ``block diagonal operator-matrix''
(i.e., decomposing $\mathcal{H}$ into the eigenspaces of $H_0$), and
conformally partitioning $H_1$ into blocks $H'_{\eta \nu}$ with $\eta,
\nu \in \sigma(H_0)$. Now we focus on the $\{ \lambda, \mu \}$-principal
block of the derivative:
\begin{align*}
0 \leq &\ C'_\beta := (P_\lambda + P_\mu) \left( \left. -\frac{d}{dt}
e^{-\beta(H_0 + t H_1)} \right|_{t=0} \right) (P_\lambda + P_\mu)\\
= &\ (P_\lambda + P_\mu) \cdot - \sum_{m=1}^\infty \frac{(-\beta)^m}{m!}
S_{m,1}(H_0, H_1) \cdot (P_\lambda + P_\mu)\\
= &\ \sum_{m=1}^\infty \begin{pmatrix} \displaystyle \frac{(-\beta
\lambda)^{m-1}}{(m-1)!} \beta H'_{\lambda \lambda} & \displaystyle -
\frac{(-\beta \lambda)^m - (-\beta \mu)^m}{m! (\lambda-\mu)} H'_{\lambda
\mu} \\
\displaystyle - \frac{(-\beta \lambda)^m - (-\beta \mu)^m}{m!
(\lambda-\mu)} H'_{\mu \lambda} & \displaystyle \frac{(-\beta
\mu)^{m-1}}{(m-1)!} \beta H'_{\mu \mu} \end{pmatrix},
\end{align*}
using~\eqref{Ecompute2x2}. It follows that
\[
0 \leq \beta^{-1} e^{\beta \lambda} C'_\beta = \begin{pmatrix}
H'_{\lambda \lambda} &
\frac{e^{\beta(\lambda-\mu)} - 1}{\beta(\lambda-\mu)} H'_{\lambda \mu}\\
\frac{e^{\beta(\lambda-\mu)} - 1}{\beta(\lambda-\mu)} H'_{\mu \lambda} &
e^{\beta(\lambda-\mu)} H'_{\mu \mu} \end{pmatrix},
\]
for all $\beta$ in an unbounded sequence $0 < \beta_k \uparrow \infty$.

From above, $H'_{\lambda \mu} = (H'_{\mu \lambda})^* \neq 0$. Thus, there
exists $v_\mu$ in the $\mu$-eigenspace of $H_0$ such that $v :=
H'_{\lambda \mu} v_\mu \neq 0$. Now set ${\bf v}_k := (-v) \oplus
e^{\beta_k(\mu-\lambda)/2} v_\mu$, and compute:
\[
0 \leq \tangle{{\bf v}_k, \beta_k^{-1} e^{\beta_k \lambda} C'_{\beta_k}
{\bf v}_k} = \tangle{v, H'_{\lambda \lambda} v} + \tangle{v_\mu, H'_{\mu
\mu} v_\mu} - 2 \left( e^{\beta_k(\lambda-\mu)/2} -
e^{\beta_k(\mu-\lambda)/2} \right) \frac{\|v\|^2}{\beta_k(\lambda-\mu)}.
\]

As $\lambda > \mu$, the right-hand side goes to $-\infty$ as $k \to
\infty$. This contradiction concludes the proof.\medskip

\noindent \textbf{Step 2:}
Now we come to the case of general $H_0 = H_0^*$. Our proof shares
several ingredients with the proof of Theorem~\ref{Thilbert} in
Section~\ref{Shilbert}. To make it easier to follow that proof, for this
proof we set
\begin{equation}
A := H_0, \qquad B := H_1, \qquad
C_\beta^\pm := \sum_{m=1}^\infty \frac{ (\pm \beta)^m}{m!} S_{m,1}(A,B)
\end{equation}
for $\beta \in (0,\infty)$. 

We now make two ``reductions''. Akin to the discussion in
Section~\ref{Sstart}, we can take $A = H_0 \geq 0$ (but $B = H_1 = H_1^*$
is not necessarily positive). This is because
\[
e^{-\beta(H_0 + c \, {\rm id}_{\mathcal{H}} + t H_1)} = e^{-\beta c}
e^{-\beta(H_0 + t H_1)}.
\]

Next, as in Step~1 we have that
$\displaystyle \left. \frac{d}{dt} e^{\pm \beta(A + t B)} \right|_{t=0} =
C^\pm_\beta$.
Thus, differentiate (for $t$) the identity $e^{-\beta (A + t B)} \cdot
e^{\beta (A + t B)} = {\rm id}_{\mathcal{H}}$, and set $t=0$ to get:
\[
C^-_\beta e^{\beta A} + e^{-\beta A} C^+_\beta = 0 \qquad \implies \qquad
C^-_\beta = - e^{-\beta A} C^+_\beta e^{-\beta A}.
\]

This shows that $-C^-_\beta \geq 0$ if and only if $C^+_\beta \geq 0$.
Thus, we assume henceforth that $A \geq 0$ and $C^+_{\beta_k} \geq 0$ for
all $0 < \beta_k \uparrow \infty$.
We are then done by Theorem~\ref{TC2} for $f(x) = e^x$.
\end{proof}

Given where the above proof stopped, the following extends
Theorem~\ref{Tthermal} to a larger class of maps:

\begin{theorem}\label{TC2}
Let $\mathcal{H}$ be any real or complex Hilbert space, and $A \geq 0, B
= B^*$ be in $\mathcal{B}(\mathcal{H})$. Given $\beta > 0$ and a real
entire map $f : \mathbb{R} \to \mathbb{R}$, define the Fr\'echet
derivative
\[
C^+_\beta := (Df)(\beta A)[\beta B] = \lim_{t \to 0}
\textstyle\frac{f(\beta A+t \beta B)-f(\beta A)}{t}.
\]
\begin{enumerate}
\item The limit $C^+_\beta$ exists in the operator norm topology of
$\mathcal{B}(\mathcal{H})$ (whence in the SOT/WOT), for all $\beta > 0$.
\item Suppose now that $f$ further satisfies:
\begin{enumerate}
\item $f$ has nonnegative Maclaurin coefficients and is not a polynomial.
\item For all scalars $0 < c < a$ and $0 < b < 2a-c$, we have
\[
\lim_{x \to \infty} \frac{x \sqrt{f'(cx) f'(bx)}}{f(ax)} = 0.
\]
\end{enumerate}
If all $C^+_{\beta_k} \geq 0$ for some increasing unbounded sequence
$0 < \beta = \beta_k \uparrow \infty$, then $AB=BA$.
\end{enumerate}
\end{theorem}

The specific example that yields Theorem~\ref{Tthermal} is $f(x) = e^x$.
In fact, we claim that a larger family of absolutely monotonic examples
satisfy the above decay conditions: $p(x) e^x$ for any polynomial $p \neq
0$ with nonnegative real coefficients. Moreover, $\sinh(x), \cosh(x)$
also work, and since they add up to $e^x$, we may increase the above
examples to a larger class:

\begin{lemma}
Let $m,n \in \mathbb{Z}_{\geq 0}$, not both zero. Fix finitely many
constants $\gamma_i, \delta_j > 0$ and nonzero polynomials $p_i(x),
q_j(x)$ with nonnegative real coefficients. Then the function
\[
f(x) = \sum_{i=1}^n p_i(x) \sinh(\gamma_i x) + \sum_{j=1}^m q_j(x)
\cosh(\delta_j x)
\]
satisfies the conditions in Theorem~\ref{TC2}.
\end{lemma}

\begin{proof}
The key is to isolate the leading terms of $f,f'$, and then replace them
in $f,f'$ respectively.
We may suppose without loss of generality that the $\gamma_i$ are
strictly decreasing and so are the $\delta_j$. We now claim that if $p_1$
has leading term $C x^{d_1}$ with $C \in (0,\infty)$ and $d_1 \geq 0$,
then
\[
\lim_{x \to \infty} \frac{p_i(x) \sinh(\gamma_i x)}{x^{d_1} e^{\gamma_1
x}} = \lim_{x \to \infty} \frac{p_i(x) e^{\gamma_i x} - p_i(x)
e^{-\gamma_i x}}{2 x^{d_1} e^{\gamma_1 x}},
\]
and this is zero if $i>1$, since ``exponential beats polynomial''. Else
if $i=1$ then the limit is easily checked to be $C/2$.

Now widen the scope of this computation to include the $\cosh$ terms too:
set
\[
\eta = \max(\gamma_1, \delta_1) > 0, \qquad d := \begin{cases}
\deg p_1, & \text{if } \gamma_1 > \delta_1,\\
\deg q_1, & \text{if } \gamma_1 < \delta_1,\\
\max(\deg p_1, \deg q_1), \quad & \text{if } \gamma_1 = \delta_1.
\end{cases}
\]

Examining term by term akin to above, it follows that
$\displaystyle \lim_{x \to \infty} \frac{f(x)}{e^{\eta x} x^d}$
is a positive constant, say $C/2$ for some $C$. Now a similar analysis
using the Leibniz rule shows that
\[
\lim_{x \to \infty} \frac{f'(x)}{e^{\eta x} x^d} = C \eta / 2.
\]

Replacing these in the desired asymptotics and canceling $C/2$,
\[
\lim_{x \to \infty} \frac{x \sqrt{f'(cx) f'(bx)}}{f(ax)} = \lim_{x \to
\infty} \frac{x \sqrt{\eta e^{\eta c x} (cx)^d \cdot \eta e^{\eta b
x} (bx)^d}}{e^{\eta ax} (ax)^d} = \lim_{x \to \infty} \frac{\eta
(bc)^{d/2} x}{a^d \exp(\eta (a - (b+c)/2)) x},
\]
and this vanishes because $2a > b+c$ by hypothesis.
\end{proof}

\begin{proof}[Proof of Theorem~\ref{TC2}]
In fact a stronger statement than the first part holds, which we will
require. Notice for every monomial $x^m$ that the Fr\'echet/Gateaux
derivative at $A$ along $B$ of $x^m$ exists, being precisely the Hurwitz
product:
\[
D(x^m)(A)[B] = \lim_{t \to 0} \frac{(A+tB)^m - A^m}{t} = S_{m,1}(A,B),
\qquad \forall m \geq 0.
\]
Hence by linearity, the limit $Dg(A)[B]$ exists in norm for every
polynomial $g$. Now a standard fact in functional analysis -- see e.g.\
\cite[(9.12.1)]{Dieu} -- says that if $T_m(f)$ is the $m$th Taylor
polynomial of $f$, then ($T_m(f) \to f$ uniformly on every compact subset
of $\mathbb{R}$, and) one has norm-convergence of the Fr\'echet
derivatives:
\[
Df(A)[B] = \lim_{m \to \infty} D(T_m(f))(A)[B].
\]
This shows part~(1). In particular, we write $f$ as a power series and
deduce $C^+_\beta$ from this:
\begin{equation}\label{ECplusbeta}
f(x) = \sum_{k=0}^\infty c_k x^k, \ c_k \in [0,\infty)\ \forall k \qquad
\implies \qquad C^+_\beta = \sum_{k=1}^\infty c_k S_{k,1}(\beta A, \beta
B).
\end{equation}

The real proof is that of part~(2). We proceed in steps, akin to the
proof of Theorem~\ref{Thilbert}(1).\medskip

\noindent \textbf{Step 1:}
Fix $\beta = \beta_k > 0$. We first replicate the analogue
of~\eqref{Ecommutator} in this setting. For this, begin with the
algebraic identity
\[
[A^m,B] = [A,S_{m,1}(A,B)] = [A, D(x^m)(A)[B]], \qquad \forall m \geq 1.
\]

Taking linear combinations, this holds for all polynomials, and in
particular:
\[
[ T_m(f)(A), B] = [A, D(T_m(f))(A)[B]], \qquad \forall m \geq 1.
\]

Take $m \to \infty$ on both sides. The right-side converges in norm to
$[A, Df(A)[B]] = [A, C^+_1]$, while the left-side clearly converges
in norm to $[f(A),B]$, since $\| f(A) - T_m(f)(A) \| \leq
\sum_{k=m+1}^\infty c_k \|A\|^k$, which goes to $0$ as $m \to \infty$.
Replacing $(A,B) \leadsto (\beta A, \beta B)$, we thus get:
\begin{equation}\label{ECbeta00}
[f(\beta A),\beta B] = [\beta A, C^+_\beta].
\end{equation}

\noindent \textbf{Step 2:}
Suppose $AB \neq BA$. Since $A \geq 0$, define $F^\sigma := F \cap
\sigma(A)$ for $F \subseteq \mathbb{R}$. Now by Lemma~\ref{Lstasi}
applied to $A' = A$ and $\lambda_0 = 0$, there exists $v_\mu \in {\rm
Ran}(E([0,c]^\sigma))$ such that $v := E(F^\sigma) B v_\mu \neq 0$.
Here if $c=0$, we may increase $c$ slightly to still satisfy $0 < c < a$
and $b+c < 2a$ (this is used at the end). Now as in the proof of
Theorem~\ref{Thilbert}(1), from~\eqref{ECbeta00} we get the desired
counterpart of~\eqref{Ecommutator}:
\begin{equation}\label{ECbeta}
|\tangle{v, f(\beta A) \beta B v_\mu}| \leq
|\tangle{v, \beta B f(\beta A) v_\mu}| +
|\tangle{v, \beta A C^+_\beta v_\mu}| +
|\tangle{v, C^+_\beta \beta A v_\mu}|.
\end{equation}

\noindent \textbf{Step 3:}
We now bound these four terms.
The left-hand side is lower-bounded as in~\eqref{Ebound1} in the proof of
Theorem~\ref{Thilbert}(1):
\[
\tangle{v, f(\beta A) \beta B v_\mu} = \beta \int_{F^\sigma} f(\beta
\lambda) d E_{v,v}(\lambda) \geq \beta f(\beta a) \|v\|^2 > 0.
\]

The first term on the right in~\eqref{ECbeta} is bounded above by
Cauchy--Schwarz and Corollary~\ref{Cborel}(1):
\[
|\tangle{v, \beta B f(\beta A) v_\mu}| \leq \|v\| \, \|B\| \, \|v_\mu\|
\, \beta \, f(\beta c).
\]

The other two terms are bounded above using that $C^+_\beta \geq 0$ for
all $\beta = \beta_k \uparrow \infty$, so that it induces a semi-inner
product on $\mathcal{H}$. Hence by Cauchy--Schwarz,
\[
|\tangle{v, A C^+_\beta v_\mu}| \leq
|\tangle{A v, C^+_\beta A v}|^{1/2}
|\tangle{v_\mu , C^+_\beta v_\mu}|^{1/2}.
\]

The first term on the right is upper-bounded by Corollary~\ref{Cborel}
and~\eqref{ECplusbeta}:
\[
|\tangle{A v, C^+_\beta A v}| \leq \sum_{k=1}^\infty c_k |\tangle{Av,
S_{k,1}(\beta A, \beta B) Av}| \leq \sum_{k=1}^\infty c_k \, \|B\| \,
\|v\|^2 \, k \beta^k b^{k+1} = \|B\| \, b^2 \, \|v\|^2 \, \beta f'(b
\beta).
\]

We similarly bound the second term on the right:
\[
|\tangle{v_\mu, C^+_\beta v_\mu}| \leq \|B\| \, \|v_\mu\|^2 \, \beta f'(c
\beta),
\]
which yields:
\[
|\tangle{v, A C^+_\beta v_\mu}| \leq \|B\| \, \|v\| \, \|v_\mu\| \, b
\cdot \beta \sqrt{f'(c \beta) f'(b \beta)}.
\]

In a similar fashion, the final term on the right-side of~\eqref{ECbeta}
is bounded above:
\[
|\tangle{v, C^+_\beta A v_\mu}| \leq \|B\| \, \|v\| \, \|v_\mu\| \, c
\cdot \beta \sqrt{f'(c \beta) f'(b \beta)}.
\]

\noindent \textbf{Step 4:}
Collecting together these bounds for the four terms in~\eqref{ECbeta},
and dividing by $\beta$,
\[
f(a \beta) \|v\|^2 \leq \|B\| \, \|v\| \, \|v_\mu\| \left( f(c \beta) +
(b+c) \beta \sqrt{f'(c \beta) f'(b \beta)} \right).
\]

Divide by $f(a \beta)$ and let $\beta = \beta_k \uparrow \infty$. Noting
that $2a > b+c$ by Lemma~\ref{Lstasi}, the given asymptotics on $f$ and
the next lemma imply $\|v\|^2 \leq 0$, which is a contradiction. Hence
$AB=BA$.
\end{proof}

\begin{lemma}\label{Ldecay}
Suppose $f(x)$ is an entire function with nonnegative Maclaurin
coefficients, and which is not a polynomial. Then for all $0 < c < a$, we
have $\lim_{x \to \infty} f(cx)/f(ax) = 0$.
\end{lemma}

\noindent As this is standard, we defer the proof to
Appendix~\ref{Appendix3}.

We close this section with

\begin{proof}[Proof of Theorem~\ref{Tmolnar}]
Here we only show the first two parts, by Moln\'ar and Wu. (The third and
fourth parts are proved in Section~\ref{Suchi}.)
\begin{enumerate}
\item The easier implication is one line, as shown by Moln\'ar: if $H_0$
is central and $\beta>0$, then
\[
e^{\beta H_0} = e^{\beta H_0/2} \cdot 1 \cdot e^{\beta H_0/2} \leq
e^{\beta H_0/2} e^{\beta (H-H_0)} e^{\beta H_0/2} = e^{\beta H_0/2 +
\beta (H-H_0) + \beta H_0/2} = e^{\beta H}.
\]

The converse is what we now show, avoiding the machinery used by
Moln\'ar: (irreducible) representations of $\mathcal{A}$ and the Kadison
transitivity theorem. The proof is in three steps:
(a)~$H_0$ commutes with every positive (self-adjoint) element,
(b)~hence with every self-adjoint element,
(c)~and hence with every element.

We begin with~(a). By Gelfand--Naimark--Segal (GNS), let $\mathcal{A}
\hookrightarrow \mathcal{B}(\mathcal{H})$ for some complex Hilbert space
$\mathcal{H}$.
Let $0 \leq H_1 = H_1^*\in \mathcal{A}$ and $\beta = \beta_k$ for some
$k$. Then $\beta H_0 \leq \beta (H_0 + t H_1)$ for all $\beta, t > 0$, so
$e^{\beta H_0} \leq e^{\beta (H_0 + t H_1)}$. Since inverting reverses
the Loewner order, we get
\[
e^{-\beta(H_0 + t H_1)} \leq e^{-\beta H_0}, \qquad \forall t > 0.
\]
Subtract and divide by $t > 0$ to get for all $k \geq 1$:
\[
\frac{e^{-\beta_k(H_0 + t H_1)} - e^{-\beta_k H_0}}{t} \leq 0, \qquad
\forall t>0.
\]
Taking $t \to 0^+$ yields precisely our hypothesis in
Theorem~\ref{Tthermal}, and hence $H_0$ commutes with $H_1$, as claimed.

The rest is standard. Using the continuous functional calculus, we see
that (b)~every self-adjoint element $H$ is a difference of two positive
elements: $2H = (|H|+H) - (|H|-H)$, so $H_0$ commutes with every $H =
H^*$. Finally, (c)~every $H \in \mathcal{A}$ is a linear combination of
two self-adjoint elements: $H = \frac{H+H^*}{2} + i \frac{H-H^*}{2i}$. So
$H_0$ commutes with all $H \in \mathcal{A}$.

\item By the first part, $H_0$ is central for every $H_0 = H_0^*$. Hence
by (c)~in the preceding step, every element is central. The converse is
as in the previous part. \qedhere
\end{enumerate}
\end{proof}

\section{Strengthening the BMV conjecture: Jack/Macdonald polynomials of
the Hurwitz-product spectrum}\label{SeventualBMV2}

Here we prove the ``unconditional'' Theorems~\ref{Tej2} and~\ref{Ttight},
which refine Theorem~\ref{Tej}. To set the tone, we first write the BMV
conjecture (now a theorem, given~\cite{Stahl,Heina}) as a series of
increasingly simpler statements:

\begin{theorem}[BMV theorem]\label{TBMViff2}
Given an integer $n \geq 1$ and a Hermitian matrix $B_{n \times n}$, the
following are equivalent:
\begin{enumerate}
\item (BMV.) For all psd matrices $A_{n \times n}$ and integers $m \geq l
\geq 1$, we have ${\rm Tr}(S_{m,l}(A,B)) \geq 0$.

\item (Eventual/asymptotic BMV.) 
For all positive definite matrices $A_{n \times n}$ and integers $l
\geq 1$, there exists a threshold $m_l = m_l(A,B)$ such that we have
${\rm Tr}(S_{m,l}(A,B)) \geq 0$ for all $m \geq m_l$.

\item (Via commuting matrices.)
For all positive definite matrices $A_{n \times n}$ with $AB=BA$, and for
any odd $l \geq 1$, we have ${\rm Tr}(S_{m,l}(A,B)) \geq 0$ for
infinitely many $m \geq l$.

\item $B$ is positive semidefinite.
\end{enumerate}
\end{theorem}

\begin{proof}
That $(4) \implies (1)$ is the BMV theorem \cite{Stahl,Heina}, and that
$(1) \implies (2) \implies (3)$ is immediate. Now write $B = U^* D U$
with $U$ unitary and $D$ diagonal. Fix $1 \leq j \leq n$ and let $A = U^*
\cdot \frac{1}{2} ({\rm Id}_n + E_{jj}) \cdot U$, where $E_{jj}$ is an
elementary (diagonal) matrix. By hypothesis~(3), there exists a sequence
$l \leq m_1 < m_2 < \cdots$ such that ${\rm Tr}(S_{m_k,l}(A,B)) \geq 0$
for all $k \geq 1$. But this trace can be explicitly computed, which
yields:
\begin{align*}
0 \leq \lim_{k \to \infty} \frac{1}{\binom{m_k}{l}} {\rm
Tr}(S_{m_k,l}(A,B)) = &\ \lim_{k \to \infty} \frac{1}{\binom{m_k}{l}}
{\rm Tr}(S_{m_k,l}(U A U^*, U B U^*))\\
= &\ \lim_{k \to \infty} \left( d_{jj}^l + \frac{1}{2^{m_k-l}}\sum_{i
\neq j} d_{ii}^l \right) = d_{jj}^l.
\end{align*}
As $l$ is odd, we get $d_{jj} \geq 0$ for all $j$, and so $B = U^* D U$
is psd.
\end{proof}

As mentioned in the opening section, our methods help prove the main
result in~\cite{FF}: the eventual-BMV property $(4) \implies (2)$ -- and
not only for all large $m \gg 0$ but also for all small $m \ll 0$. We now
prove the stronger result mentioned in the Introduction:
Theorem~\ref{Tej2}. Even more strongly, we show:

\begin{theorem}\label{Tstronger}
Fix integers $n,l \geq 1$ and a positive definite matrix $A_{n \times
n}$. Let $\mathbb{H}_n$ denote the Hermitian matrices in $\mathbb{C}^{n
\times n}$, and let $f : \mathbb{H}_n \to \mathbb{R}$ be a continuous map
satisfying:
\begin{enumerate}
\item $f(U^* M U) = f(M)$ for all Hermitian $M$ and unitary $U$.
\item If $0 \neq M$ is psd of rank at least $d_A$ (see
Theorem~\ref{Tej2}), then $f(M) > 0$.
\item There exists a continuous map $g_0 : (0,\infty) \to (0,\infty)$
such that $f(cM) = g_0(c)f(M)$ for all scalars $c \in (0,\infty)$ and
Hermitian $M$.
\end{enumerate}
Then for any positive definite $B_{n \times n}$ and real exponent
$\alpha$, there exists a threshold $m_l > 0$ such that $f(A^\alpha
S_{m,l}(A,B) A^\alpha) > 0$ for $m \geq m_l$.
\end{theorem}

\begin{proof}
First suppose $m>0$, and let $A$ have maximum eigenvalue
$\lambda_1 = \lambda_A$ as in the previous sections. Since
$f$ is unitarily invariant, we may assume without loss of generality that
$A$ is a block diagonal matrix $\lambda_1 {\rm Id}_{n_1} \oplus A'$, with
$\lambda_1 > 0$ the maximum eigenvalue of $A$. Let $P_1 :=
\begin{pmatrix} {\rm Id}_{n_1} & 0 \\ 0 & 0 \end{pmatrix}$ be the
projection onto the $\lambda_1$-eigenspace of $A$, and
recall~\eqref{Egrowth11}, for the sequence $m_k = k$ if $l\geq3$, as well
as the asymptotics in~\eqref{Ecompute2x2} for $l=1$ and in Step~5 of the
proof of Theorem~\ref{Tnxn} for $l=2$. Pre- and post- multiplying these
by $A^\alpha$ yields:
\[
\lim_{m \to \infty} \frac{A^\alpha S_{m,l}(A,B) A^\alpha}{m^l
\lambda_1^{m-l}} = \frac{1}{l!} A^\alpha (P_1 B P_1)^l A^\alpha.
\]
But $A$ commutes with its spectral projection $P_1 = P_1^2$, whence so does
$A^\alpha$ for any $\alpha$, and so
\[
A^\alpha P_1 = A^\alpha P_1^2 = P_1 A^\alpha P_1 = \lambda_1^\alpha P_1.
\]

This implies:
\[
\lim_{m \to \infty} \frac{A^\alpha S_{m,l}(A,B) A^\alpha}{m^l
\lambda_1^{m-l}} = \frac{\lambda_1^{2\alpha}}{l!} (P_1 B P_1)^l.
\]
As $B$ is positive definite, $P_1 B P_1$ is a positive definite operator
on the $\lambda_1$-eigenspace of $A$, hence of rank $d_A$. Thus, the same
holds for the right-hand side, and so by the hypotheses,
\[
0 < f \left( \frac{\lambda_1^{2 \alpha}}{l!} (P_1 B P_1)^l \right) =
\lim_{m \to \infty} g_0(m^{-l} \lambda_1^{l-m}) f(A^\alpha S_{m,l}(A,B)
A^\alpha),
\]
and the positivity of $g_0$ now concludes the proof.
\end{proof}

From this, we deduce

\begin{proof}[Proof of Theorem~\ref{Tej2}]\hfill
\begin{enumerate}
\item For the first assertion in here, write the entries of $\lambda$ as
$\lambda_1 \geq \cdots \geq \lambda_{d'} > \lambda_{d'+1} = 0$, with $d'
\leq d_A$. Also set $f(M) = m_\lambda(\sigma(M))$; then $f$ satisfies
Theorem~\ref{Tstronger}(1) as the spectrum is unitarily invariant. Next,
if $0 \leq M$ is psd of rank $\geq d_A > 0$, then write $M = U^* D U$
with diagonal $D$ whose nonnegative diagonal entries $(a_{jj})$ are in
non-increasing order. Now by assumption~(1),
\[
m_\lambda(\sigma(M)) = m_\lambda(\sigma(D)) \geq \prod_{j=1}^{d'}
a_{jj}^{\lambda_j} > 0.
\]

This shows Theorem~\ref{Tstronger}(2), and the next assumption is obvious
by inspection with $g(c) := c^{|\lambda|}$, where $|\lambda| := \lambda_1
+ \lambda_2 + \cdots$ counts the total number of cells in $\lambda$.
Thus the hypotheses of Theorem~\ref{Tstronger} hold, and the first
assertion in Theorem~\ref{Tej2} follows from Theorem~\ref{Tstronger} upon
setting $\alpha=0$.

We next consider the second assertion -- i.e., with $\sigma( (-1)^l
S_{-m,l}(A,B))$ for $A,B$ positive definite and $m \gg 0$. This
essentially follows from (the above proof for) the first assertion, using
Proposition~\ref{Psml} and the specializations
\[
(A, m, g(c), B, A^\alpha) \ \leadsto \ (A^{-1}, m-1+l \geq m,
c^{|\lambda|}, B' = A^{-1/2} B A^{-1/2}, A^{-1/2}).
\]

\item If $AB=BA$, $m_+ \geq l$, and $m_- \geq 1$, then both
\begin{align*}
S_{m_+,l}(A,B) = &\ \binom{m_+}{l} B^{l/2} A^{m_+-l} B^{l/2} \qquad
\text{and}\\
(-1)^l S_{-m_-,l}(A,B) = &\ \binom{m_-+l-1}{l} B^{l/2} A^{-(m_-+l)}
B^{l/2}
\end{align*}
are positive definite, so every monomial function $m_\lambda(\cdot)$ in
their spectra is positive.

\item This follows immediately from the first part, since any
monomial-positive symmetric function is a finite $\mathbb{R}_+$-linear
combination of monomial symmetric functions $m_\lambda$.

\item This does not follow directly from the preceding part, since
Schur/Jack/Macdonald polynomials are linear combinations of monomials
$m_\mu$ with $\mu$ majorized by $\lambda$ -- so $\mu$ can have more than
$d_A$ parts even if $\lambda$ does not. Thus, an argument is needed.

The key is to again invoke Theorem~\ref{Tstronger}. Thus, fix $A,B > 0$
and let $f(M) = m_\lambda(\sigma(M)) + \sum_{\mu \prec \lambda} c_\mu
m_\mu(\sigma(M))$ be any symmetric function with all $c_\mu \in
[0,\infty)$. (This subsumes the Jack/Macdonald polynomials in the
assertion, and hence all of their specializations, as desired.) We claim
that $f$ satisfies the hypotheses of Theorem~\ref{Tstronger}, which would
prove this assertion (with $\alpha=0$).

Indeed, $f$ is unitarily invariant as in~(1), and the leading term
$m_\lambda(\cdot)$ helps verify~(2) since for all other partitions $\mu
\prec \lambda$, $m_\mu(\sigma(M)) \geq 0$. Finally, $\mu \prec \lambda$
implies $|\mu| = |\lambda|$, and so~(3) holds with $g_0(c) =
c^{|\lambda|}$. This concludes the proof. \qedhere
\end{enumerate}
\end{proof}

It remains to show that Theorem~\ref{Tej2}(1) is tight:

\begin{proof}[Proof of Theorem~\ref{Ttight}]
For both parts, we first denote $\lambda := \lambda_{\max}(A), d := d_A <
n$, and use the Spectral theorem to write
\[
U A U^* = \lambda {\rm Id}_d \oplus {\rm diag}(\mu_1, \dots,
\mu_{n-d}),
\]
where $\lambda = \lambda_{\max}(A) > \mu_1 \geq \cdots \geq \mu_{n-d} >
0$. Define $D := {\rm diag}(\mu_2, \dots, \mu_{n-d})$, and denote
$e_k(\sigma(M)) =: e_k(M)$ for convenience. (As usual, $e_0(\cdot) :=
1$.)
\begin{enumerate}
\item 
For $b > 0,\epsilon \in [0,1)$ we define
\[
B_{b,\epsilon} := U^* \begin{pmatrix} {\rm Id}_{d-1} & 0 & 0 \\ 0 & b
{\bf 1}_{2 \times 2} + \epsilon {\rm Id}_2 & 0 \\ 0 & 0 & {\rm
Id}_{n-d-1} \end{pmatrix} U.
\]
Then
\[
U S_{m,1}(A,B_{b,\epsilon}) U^* = \begin{pmatrix} m \lambda^{m-1} {\rm
Id}_{d-1} & 0 & 0 & 0 \\ 0 & m \lambda^{m-1} (b + \epsilon) &
\frac{\lambda^m - \mu_1^m}{\lambda - \mu_1} \, b & 0\\
0 & \frac{\lambda^m - \mu_1^m}{\lambda - \mu_1} \, b & m \mu_1^{m-1} (b +
\epsilon) & 0 \\
0 & 0 & 0 & m D^{m-1} \end{pmatrix}.
\]
For convenience, we set $\widetilde{D} := \begin{pmatrix} m \lambda^{m-1}
{\rm Id}_{d-1} & 0 \\ 0 & m D^{m-1} \end{pmatrix}$. 

Now $e_k(M)$ is the sum of the principal $k \times k$ minors of any
square matrix $M$. Each such minor includes zero, one, or two rows (and
columns) from the $d$th and $(d+1)$th ones. Using this, for all $2 \leq k
\leq n$ we compute:
\begin{align*}
&\ e_k(S_{m,1}(A,B_{b,\epsilon}))\\
= &\ e_k(U S_{m,1}(A,B_{b,\epsilon}) U^*)\\
= &\ e_k(\widetilde{D}) + (b + \epsilon) (m \lambda^{m-1} + m
\mu_1^{m-1}) e_{k-1}(\widetilde{D})\\
&\ + e_{k-2}(\widetilde{D}) \left( \epsilon^2 \cdot m^2 \lambda^{m-1}
\mu_1^{m-1} + 2 m^2 b \epsilon (\lambda \mu_1)^{m-1} + b^2
t'_m \right),
\end{align*}
where we define and claim
\begin{equation}\label{Eamgm}
t'_m := m^2 (\lambda \mu_1)^{m-1} - \left( \frac{\lambda^m -
\mu_1^m}{\lambda - \mu_1} \right)^2 < 0, \quad \forall m \geq 2.
\end{equation}

We show the claim later; for now we complete the proof. As each
$e_{k-2}(\widetilde{D}) > 0$, choose sufficiently large $b_0 \gg 0$ such
that
\[
e_k(S_{m,1}(A,B_{b_0,0})) < 0, \qquad \forall 2 \leq k \leq n, \
2 \leq m \leq m_0.
\]
Hence by perturbation, there exists small $0 < \epsilon_0 \ll 1$ such
that
\[
e_k(S_{m,1}(A,B_{b_0,\epsilon_0})) < 0, \qquad \forall 2 \leq k \leq n, \
2 \leq m \leq m_0
\]
and we are done as $B_{b_0,\epsilon_0}$ is positive definite.

It remains to show the claim~\eqref{Eamgm}. This is clear via ``pairing
terms'', since $\lambda > \mu_1 > 0$:
\begin{align*}
-t'_m = &\ (\lambda^{m-1} + \cdots + \mu_1^{m-1})^2 - m^2 \lambda^{m-1}
\mu_1^{m-1}\\
= &\ \sum_{j=0}^{m-1} k_j (\lambda^{m-1+j} \mu_1^{m-1-j} +
\lambda^{m-1-j} \mu_1^{m-1+j} - 2 \lambda^{m-1} \mu_1^{m-1})\\
= &\ \sum_{j=0}^{m-1} k_j (\lambda\mu_1)^{m-1-j} (\lambda^j - \mu_1^j)^2 > 0,
\end{align*}
where $k_0, \dots, k_{m-1}$ are certain positive integers that add up to
$m^2$.
\medskip

\item For convenience, we split the proof into steps.\medskip

\noindent \textbf{Step 1:}
With $U A U^*$ as above, apply a further permutation matrix $\eta$ such
that
\[
\eta (U A U^*) \eta^{-1} = \begin{pmatrix} \lambda & 0 \\ 0 & \mu_1
\end{pmatrix} \oplus \lambda {\rm Id}_{d-1} \oplus D.
\]

Now set $U_1 := \eta U$, which is again unitary, and also choose any psd
matrices
\[
\begin{pmatrix} a & b \\ \overline{b} & c \end{pmatrix}_{2 \times 2},
\quad B_1 \in \mathbb{C}^{(d-1) \times (d-1)}, \quad
D_2 \in \mathbb{C}^{(n-d-1) \times (n-d-1)},
\]
with $b \neq 0$, $ac > |b|^2$, $D_2$ diagonal, and $B_1, D_2$
nonsingular. Now define
\[
B' := \begin{pmatrix} a & b \\ \overline{b} & c \end{pmatrix}_{2 \times 2}
\oplus B_1 \oplus D_2 \geq 0,
\]
and $B := U_1^* B' U_1 \geq 0$. Then $S_{m,l}(A,B)$ is unitarily
conjugate to the block diagonal matrix
\begin{equation}\label{Esigma}
S_{m,l}(U_1 A U_1^*, B') = S_{m,l} \left( \begin{pmatrix} \lambda & 0 \\
0 & \mu_1 \end{pmatrix}, \begin{pmatrix} a & b \\ \overline{b} &
c\end{pmatrix} \right) \oplus \binom{m}{l} \lambda^{m-l} B_1^l \oplus
\binom{m}{l} D_2^l D^{m-l},
\end{equation}
and hence their $e_k \circ \sigma$-values coincide for all $k$.\medskip

\noindent \textbf{Step 2:}
Now it suffices to prove the following asymptotic:
\begin{equation}
\lim_{m \to \infty} \left(\frac{\lambda^{k-d-1}}{\mu_2 \cdots \mu_{k-d}}
\right)^{m-l} \, m^2 \, e_k \left( \textstyle{\frac{l!}{m^l
\lambda^{m-l}} } S_{m,l}(U_1 A U_1^*, B') \right) < 0,
\qquad \forall k \in (d, n].
\end{equation}

In fact, we compute the above limit. For this, 
let $D_2 = {\rm diag}(p_1, \dots, p_{n-d-1})$ and let $r \geq 1$; then
\begin{equation}\label{Eer}
e_r(D_2^l (\lambda^{-1} D)^{m-l}) = e_r( p_1^l (\lambda^{-1}
\mu_2)^{m-l}, \dots, p_{n-d-1}^l (\lambda^{-1} \mu_{n-d})^{m-l} )
\end{equation}

Since $\lambda > \mu_2 \geq \cdots \geq \mu_{n-d} > 0$, we see that this
can be compared with its slowest-decaying term, which is $\prod_{j=1}^r
(\lambda^{-1} \mu_{j+1})^{m-l}$. Thus, every monomial in~\eqref{Eer},
upon dividing by $\prod_{j=1}^r (\lambda^{-1} \mu_{j+1})^{m-l}$, has a
finite nonnegative limit as $m \to \infty$. Hence the following limit
exists:
\begin{equation}\label{Egamma}
\lim_{m \to \infty} \frac{e_r(D_2^l (\lambda^{-1}
D)^{m-l})}{\prod_{j=1}^r (\lambda^{-1} \mu_{j+1})^{m-l}} =:
\lim_{m \to \infty} \frac{e_r(D_2^l D^{m-l})}{\prod_{j=1}^r
\mu_{j+1}^{m-l}} =:
\gamma_{r,l} \in (0,\infty), \qquad 0 \leq r \leq n-d-1.
\end{equation}

Also set $\gamma_{0,l} := 1$.
We now claim that for all $l \geq 1$ and $d < k \leq n$,
\begin{align}\label{Eets}
\begin{aligned}
&\ \lim_{m \to \infty} \frac{\lambda^{(m-l)(k-d-1)}}{(\mu_2 \cdots
\mu_{k-d})^{m-l}} m^2 \, e_k \left( \textstyle{\frac{l!}{m^l
\lambda^{m-l}} } S_{m,l}(U_1 A U_1^*, B') \right)\\
= &\ -\gamma_{k-d-1,l} l \left( \frac{a^{l-1} |b| \lambda}{\lambda -
\mu_1} \right)^2 \det(B_1)^l < 0.
\end{aligned}
\end{align}
This would prove the result.\medskip

\noindent \textbf{Step 3:}
We begin by computing the argument of $e_k$ in~\eqref{Eets}. Define $p_m
:= 1 (1 - \frac{1}{m}) \cdots (1 - \frac{l-1}{m}) = \frac{l!}{m^l}
\binom{m}{l}$ (so $p_m = 1$ if $l=1$). Then,
\[
\frac{l!}{m^l \lambda^{m-l}} S_{m,l}(U_1 A U_1^*, B')
= C_{m,l} \oplus p_m B_1^l \oplus p_m D_2^l (\lambda^{-1} D)^{m-l},
\]
where $C_{m,l}$ was defined and its asymptotics computed in
Proposition~\ref{Pasymptotics}. We use the same strategy as in the
previous part to compute $e_k$: we take zero, one, or two rows/columns
from the first two (i.e.\ from $C_{m,l}$), then take $i$ out of the $d-1$
rows/columns from the block $p_m B_1^l$, and the rest from $p_m D_2^l
(\lambda^{-1}D)^{m-l}$. Summing over all such $k \times k$ minors yields:
\begin{align}\label{Ealigned}
e_k \left( \textstyle{\frac{l!}{m^l \lambda^{m-l}} } S_{m,l}(U_1 A
U_1^*, B') \right)
= &\ p_m^k \sum_{i=0}^{d-1} e_i(B_1^l) e_{k-i}(D_2^l (\lambda^{-1}
D)^{m-l})\\
&\ + p_m^{k-1} \sum_{i=0}^{d-1} e_i(B_1^l) e_{k-1-i}(D_2^l (\lambda^{-1}
D)^{m-l}) \cdot {\rm Tr} (C_{m,l}) \notag\\
&\ + p_m^{k-2} \sum_{i=0}^{d-1} e_i(B_1^l) e_{k-2-i}(D_2^l (\lambda^{-1}
D)^{m-l}) \cdot \det C_{m,l}. \notag
\end{align}

\noindent \textbf{Step 4:}
We now multiply~\eqref{Ealigned} by $\displaystyle
\frac{\lambda^{(m-l)(k-d-1)}}{(\mu_2 \cdots \mu_{k-d})^{m-l}} m^2$, then
take $m \to \infty$, and use~\eqref{Egamma} to examine what happens to
each summand. This involves working with $e_{q-i}(D_2^l (\lambda^{-1}
D)^{m-l})$ for $q = k, k-1, k-2$. Note also that $p_m \to 1$ for all $l
\geq 1$; and $e_i(B_1^l)$ is independent of $m$.

First, for each term in the first sum -- so $q=k$ and $i \leq d-1$ -- we
have:
\begin{align*}
&\ \lim_{m \to \infty} \left(\frac{\lambda^{k-d-1}}{\mu_2 \cdots
\mu_{k-d}} \right)^{m-l} \, m^2 \cdot p_m^k e_i(B_1^l) e_{k-i}(D_2^l
(\lambda^{-1} D)^{m-l})\\
= &\ \lim_{m \to \infty} \left(\frac{\lambda^{k-d-1}}{\mu_2 \cdots
\mu_{k-d}} \right)^{m-l} \, m^2 \cdot p_m^k e_i(B_1^l) \cdot
\gamma_{k-i,l} \prod_{j=1}^{k-i} (\lambda^{-1} \mu_{j+1})^{m-l}\\
= &\ \lim_{m \to \infty} \left(\frac{\mu_{k-d+1} \cdots
\mu_{k-i+1}}{\lambda^{d-i+1}} \right)^{m-l} \, m^2 \cdot p_m^k e_i(B_1^l)
\gamma_{k-i,l},
\end{align*}
and since ``exponential beats polynomial'', this limit is zero.

Similarly, for each term in the second sum -- so $q=k-1$ and $i \leq d-1$
-- we have:
\begin{align*}
&\ \lim_{m \to \infty} \left(\frac{\lambda^{k-d-1}}{\mu_2 \cdots
\mu_{k-d}} \right)^{m-l} \, m^2 \cdot p_m^{k-1} e_i(B_1^l)
e_{k-1-i}(D_2^l (\lambda^{-1} D)^{m-l}) \cdot {\rm Tr} (C_{m,l})\\
= &\ \lim_{m \to \infty} \left(\frac{\lambda^{k-d-1}}{\mu_2 \cdots
\mu_{k-d}} \right)^{m-l} \, m^2 \cdot p_m^{k-1} e_i(B_1^l) \cdot
\gamma_{k-1-i,l} \prod_{j=1}^{k-1-i} (\lambda^{-1} \mu_{j+1})^{m-l} \cdot
{\rm Tr} (C_{m,l})\\
= &\ \lim_{m \to \infty} \left(\frac{\mu_{k-d+1} \cdots
\mu_{k-i}}{\lambda^{d-i}} \right)^{m-l} \, m^2 \cdot p_m^{k-1}
e_i(B_1^l) \gamma_{k-1-i,l} \cdot {\rm Tr} (C_{m,l}),
\end{align*}
and this limit is again zero, by Proposition~\ref{Pasymptotics}(2).

We come to the terms in the third sum -- so $q = k-2$ and $i \leq d-1$.
Now we again use Proposition~\ref{Pasymptotics}(2) to compute:
\begin{align*}
&\ \lim_{m \to \infty} \left(\frac{\lambda^{k-d-1}}{\mu_2 \cdots
\mu_{k-d}} \right)^{m-l} \, m^2 \cdot p_m^{k-2} e_i(B_1^l)
e_{k-2-i}(D_2^l (\lambda^{-1} D)^{m-l}) \cdot \det (C_{m,l})\\
= &\ \lim_{m \to \infty} \left(\frac{\lambda^{k-d-1}}{\mu_2 \cdots
\mu_{k-d}} \right)^{m-l} \ m^2 \cdot p_m^{k-2} e_i(B_1^l) \cdot
\gamma_{k-2-i,l} \prod_{j=1}^{k-2-i} (\lambda^{-1} \mu_{j+1})^{m-l} \cdot
\det (C_{m,l})\\
= &\ \lim_{m \to \infty} \left(\frac{\mu_{k-d+1} \cdots
\mu_{k-i-1}}{\lambda^{d-i-1}} \right)^{m-l} \cdot p_m^{k-2}
e_i(B_1^l) \gamma_{k-2-i,l} \cdot m^2 \det (C_{m,l})\\
= &\ - \left( \frac{\lambda}{\lambda - \mu_1} \right)^2 |b|^2 a^{2l-2} l
\cdot e_i(B_1^l) \gamma_{k-2-i,l} \cdot \lim_{m \to \infty}
\left(\frac{\mu_{k-d+1} \cdots \mu_{k-i-1}}{\lambda^{d-i-1}}
\right)^{m-l}.
\end{align*}

Notice that the final limit is the indicator ${\bf 1}(i = d-1)$, and if
$i=d-1$ then $e_i(B_1^l) = \det(B_1)^l$. This concludes the proof.
\qedhere
\end{enumerate}
\end{proof}

For completeness, we also show that the bound of $d_{A^{-1}}$ is tight
for negative $m$:

\begin{theorem}\label{Tneg}
Given positive definite $A_{n \times n}$ with $d_A < n$ (so $d_{A^{-1}} <
n$) and an integer $l \geq 1$, there exists a positive definite $B_{n
\times n}$ such that we have ``eventual negativity'' above
the BMV-degree $d_{A^{-1}}$:
\[
e_k(\sigma((-1)^l S_{m,l}(A,B))) < 0, \qquad \forall k \in
(d_{A^{-1}},n],
\]
simultaneously for all $m < 0$ small enough.
\end{theorem}

\begin{proof}
Akin to the previous proof, let $\lambda := \lambda_{\max}(A^{-1}) =
\lambda_{\min}(A)^{-1}$ and $d := d_{A^{-1}}$. Let the second-largest
eigenvalue of $A^{-1}$ be $\mu$; now via simultaneous unitary
conjugation, suppose
\[
A^{-1} = \begin{pmatrix} \lambda & 0 \\ 0 & \mu \end{pmatrix} \oplus
\lambda {\rm Id}_{d-1} \oplus D,
\]
with $D$ positive definite and diagonal. Moreover, $0 < \mu < \lambda$
and $0 < D < \lambda \, {\rm Id}_{n-d-1}$.

Next, choose $a \in (0,\infty)$, $0 \neq b \in \mathbb{C}$, $c \in (
|b|^2/a, \infty)$, and any positive definite diagonal matrices $B_1 \in
\mathbb{C}^{(d-1) \times (d-1)}$ and $B_2 \in \mathbb{C}^{(n-d-1) \times
(n-d-1)}$. Now define
\[
B' := \begin{pmatrix} a & b \\ \overline{b} & c \end{pmatrix} \oplus B_1
\oplus B_2, \qquad B := \sqrt{A} B' \sqrt{A}.
\]
Then $A^{\pm 1}, B', B$ are all positive definite, and as in
Proposition~\ref{Psml},
\begin{equation}
(-1)^l S_{-m,l}(A,B) = T_{m-1+l,l} \oplus 
\textstyle\binom{m-1+l}{l} \lambda^m B_1^l \oplus
\textstyle\binom{m-1+l}{l} D^m B_2^l
\end{equation}
for $m>0$, where
\[
T_{m-1+l,l} := \begin{pmatrix} \sqrt{\lambda} & 0 \\ 0 & \sqrt{\mu}
\end{pmatrix} S_{m-1+l,l} \left( \begin{pmatrix} \lambda & 0 \\ 0 & \mu
\end{pmatrix}, \begin{pmatrix} a & b \\ \overline{b} & c \end{pmatrix}
\right) \begin{pmatrix} \sqrt{\lambda} & 0 \\ 0 & \sqrt{\mu}
\end{pmatrix}.
\]

We now compute $e_k$ of the block diagonal matrix $(-1)^l S_{-m,l}(A,B)$,
as in the proofs above. This is a sum of principal $k \times k$ minors,
and each minor includes zero, one, or two rows/columns from the first
two. Thus,
\begin{align*}
e_k((-1)^l S_{-m,l}(A,B)) = \sum_{i=0}^{d-1} &\ e_i(
\textstyle\binom{m-1+l}{l} \lambda^m B_1^l) \left[
e_{k-i-2}(\textstyle\binom{m-1+l}{l} D^m B_2^l) \det T_{m-1+l,l}
\right.\\ &\ \left. + e_{k-i-1}(\textstyle\binom{m-1+l}{l} D^m B_2^l)
{\rm Tr}(T_{m-1+l,l}) + e_{k-i}(\textstyle\binom{m-1+l}{l} D^m B_2^l)
\right]
\end{align*}

We now multiply and divide both appearances of the matrix $T_{m-1+l,l}$ by
$c_m := (m-1+l)^l \lambda^{m-1} / l!$. By
Proposition~\ref{Pasymptotics}(1), we get the matrix
\[
C'_{m-1+l,l} := \begin{pmatrix} \lambda a_{m-1+l,l} & \sqrt{\lambda \mu}
b_{m-1+l,l} \\ \sqrt{\lambda \mu} \, \overline{b}_{m-1+l,l} & \mu
c_{m-1+l,l} \end{pmatrix};
\]
note that $\det C'_{m-1+l,l} = \lambda \mu \cdot \det C_{m-1+l,l}$.
Moreover, this rescaling by $c_m$ rescales the trace-term in the same
way, and the $\det$-term by $c_m^2$.

Now divide $e_k((-1)^l S_{-m,l}(A,B))$ by
\[
m^{-2} (m^l/l!)^k \cdot \lambda^{2m-2 + m(d-1)} \cdot e_{k-d-1}(D^m
B_2^l)
\]
and take $m \to \infty$. This yields:
\begin{align*}
&\ \lim_{m \to \infty} \frac{e_k((-1)^l S_{-m,l}(A,B))}{m^{-2} (m^l/l!)^k
\cdot \lambda^{2m-2 + m(d-1)} \cdot e_{k-d-1}(D^m B_2^l)}\\
= &\  \sum_{i=0}^{d-1} e_i(B_1^l) \left[ \lim_{m \to \infty}
\frac{e_{k-i-2}( D^m B_2^l)}{e_{k-d-1}(D^m B_2^l) \lambda^{m(d-1-i)}}
\cdot \lim_{m \to \infty} \lambda \mu \, m^2 \det C_{m-1+l,l} \right.\\
&\ + \lim_{m \to \infty} \frac{e_{k-i-1}( D^m B_2^l) m^2}{e_{k-d-1}(D^m
B_2^l) \lambda^{m(d-i)-1}} \cdot \lim_{m \to \infty} (\lambda a_{m-1+l,l}
+ \mu c_{m-1+l,l})\\
&\ + \left. \lim_{m \to \infty} \frac{e_{k-i}( D^m B_2^l)
m^2}{e_{k-d-1}(D^m B_2^l) \lambda^{m(d+1-i)-2}} \right].
\end{align*}

Now consider a typical limit of the form
\[
\lim_{m \to \infty} \frac{m^\alpha e_a(D^m B_2^l)}{e_b(D^m B_2^l)
\lambda^{m(a-b)}}, \qquad 0 \leq b < a \leq n-d-1, \quad 0 \leq \alpha.
\]

By considering the growth rates (as in e.g.\ the proof of
Theorem~\ref{Ttight}), and since $0 < D < \lambda \, {\rm Id}_{n-d-1}$,
it follows that these limits all vanish. In particular, the only nonzero
limit in the previous calculation is for $e_{k-i-2}$ in the numerator,
when $i=d-1$. Now using Proposition~\ref{Pasymptotics}, we obtain:
\begin{align*}
&\ \lim_{m \to \infty} \frac{e_k((-1)^l S_{-m,l}(A,B))}{m^{-2} (m^l/l!)^k
\cdot \lambda^{2m-2 + m(d-1)} \cdot e_{k-d-1}(D^m B_2^l)}\\
= &\ e_{d-1}(B_1^l) \cdot \lambda \mu \cdot - \left(
\frac{\lambda}{\lambda - \mu} \right)^2 |b|^2 a^{2l-2} l < 0,
\end{align*}
and this is true for every $d_{A^{-1}} < k \leq n$, so one can find a
simultaneous threshold $m_l \gg 0$.
\end{proof}

Since we are exploring unconditional results in this section, we conclude
with one more, for completeness. Namely, for small $l$ we provide an
alternate proof (from the ones in the literature) of the BMV property, as
an application of our analytic combinatorics methods above.

\begin{prop}\label{PBMV12}
The BMV property holds for $l=1,2$ for all $m \geq l$.
\end{prop}

\begin{proof}
We first work with general $l,m$ and take a closer look at the
$\binom{m}{l}$ formal words of length $m$ in two positive letters $A,B$:
\[
\mathcal{W}_{m,l} := \{ W_{\bf j} := A^{j_1} B A^{j_2} B \cdots B
A^{j_{l+1}} \, : \, {\bf j} \in \mathbb{Z}_{\geq 0}^{l+1}, \, |{\bf j}| =
m-l \}.
\]

Since we only care about the traces of such words, define an equivalence
relation on $\mathcal{W}_{m,l}$ as the transitive closure of every word
$W_{\bf j} = M_1 \cdots M_m$ being equivalent to $M_m M_1 \cdots
M_{m-1}$. This partitions $\mathcal{W}_{m,l}$ into orbits. We now claim:

\begin{lemma}\label{Lcombin}
Inside every orbit, the proportion of words starting with $B$ is
precisely $l/m$.
\end{lemma}

\begin{proof}
Let $\rho : \mathcal{W}_{m,l} \to \mathcal{W}_{m,l}$ act on any word by
the above ``r(h)otation'' by one letter. Then the cyclic group
$\tangle{\rho} \cong \mathbb{Z} / m \mathbb{Z}$ acts on
$\mathcal{W}_{m,l}$. Fix an orbit $\mathcal{O}$ and a word $W_{\bf j}$ in
it, and let $1 \leq k \leq m$ be the smallest positive integer such that
$\rho^k(W_{\bf j}) = W_{\bf j}$. (For instance, if $W_{\bf j} = BAABAA$
then $m=6$, $l=2$, and $k=3$.) Now we claim that $k|m$. Indeed, for all
$a \geq 1$ we have $\rho^{ak}(W_{\bf j}) = \rho^{(a-1)k}(W_{\bf j})$,
which together with $\rho^{bm}(W_{\bf j}) = W_{\bf j}$ (for all $b \in
\mathbb{Z}$) yields $\rho^{\gcd(k,m)}(W_{\bf j}) = W_{\bf j}$. By
minimality we get $k = \gcd(k,m)$, so $k|m$.

Now let $W_{\bf j} = M_1 \cdots M_m$ and set $W' := M_1 \cdots M_k$. Then
$W_{\bf j} = (W')^{m/k}$. Moreover,
$\mathcal{O} = \{ W_{\bf j}, \rho(W_{\bf j}), \dots, \rho^{k-1}(W_{\bf
j}) \}$ by minimality of $k$. Thus, the number of words in $\mathcal{O}$
that start with $B$ is precisely the number of $B$'s that occur inside
$W'$, which is $l \cdot (k/m)$. Hence the desired proportion is
$\frac{lk/m}{k} = l/m$.
\end{proof}

Now we proceed with the proof using resolvents. We stick to the notation
and methods used in the proof of Theorem~\ref{Tnxn}. Namely, $A = {\rm
diag}(\lambda_i {\rm Id}_{n_i})_{i=1}^k$ and $B = (B_{ij})_{i,j=1}^k$,
with the eigenvalues $\lambda_1 > \cdots > \lambda_k \geq 0$.

Suppose first that $B_{11} = 0$. Then by positivity, $B_{1j}, B_{j1}$
both must vanish as well. In particular, $A,B$ become block diagonal
matrices. Write $A = \lambda_1 {\rm Id}_{n_1} \oplus A'$, $B = {\bf
0}_{n_1 \times n_1} \oplus B'$; then $S_{m,l}(A,B) = {\bf 0}_{n_1 \times
n_1} \oplus S_{m,l}(A',B')$ for $m \geq l \geq 1$. Thus we reduce to
working with $A',B'$. Continuing in this way, either $B = {\bf 0}_{n
\times n}$ (in which case the result is trivial), or we can assume
without loss of generality that $B_{11} \neq 0$.

We proceed in the latter case. Again if $\lambda_1 = 0$, then $A = {\bf
0}_{n \times n}$ and the result is trivial. So $\lambda_1 > 0$ and
$B_{11} \neq 0$ is psd. We once again compute asymptotics. By a result
similar to Proposition~\ref{Pgenfn}, and setting $R(z) = ({\rm Id}_n - z
A)^{-1}$,
\[
z^l (B R(z))^l = \sum_{m=l}^\infty z^m \sum_{|{\bf k}| = m-l} B A^{k_1}
\cdots B A^{k_l}.
\]

Note that this series converges on $D(0, 1/\|A\|)$. Taking traces at any
$z$ in this disc and invoking Lemma~\ref{Lcombin}, the coefficient of
$z^m$ is precisely $l/m$ times ${\rm Tr}(S_{m,l}(A,B))$. Thus, eventual
BMV positivity is equivalent to the eventual positivity of the
coefficients in ${\rm Tr}(B R(z))^l$. This is equivalent to working with
\[
T_l(z) := {\rm Tr}((\sqrt{B} R(z) \sqrt{B})^l) = {\rm Tr} \left(
\sum_{k=0}^\infty \sqrt{B} (zA)^k \sqrt{B} \right)^l.
\]

We now consider the two cases of interest. First if $l=1$, then
\[
[z^m]T_1(z) = {\rm Tr}(K_m), \qquad \text{where} \qquad K_m := \sqrt{B}
A^m \sqrt{B} \ \, \forall m \geq 0,
\]
and all $K_m$ are clearly psd if $A,B$ are. This shows the BMV property
for $l=1 \leq m$.

Next if $l=2$, then
\[
[z^m]T_2(z) = \sum_{j=0}^m {\rm Tr}(K_j K_{r-j}) = \sum_{k=0}^m {\rm Tr}
(\sqrt{K_j} K_{r-j} \sqrt{K_j}).
\]
Once again, each summand-matrix is psd, showing the BMV property for $l=2
\leq m$.
\end{proof}

\section{Application 5: Stronger Uchiyama result; answering
Ricker--V\"ath's 2002~question}\label{Suchi}

Our goal here is to prove the Uchiyama and Ricker--V\"ath results in
their vastly strengthened form -- followed by the Ogasawara and Pedersen
results, and then Theorem~\ref{Tadmissible} which similarly generalizes
our original characterization in Theorem~\ref{Thilbert}(1).

\subsection{Proof of Theorem~\ref{Tuchi}(1)}

We adapt here the proof of Theorem~\ref{Thilbert}. Thus, by ``Step~$l$''
we will mean the corresponding step in that proof.

Suppose $A^2B \neq BA^2$. Then Steps~1 and~2 go through verbatim, with
$\lambda_0 = 0$ in Lemma~\ref{Lstasi}. Next, the telescoping commutation
identity in Step~3 (of our ``BMV / Hurwitz-product'' case) is replaced in
this case by the ``brute force'' equation:
\begin{align}\label{Ebrute}
\begin{aligned}
A^n f(A) B = &\ W_n(A,B) - B A^n f(A) - 2 A^n r(A),\\
\text{where} \quad W_n(A,B) := &\ W_{\mathcal{P}, {\bf 0}, f, n}(A,B) +
W_{\mathcal{P}, {\bf 0}, f, n}(A,B)^*
\end{aligned}
\end{align}
is assumed to be positive for all $0 < n = n_k \uparrow \infty$.

We lower-bound the left-hand side in the same fashion as in Step~3,
using the assumptions on $f$:
\begin{equation}\label{Euchi1}
|\tangle{v, A^n f(A) B v_\mu}| = \left| \int_{(\varepsilon F)^\sigma}
\lambda^n f(\lambda) \, d E_{v,v}(\lambda) \right| = \int_{(\varepsilon
F)^\sigma} |\lambda^n f(\lambda)| \, d E_{v,v}(\lambda) \geq
|f(\varepsilon a)| \, a^n \|v\|^2 > 0.
\end{equation}

Next, upper-bound each term on the right -- except the first one -- using
Cauchy--Schwarz and Corollary~\ref{Cborel}(1), as in Step~4:
\begin{align}\label{Euchi2}
| \tangle{v, B A^n f(A) v_\mu} | \leq &\ \|v\| \, \|B\| \, \|f\|_{\infty,
[-c,c]^\sigma} \, \|v_\mu\| \, c^n,\\
| \tangle{v, 2 r(A) A^n v_\mu} | \leq &\ 2 \, \|v\| \cdot \| r(x) x^n
\|_{\infty, [-c,c]^\sigma} \cdot \|v_\mu\| \leq 2 \, \|v\| \, \|v_\mu\|
\, \|r\|_{\infty,[-c,c]^\sigma} \, c^n.\notag
\end{align}

Finally, if $W_n(A,B) \geq 0$ then it defines a semi-inner product for
which Cauchy--Schwarz holds. Thus,
\[
|\tangle{v, W_n(A,B) v_\mu}| \leq |\tangle{v,W_n(A,B) v}|^{1/2}
|\tangle{v_\mu, W_n(A,B) v_\mu}|^{1/2}.
\]
The first factor is again bounded above using~\eqref{Euchi},
Cauchy--Schwarz, and Corollary~\ref{Cborel}(1):
\begin{align*}
|\tangle{v,W_n(A,B) v}| \leq &\ |\tangle{A^n f(A)v, B v}| +
|\tangle{v, B A^n f(A)v}| + 2 |\tangle{v, r(A)A^nv}|\\
\leq &\ 2\, \|v\|^2 \left( \|B\| \, b^n |f(\varepsilon b)| +
\|r\|_{\infty, (\varepsilon F)^\sigma} b^n \right).
\end{align*}

Similarly, the second factor is bounded above by:
\[
|\tangle{v_\mu,W_n(A,B) v_\mu}| \leq 2 \, \|v_\mu \|^2 \left( \|B\|
\max(|f(c)|,|f(-c)|) c^n + \|r\|_{\infty, [-c,c]^\sigma} c^n \right).
\]

Putting together these estimates using~\eqref{Ebrute}, akin to Step~5 we
get:
\begin{equation}\label{Eestimate}
|f(\varepsilon a)| \, a^n \|v\|^2 \leq C_1 c^n + C_2 c^n + C_3
(bc)^{n/2},
\end{equation}
for some constants $C_1, C_2, C_3 \geq 0$ not depending on $n$, and for
infinitely many $n = n_k \uparrow \infty$ by hypothesis. Dividing by
$|f(\varepsilon a)| a^{n_k} > 0$ and taking $k \to \infty$ yields
$\|v\|^2 \leq 0$ (noting that $a > (b+c)/2 > \sqrt{bc}$). This
contradicts the choice of $v \neq 0$, showing that $A^2$ commutes with
$B$. \qed

\subsection{Proof of Theorem~\ref{Tuchi}(2)}

This is the sole instance in this work where we will invoke
Lemma~\ref{Lstasi} with $\lambda_0 > 0$. Indeed, by the hypotheses we can
choose $\lambda_0 = 1$. Now suppose $AB \neq BA$ and apply
Lemma~\ref{Lstasi} to $A' = A$, followed by Step~2, to obtain
\[
\lambda_0 = 1 \leq c < a \leq b \leq \| A \|, \qquad v_\mu \in {\rm
Ran}(E([\lambda_0, c]^\sigma)), \qquad 0 \neq v = E([a,b]^\sigma) B
v_\mu.
\]

Next, setting $\beta = \beta_k$ for general $k$, the analogue
of~\eqref{Ebrute} is:
\begin{align}\label{Ebrute2}
A^\beta f(A) B = &\ W_\beta(A,B) - B A^\beta f(A) - \sum_{j=1}^m (p_j(A)
C_j q_j(A) + q_j(A) C_j^* p_j(A)) - 2 A^\beta r(A),\notag\\
\text{where} \quad W_\beta(A,B) := &\ W_{\mathcal{P}, {\bf C}, f,
\beta}(A,B) + W_{\mathcal{P}, {\bf C}, f, \beta}(A,B)^*
\end{align}
is assumed positive for all $0 < \beta = \beta_k \uparrow \infty$. Apply
$\tangle{v, - \cdot v_\mu}$ to both sides and then the triangle
inequality, as in the previous part. This again yields, akin to the three
bounds in~\eqref{Euchi1} and~\eqref{Euchi2}:
\begin{align*}
|\tangle{v, A^\beta f(A) B v_\mu}| \geq &\ |f(a)| \, a^\beta \, \|v\|^2 >
0,\\
| \tangle{v, B A^\beta f(A) v_\mu} | \leq &\ \|v\| \, \|B\| \, |f(c)| \,
\|v_\mu\| \, c^\beta,\\
| \tangle{v, 2 r(A) A^\beta v_\mu} | \leq &\ 2 \, \|v\| \,
\|v_\mu\| \, \|r\|_{\infty,[1,c]^\sigma} \, c^\beta.
\end{align*}

We next bound the other terms on the right, starting with
\begin{align*}
&\ \left| \left\langle v, \sum_{j=1}^m (p_j(A) C_j q_j(A) + q_j(A) C^*_j
p_j(A)) \cdot v_\mu \right\rangle \right|\\
\leq &\ \sum_{j=1}^m |\tangle{p_j(A)v, C_j q_j(A) v_\mu}| + \sum_{j=1}^m
|\tangle{C_j q_j(A)v, p_j(A) v_\mu}|\\
\leq &\ \sum_{j=1}^m \|v\| \, \|C_j\| \, \|v_\mu\| \, ( 
\|p_j\|_{\infty, [a,b]^\sigma} \, \|q_j\|_{\infty, [1,c]^\sigma} +
\|q_j\|_{\infty,[a,b]^\sigma} \, \|p_j\|_{\infty, [1,c]^\sigma} ) =: C_0,
\end{align*}
say. As before, this is obtained using Corollary~\ref{Cborel}(1) and
Cauchy--Schwarz.

Finally, since $W_\beta(A,B) \geq 0$, we come to the semi-inner product
and Cauchy--Schwarz usage (as above) in order to upper-bound
\[
|\tangle{v, W_\beta(A,B) v_\mu}| \leq |\tangle{v,W_\beta(A,B) v}|^{1/2}
|\tangle{v_\mu, W_\beta(A,B) v_\mu}|^{1/2}.
\]

Akin to above, the first factor on the right (without the square root) is
bounded above by
\[
2 \, \|v\|^2 \left( \|B\| \, b^\beta \, |f(b)| + \sum_{j=1}^m \|C_j\| \,
\|p_j\|_{\infty, [a,b]^\sigma} \, \|q_j\|_{\infty, [a,b]^\sigma} +
\|r\|_{\infty, [a,b]^\sigma} \, b^\beta \right) \leq C_1 b^\beta
\]
for some $C_1 \geq 0$. Here we use the fact from Lemma~\ref{Lstasi} that
$\lambda_0 = 1 \leq c < b$. This same fact implies similarly that the
second factor on the right (without the square root) is bounded above by
$C_2 c^\beta$ for some $C_2 \geq 0$.

Collecting all of this together,
\[
0 < |f(a)| \, a^\beta \, \|v\|^2 \leq (\|B\| \, |f(c)| + 2 \,
\|r\|_{\infty,[1,c]^\sigma}) \|v\| \, \|v_\mu\| \, c^\beta + C_0 \cdot
1^\beta + \sqrt{C_1 C_2} (bc)^{\beta/2},
\]
for all $\beta = \beta_k \uparrow \infty$. Divide by $a^\beta \, |f(a)| >
0$ and let $\beta = \beta_k \to \infty$. As $a > (b+c)/2 > \sqrt{bc} > c
\geq 1 = \lambda_0$, we get $\|v\|^2 \leq 0$, which is the desired
contradiction. \qed

\subsection{Proof of Theorem~\ref{Tuchi}(3)}

This is omitted as it is very similar to the preceding proof.

\subsection{Proof of Theorem~\ref{Tuchi}(4)}

Theorem~\ref{Trv} follows from part~(1) by setting ${\bf C} =
{\bf 0}$, $f(\lambda) \equiv 1$, and $r(\lambda) \equiv d/2$.
Indeed, this gives $[A^2,B]=0$, and since $A \geq 0$, we get $[A,B]=0$ as
in the proof of Corollary~\ref{Ccstar}. Note here that we do not require
$B,d$ to be nonnegative.

Next, we claim that Theorem~\ref{Trv} also follows from the second part,
with the same ${\bf C}, f, r$ and $\beta_k = Nk$. This is achieved as
follows: set $A_1 := A^N + {\rm id}_{\mathcal{H}}$, so that $A_1$
satisfies the hypotheses in part~(2). Now compute:
\begin{align*}
&\ W_{\mathcal{P}, {\bf 0}, f, n}(A_1,B) + W_{\mathcal{P}, {\bf 0}, f,
n}(A_1,B)^* = A_1^n B + B A_1^n + d A_1^n\\
= &\ (A^N + {\rm id}_{\mathcal{H}})^n B + B (A^N + {\rm
id}_{\mathcal{H}})^n + d (A^N + {\rm id}_{\mathcal{H}})^n\\
= &\ \sum_{j=0}^n \binom{n}{j} ( A^{Nj} B + B A^{Nj} + d A^{Nj} ),
\end{align*}
and each term for $j>0$ is positive by the hypotheses of
Theorem~\ref{Trv}. So is the $j=0$ term, since $2B + d \, {\rm
id}_{\mathcal{H}} \geq 0$. Hence by Theorem~\ref{Tuchi}(2) above, $0 =
[A_1,B] = [A^N,B]$. Finally, approximate $x^{1/N}$ uniformly on
$[0,\|A^N\|]$ by a sequence of polynomials, and use this to conclude
(since $A \geq 0$) that $0 = [A,B]$.

Finally, setting $f(\lambda) \equiv 1$ and $r(\lambda) \equiv d/2$, we
recover \cite[Theorem~3]{RickerVath}.
\qed

\subsection{Ogasawara and Pedersen}

We next show the aforementioned results of Ogasawara and Pedersen, as a
consequence of Theorem~\ref{Tuchi}:

\begin{proof}[Proof of Theorem~\ref{Tmolnar}]\hfill
\begin{enumerate}
\setcounter{enumi}{2}
\item If $0 \leq H_0 \leq H$ in $\mathcal{A}$ and $[H_0, H] = 0$, then
\[
H^2 - H_0^2 = (H - H_0)^2 + 2 (H-H_0) H_0 = (H - H_0)^2 + 2 \sqrt{H_0} (H
- H_0) \sqrt{H_0} \geq 0.
\]

Conversely, by the GNS construction we may assume $\mathcal{A}
\hookrightarrow \mathcal{B}(\mathcal{H})$ for some complex Hilbert space
$\mathcal{H}$. Now for any $n \geq 1$, we have $H_1^n \geq 0$, and
\[
0 \leq (H_0 + \varepsilon H_1^n)^2 - H_0^2, \qquad \forall \varepsilon > 0.
\]
Dividing by $\varepsilon$ and then taking $\varepsilon \to 0^+$ yields:
$H_0 H_1^n + H_1^n H_0 \geq 0$ for all $n \geq 1$. But now
Theorem~\ref{Tuchi} (or even Uchiyama's original result) applies, to
give: $[H_0, H_1] = 0$. So $H_0$ commutes with all positive operators in
$\mathcal{A}$, hence with all self-adjoint operators, hence with all $H
\in \mathcal{A}$ -- by the same reasoning as in the proof of
Theorem~\ref{Tmolnar}(1).

The final assertion here is proved similarly: if every positive element
in $\mathcal{A}$ is central, then so is every self-adjoint element, and
hence every element of $\mathcal{A}$.

\item One way is clear: if $0 \leq H_0 \leq H$ in $\mathcal{A}$ and
$[H_0,H_1] = 0$, then setting $\Delta := H - H_0$, the functional
calculus (aka Gelfand--Naimark) for the commuting pair $H_0,\Delta \geq
0$ yields $\Psi : C(\sigma_{joint}(H_0,\Delta), \mathbb{R}) \to
\mathcal{B}(\mathcal{H})$. Letting $f(\lambda_1, \lambda_2) := (\lambda_1
+ \lambda_2)^p - \lambda_1^p$, we have by Spectral mapping that
\[
\sigma(\Psi(f)) = f(\sigma_{joint}(H_0,\Delta)) \subset [0,\infty),
\]
whence $\Psi(f) = (H_0 + \Delta)^p - H_0^p = H^p - H_0^p \geq 0$.
Conversely, by the Loewner--Heinz inequality (or simply Loewner's
monotonicity theorem \cite{Loewner34}), $\lambda \mapsto \lambda^{2/p}$
is Loewner monotone. Hence if $H^p \geq H_0^p$ then $H^2 \geq H_0^2$ for
all $0 \leq H_0 \leq H$ in $\mathcal{A}$. The working in the previous
part now shows $H_0$ is central in $\mathcal{A}$. \qedhere
\end{enumerate}
\end{proof}

\subsection{Proof of Theorem~\ref{Tadmissible}}

We conclude the main body of the paper by strengthening our original main
Theorem~\ref{Thilbert}(1) to Theorem~\ref{Tadmissible}. The proof is
along the same lines, but we first justify (1)~the convergence in norm of
the sums
\[
\sum_{i \geq 0} c_i S_{m+i,1}(A,B), \qquad
\sum_{i \geq 0} |c_i| S_{m+i,1}(A,B).
\]

In fact, we show that both  series are absolutely convergent. This is
immediate if $A=0$, else:
\begin{align}\label{Eabsconv}
\left\| \sum_{i \geq 0} |c_i| S_{m+i,1}(A,B) \right\| \leq &\ \sum_{i
\geq 0} |c_i| (m+i) \, \|B\| \, \|A\|^{m+i-1}\\
= &\ m \|A\|^{m-1} \|B\| \, \widetilde{f}(\|A\|) + \|A\|^m \, \|B\| \,
\widetilde{f}'(\|A\|) < \infty.\notag
\end{align}

We now show (2)~the main assertion -- we will only write the steps, as
these arguments have already shown up in the paper. Throughout, we will
assume $m \geq 1$. First, Steps~1 and~{2} in the proof of
Theorem~\ref{Thilbert}(1) go through verbatim, with $\lambda_0 = 0$ in
Lemma~\ref{Lstasi}. Now we use the algebraic identity in that Step~3,
$[A^m,B] = [A, S_{m,1}(A,B)]$, to obtain as in~\eqref{ECbeta00} (with
$\beta=1$):
\[
[A^m f(A),B] = [A, C^+_m], \qquad \text{where} \qquad
C^+_m := \sum_{i \geq 0} c_i S_{m+i,1}(A,B) = D(x^mf(x))[A](B).
\]

Then as in~\eqref{Ecommutator}
\[
|\tangle{v, A^m f(A) B v_\mu}| \leq |\tangle{v, B A^m f(A) v_\mu}| +
|\tangle{Av, C_m^+ v_\mu}| + |\tangle{v, C_m^+ A v_\mu}|.
\]

As usual, we lower-bound the left-side -- using the properties of
$\mathcal{F}_{adm}(A)$:
\[
| \tangle{v, A^m f(A) B v_\mu}| = \left| \int_{(\varepsilon F)^\sigma}
f(\lambda) \lambda^m \, d E_{v,v}(\lambda) \right| = \int_{(\varepsilon
F)^\sigma} |f(\lambda) \lambda^m | \, d E_{v,v}(\lambda) \geq
|f(\varepsilon a)| \, a^m \, \|v\|^2 > 0.
\]

Using Corollary~\ref{Cborel}(1) and Cauchy--Schwarz, the first term on
the right has upper bound
\[
| \tangle{v, B A^m f(A) v_\mu} | \leq \|B\| \, \|v\| \, \|v_\mu\| \,
\|f\|_{\infty,  [-c,c]^\sigma} c^m.
\]

The other two terms are bounded above as usual, using the semi-inner
product $\tangle{-, C^+_m \cdot -}$ on $\mathcal{H}$. For instance,
\[
|\tangle{Av, C^+_m v_\mu}| \leq |\tangle{Av, C^+_m Av}|^{1/2}
|\tangle{v_\mu, C^+_m v_\mu}|^{1/2}.
\]

By Corollary~\ref{Cborel}(2) and Cauchy--Schwarz,
\begin{align*}
|\tangle{Av, C^+_m Av}| \leq &\ \sum_{i \geq 0} |c_i| \cdot |\tangle{ Av,
S_{m+i,1}(A,B) Av}| \leq \sum_{i \geq 0} |c_i| (m+i) b^{m+i+1} \, \|B\|
\, \|v\|^2\\
= &\ m b^{m+1} \|B\| \, \widetilde{f}(b) \, \|v\|^2 + b^{m+2} \, \|B\| \,
\widetilde{f}'(b) \, \|v\|^2\\
\leq &\ \|B\| \, \|v\|^2 \, m b^{m+1} ( \widetilde{f}(b) + b
\widetilde{f}'(b)) < \infty.
\end{align*}

Similarly, the other term is upper-bounded by:
\[
|\tangle{v_\mu, C^+_m v_\mu}| \leq \|B\| \, \|v_\mu\|^2 \, m c^{m-1} (
\widetilde{f}(c) + c \widetilde{f}'(c)).
\]

From this, it follows that
\begin{align*}
|\tangle{Av, C^+_m v_\mu}| \leq &\ \|B\| \, \|v\| \, \|v_\mu\| \, m b
(bc)^{(m-1)/2} \sqrt{(\widetilde{f}(b) + b \widetilde{f}'(b))
(\widetilde{f}(c) + c \widetilde{f}'(c))}\\
= &\ mb (bc)^{(m-1)/2} \cdot C
\end{align*}
for some constant $C>0$ not depending on $m$. Using this $C$, one
estimates:
\[
|\tangle{v, C^+_m A v_\mu}| \leq mc (bc)^{(m-1)/2} \cdot C.
\]

Putting all of these estimates together,
\[
|f(\varepsilon a)| \, a^m \, \|v\|^2 \leq \|B\| \, \|v\| \, \|v_\mu\| \,
\|f\|_{\infty,  [-c,c]^\sigma} c^m + C (b+c) m (bc)^{(m-1)/2}.
\]

Finally, divide by $|f(\varepsilon a) \,|a^m$ and let $m = m_k \uparrow
\infty$, recalling that $a > (b+c)/2 > \sqrt{bc}$. As earlier, this
yields $\|v\|^2 \leq 0$, and this contradiction completes the proof. \qed

\addtocontents{toc}{\protect\vspace{1ex}}

\appendix
\section{(By Vamsi P.~Pingali)
Geometric application 6: Griffiths negativity of some constructions of
Hermitian metrics}\label{sec:geom}

In this section, we provide additional motivation from geometry for the
above results, and then show some geometric applications of them that can
be viewed as ``negative'' results about metrics on complex bundles.

The notion of positively curved Hermitian metrics on holomorphic line
bundles (or more generally, quasi-plurisubharmonic (qpsh) functions) has
played an important role in complex geometry, starting from the Kodaira
embedding and vanishing theorems. The algebro-geometric notion of
ampleness of a vector bundle (i.e., that the
$\mathcal{O}_{\mathbb{P}(E^*)}(1)$ bundle over $\mathbb{P}(E^*)$ is
ample) was pioneered and used to prove vanishing theorems
\cite{Hartshorne}. There is no unique differentio-geometric counterpart
to this notion of ampleness. Griffiths \cite{Griffiths} introduced a
notion of positivity (Griffiths (strict) positivity), which demands the
(strict) positivity of the real $(1,1)$-form $\langle v,\sqrt{-1}\Theta
v\rangle$ for all $v\neq 0$. The dual of a Griffiths positively curved
metric can be seen easily to be Griffiths negatively curved. Griffiths
also conjectured that ampleness implied the existence of a Griffiths
positively curved metric. This conjecture had some evidence in its favour
\cite{Mori, Ume, CF, Berndtsson, Liu, P, Xiao}, but was disproven
recently~\cite{DuXie}.

There are line bundle analogues of a purely algebraic criterion implying
a differentio-geometric version of positivity. The ``classical''
Nakai--Moishezon criterion (later improved by Demailly--Paun \cite{DP})
is an example. The proof of the Demailly--Paun criterion involved two
parts -- a PDE result (concentration of mass) and a regularization result
(involving Demailly approximation and gluing). This broad strategy was
successfully applied to other cases such as the pioneering work of Chen
\cite{Chen} on the J-equation and the deformed Hermitian--Yang--Mills
equation (and later by others \cite{Datar, Song, CLT, Ballal, forms}). If
one hopes that such a strategy can be used for equations involving
metrics on vector bundles (and related results generalising the
pluripotential theory of line bundles to higher-rank vector bundles), one
crucially needs the ability to regularize metrics on vector bundles.
Locally, Griffiths (strictly) negatively curved metrics can be
regularized using convolution whilst preserving Griffiths (strict)
negativity \cite{Raufi} and globally on compact manifolds due to
\cite{Mourougane}. However, neither of these constructions gives a
Richberg-style gluing technique on compact manifolds which is paramount
in pluripotential theory. The aim of this appendix is to study natural
constructions that are motivated from the line bundle case.

In the line bundle case, here are some commonly used ingredients for
regularization.

\begin{enumerate}
\item The \textit{maximum construction:}
If $\phi_1, \phi_2$ are $\omega$-qpsh, then so is $\max(\phi_1,\phi_2)$.

\item \textit{Demailly approximation:}
Roughly speaking, qpsh functions can be approximated using Bergman
kernels, i.e., functions of the form $\frac{1}{k}\ln (\sum_i \vert
f_{i,k} (z)\vert^2)$ where $f_{i,k}$ are holomorphic functions. Using
either the maximum construction or the partition-of-unity construction,
this result was globalized (the precise statements can be found in
\cite{Dem} for this method). Suppose $u_i$ are qpsh functions defined on
unit-balls covering a manifold such that the corresponding half-balls
continue to cover. Let $\rho_i$ be a partition-of-unity. If there exists
$\delta>0$ such that $u_i(p)\leq u_j(p)-\delta$ for any $i,j$ such that
$p\in U_i \cap U_j$, and $p$ is in the corresponding half-ball of $U_j$
and is outside the half-ball of $U_i$, then $\frac{1}{k} \ln
\left(\displaystyle \sum_i \rho_i e^{ku_i} \right)$ is qpsh with a small
loss of positivity when $k$ is large.

Thus akin to the method of Richberg, one can regularize locally, subtract
a small amount, and then use this lemma. This method was also used in a
somewhat different context in \cite{Plis}, for instance.
\end{enumerate}

We show that these particular constructions for vector bundles have
limitations. For instance, a maximum construction can be attempted, but
unsuccessfully so:

\begin{theorem}[The maximum construction, see
e.g.~\cite{Chen}]\label{thm:max}
Let $V$ be a finite-dimensional complex inner product space and $A,B$ be
inner products on $V$. Then there exists a unique inner product
$M=\max(A,B)$ satisfying the following properties.
\begin{enumerate}
\item $M$ is symmetric in $A,B$.
\item $M\geq A$ and $M\geq B$.
\item The inner product $A$ induces a Hilbert--Schmidt-type norm (that we
denote as $\| \cdot \|$) on endomorphisms. In this norm, $M$ minimises
$f_{A,B}(C) :=\| A^{-1} C- {\rm Id}\|_A^2 + \| A^{-1} C-A^{-1}B\|_A^2$
among all inner products $C$ such that $C\geq A$ and $C\geq B$.
\end{enumerate}

Moreover, for every integer $r\geq 2$, there exists a holomorphically
trivial complex bundle of rank $r$ over a neighbourhood of the origin in
$\mathbb{C}$ and two Griffiths strictly negatively curved Hermitian
metrics $h_1,h_2$ such that $\max(h_1,h_2)$ is not Griffiths negatively
curved.
\end{theorem}

This construction has appeared in the literature in the context of the
J-equation~\cite{Chen}. We can in fact generalize Theorem~\ref{thm:max}
to the following ``no-go'' result by reducing it to
Theorem~\ref{thm:max}.

\begin{theorem}[No general maximum works]\label{thm:genmax}
Let $\mathbb{P}_r^+ \subset \mathbb{C}^{r \times r}$ denote the set of
(Hermitian) positive definite $r \times r$ matrices (with $r \geq 2$).
Then there is no continuous function $G : \mathbb{P}_r^+ \times
\mathbb{P}_r^+ \to \mathbb{P}_r^+$ satisfying the following properties
(for all $A,B \in \mathbb{P}_r^+$):
\begin{enumerate}
\item $G(A,B)=G(B,A)$.
\item $G(PAP^*, PBP^*) = P G(A,B) P^*$ for all $P \in GL(r,\mathbb{C})$.
\item $G(A,B)\geq A$, with equality if $A\geq B$ in the Loewner order
(i.e., if $A-B$ is psd).
\item Let $\Omega\subseteq \mathbb{C}^d$ be an open set and $A,B: \Omega
\rightarrow \mathbb{P}_r^+$ be smooth metrics on the trivial bundle
$\Omega\times \mathbb{C}^r$ such that $A$ and $B$ are Griffiths
negatively curved on $\Omega$. Then $G(A,B)$ is a Griffiths negatively
curved continuous Hermitian metric on $\Omega$ in the sense that for
every local holomorphic section $s$, $s^* G(A,B) s$ is a
plurisubharmonic function.
\end{enumerate}
\end{theorem}

Theorem~\ref{thm:genmax} is proven with assistance from ChatGPT 5.5 Pro.
It implies that the method of Richberg has severe limitations when applied to vector bundles.

On the other hand, the Bergman kernel approaches suffer from the problem
of the inability to take ``roots'' of metrics on tensor powers of vector
bundles. Moreover, even if one attempts a more naive construction of
taking powers of endomorphisms, one runs into the following problem.

\begin{theorem}[Powers of metrics]\label{thm:powers}
There exists a holomorphically trivial complex bundle over a
neighbourhood of the origin in $\mathbb{C}$, and a Griffiths negatively
curved Hermitian metric $H$ on it such that the matrix product $H^m$ is
not Griffiths negatively curved for any $m \geq 2$.
\end{theorem}

\subsection{The proofs}

We begin by proving Theorem~\ref{thm:max}:

\begin{proof}[Proof of Theorem~\ref{thm:max}]
Choose a basis for $V$ that is simultaneously an orthonormal basis
$\mathcal{B}$ for $A$ and a basis where $B$ is diagonal. Then define
$M_{ij}=\max\left(1,b_{ii}\right)\delta_{ij}$. Obviously $M\geq A, B$. We
shall now prove that the $M$ attains the minimum of the ``cost function''
$f_{A,B}(C)$ described in the theorem. Indeed, note that in the basis
$\mathcal{B}$,
\[
f_{A,B}(C) = \sum_{i,j} \left( \vert C_{ij}-\delta_{ij}\vert^2+ \vert
C_{ij}-b_{ii} \delta_{ij}\vert^2\right) \geq f_{A,B}(M),
\]
with equality if and only if $C=M$. To prove that $M$ is symmetric in
$A,B$, choose another basis $f_i = \frac{e_i}{\sqrt{b_{ii}} }$. Then
$f_i$ are orthonormal for $B$ and the matrix for $A$ in this basis is
$\tilde{A}_{ij}=\frac{1}{b_{ii}}\delta_{ij}$. In this basis,
$\tilde{M}_{ij}=\max\left(1,\frac{1}{b_{ii}}\right)\delta_{ij}$. Hence
$\tilde{M}$ and $M$ represent the same inner product in different bases.

For the second part of the theorem, let for $i=1,2$, define real
positive definite matrices
\[
B_i= \begin{pmatrix} a & b_i \\ b_i & c_i \end{pmatrix}
\qquad \text{and} \qquad
A = \begin{pmatrix} \lambda & 0 \\ 0 & \mu \end{pmatrix},
\]
where $\lambda>1>\mu>0$. Thus, $a, c_i, ac_i-b_i^2 > 0$.
Motivated by the proof of Theorem~\ref{thm:Tposdef}, suppose
\begin{gather}
\frac{a c_2}{b_2^2} \frac{4\lambda}{\mu} - \left(\frac{\lambda}{\mu}+1
\right)^2>0.
\label{ineq:first}
\end{gather}

The proof of Theorem~\ref{thm:Tposdef} shows that $AB_2+B_2A$ is
positive definite. Assume that the parameters are chosen so that the
following holds:
\begin{gather}
4a c_1\lambda-\left (
b_1(1+\lambda)+\frac{(1-\lambda)(\mu+\lambda)(b_2-b_1)}{\mu-\lambda}
\right)^2<0.
\label{ineq:second}
\end{gather}
For instance, the choice $\lambda =4, \mu=\frac{1}{4}, B_1 = {\rm Id}_2,
B_2 = \begin{pmatrix} 1 & 2 \\ 2 & 20 \end{pmatrix}$
works.

Define $h_1=({\rm Id}+B_1 \vert z \vert^2)^2$ and $h_2=({\rm Id}+B_2
\vert z \vert^2)A ({\rm Id}+B_2 \vert z \vert^2)$. Note that
\begin{gather}
h_1 \Theta_{h_1} (0)=-2B_1 dz\wedge d\bar{z} \nonumber \\
h_2 \Theta_{h_2} (0)=-(B_2A+AB_2) dz\wedge d\bar{z}.
\end{gather}

By the choices made above, $h_1$ and $h_2$ are Griffiths negatively
curved on some neighbourhood of the origin. Moreover, $h_2(0)$ has
distinct eigenvalues that are distinct from $1$. Now define
\begin{gather}
h:= \max(h_1,h_2) \nonumber \\
=({\rm Id}+B_1 \vert z \vert^2) \cdot \max\left({\rm Id}, ({\rm Id}+B_1
\vert z \vert^2)^{-1} h_2 ({\rm Id}+B_1 \vert z \vert^2)^{-1} \right)
\cdot ({\rm Id}+B_1 \vert z \vert^2).
\end{gather}

Note that the metric $h_3=({\rm Id}+B_1 \vert z \vert^2)^{-1} h_2 ({\rm
Id}+B_1 \vert z \vert^2)^{-1} $ has the same eigenvalues as $h_2(0)$ at
the origin. Thus, by the implicit function theorem, there is a
neighbourhood around the origin wherein a smooth real orthogonal matrix
$U$ can diagonalise $h_3$ to get a smooth diagonal matrix $D =
\begin{pmatrix} \alpha & 0 \\ 0 & \beta \end{pmatrix} $ such that
$h_3=UDU^T$. It is easy to see that
\[
U={\rm Id}+V\vert z \vert^2+ \mathrm{higher \ order \ terms}
\]
and
\[
D=A+C \vert z \vert^2+\mathrm{higher \ order \ terms},
\]
where $C$ is diagonal and $V=-V^T$. Comparing the power series for $h_3$
and $UDU^T$, we see that
\[
V= \begin{pmatrix} 0 & \frac{(\mu+\lambda)(b_2-b_1)}{\mu-\lambda} \\
- \frac{(\mu+\lambda)(b_2-b_1)}{\mu-\lambda} & 0 \end{pmatrix}, \qquad
C= \begin{pmatrix} 0 & 0 \\ 0 & 2\mu (c_2-c_1)
\end{pmatrix}.
\]

It follows that
\[
h=({\rm Id}+B_1\vert z \vert^2) U \tilde{D} U^T ({\rm Id}+B_1 \vert z
\vert^2),
\]
for some diagonal smooth function $\tilde{D}$ such that
$\tilde{D}(0)=\begin{pmatrix} \lambda & 0 \\ 0 & 1
\end{pmatrix}$.

Hence,
\begin{gather}
-h\Theta_h (0)=
 \begin{pmatrix} 2a \lambda &
(\lambda+1)b_1+\frac{(\mu+\lambda)(b_2-b_1)}{\mu-\lambda}(1-\lambda) \\
(\lambda+1)b_1+\frac{(\mu+\lambda)(b_2-b_1)}{\mu-\lambda}(1-\lambda) &
2c_1 \end{pmatrix} dz\wedge d\bar{z}.
\label{eq:curvofmax}
\end{gather}

From this and the choice of the parameters in \eqref{ineq:second}, we see
that $h$ is not Griffiths negatively curved. One can upgrade this
construction to arbitrary rank $r \geq 2$ trivial bundles by taking a
direct sum of $h_1, h_2$ with $e^{\vert z \vert^2} {\rm Id}_{r-2}$.
\end{proof}

We now prove Theorem~\ref{thm:genmax} by reducing it to
Theorem~\ref{thm:max}.

\begin{proof}[Proof of Theorem~\ref{thm:genmax}]
Suppose $G$ exists. Let $g(X)=G({\rm Id},X)$. Our first aim is to prove
that $g$ is monotone in the Loewner order, i.e., if $0<X\leq Y$, then
$g(X)\leq g(Y)$.

Consider $B(z)=X+\vert z \vert^2(Y-X)$. Note that $B$ is smooth and is
Griffiths negatively curved on the closed unit disc. Hence so is
$H=g(B(z))$ by the hypotheses.
This implies that for any fixed vector $v$, $v^{*}H(z) v$ is subharmonic.
On $\vert z \vert =1$, $B(z)=Y$, and on $z=0$, $B(z)=X$. Thus by the
maximum principle, $v^* g(Y)v \geq v^* g(X) v$ for all $v$ and hence
$g(Y)\geq g(X)$.

Now suppose $D$ is a diagonal positive-definite matrix with diagonal
entries $D_{jj}>0$ for $1 \leq j \leq r$. Let $M$ be a diagonal matrix
with entries $M_{jj}=\max(D_{jj},1)$. Then $D\leq M, {\rm Id}_r \leq M$.
We want to prove that indeed $g(D)=M$. Using the
invariance-under-conjugation assumption, this will prove that $G$
coincides with the construction in the proof of Theorem~\ref{thm:max},
and hence such a function $G$ does not exist.

Indeed, since $M=g(M)\geq g(D)\geq D, {\rm Id}_r$, we see that
$(g(D))_{jj}=M_{jj}$. Now $M-g(D)\geq 0$ and has zero diagonal entries,
so $g(D)=M$.
\end{proof}

Lastly, we prove the outstanding result:

\begin{proof}[Proof of Theorem~\ref{thm:powers}]
Let $H=A+B \vert z \vert^2$ where $A$ and $B$ are chosen using the second
part of Theorem~\ref{thm:Tposdef}. Thus $H$ is Griffiths negatively
curved. Then $H^m \Theta_{H^m}(0)=-S_{m,1}(A,B) dz\wedge d\bar{z}$.
Hence, for all $m\geq 2$, $H^m$ is not Griffiths negatively curved at the
origin.
\end{proof}

\section{Standard notions and basic proofs}\label{Appendix}

\subsection{Spectral measures and the Spectral theorem}\label{Sspectral}

We recall here the notion of a spectral measure, followed by the Spectral
theorem. For a reference, see \cite[Chapter~5]{Theo}.

\begin{defn}\label{Dspectralproj}
Fix a set $X$, a $\sigma$-algebra $\mathcal{F} \subseteq 2^X$, and a
complex Hilbert space $(\mathcal{H}, \tangle{\cdot,\cdot})$ with the
\textbf{global convention} that the inner product is $\mathbb{C}$-linear
in the second argument.

A \textit{spectral measure} or \textit{projection-valued measure} for the
triple $(X, \mathcal{F}, \mathcal{H})$ is a map $E: \mathcal{F} \to
\mathcal{B}(\mathcal{H})$ satisfying:
\begin{enumerate}
\item For each $F \in \mathcal{F}$, $E(F)$ is an orthogonal projection
(i.e., it is idempotent and self-adjoint).

\item $E(\emptyset) = 0$ and $E(X) = {\rm id}_{\mathcal{H}}$.

\item $E(F_1 \cap F_2) = E(F_1) E(F_2)$ for $F_1, F_2 \in \mathcal{F}$.
(In particular, the projections $E(F), F \in \mathcal{F}$ pairwise
commute.)

\item If $F_1, F_2, \dots$ are pairwise disjoint sets in $\mathcal{F}$,
then the partial sums $\sum_{n=1}^N E(F_n)$ converge to
$E(\coprod_{n=1}^\infty F_n)$ in the SOT, as $N \to \infty$.
\end{enumerate}
\end{defn}

A spectral measure yields a collection of countably additive
$\mathbb{C}$-valued measures on $\mathcal{F}$ as follows:

\begin{prop}
Fix $(X, \mathcal{F}, \mathcal{H})$ as above and let $E$ be a spectral
measure. Given vectors $g,h \in \mathcal{H}$, the map
\[
E_{g,h} : \mathcal{F} \to \mathbb{C}, \qquad F \mapsto \tangle{E(F)g, h}
\]
defines a countably additive $\mathbb{C}$-valued measure on
$(X,\mathcal{F})$ with total variation at most $\| g \| \cdot \| h \|$.
If moreover $g=h$, then $E_{g,g}$ is a positive measure.
\end{prop}

We will also require the notion of the continuous functional calculus,
which requires some notation:

\begin{defn}\label{Djointspec}
Suppose $\mathcal{H}$ is a complex Hilbert space, $d \geq 1$, and $A_1,
\dots, A_d \in \mathcal{B}(\mathcal{H})$ are pairwise commuting
self-adjoint operators. Define the {\em spectrum} of any of the $A_j$ to
be:
\begin{equation}\label{Espec}
\sigma(A_j) := \{ \lambda \in \mathbb{C} : A_j - \lambda \, {\rm
id}_{\mathcal{H}} \text{ is not invertible} \}.
\end{equation}

Also define the {\em joint spectrum} $\sigma_{joint}({\bf A})$ of the
tuple ${\bf A} := (A_1, \dots, A_d)$ to consist of those $\bla =
(\lambda_1, \dots, \lambda_d) \in \mathbb{C}^d$ such that there exists a
sequence of ``simultaneous unit approximate-eigenvectors'' $x_n \in
\mathcal{H}$, satisfying $\lim_{n \to \infty} \| (A_j - \lambda_j {\rm
id}_{\mathcal{H}})x_n \| = 0$ for all $1 \leq j \leq d$.
\end{defn}

\begin{theorem}\label{Tfnalcalc}
Let $d \geq 1, \mathcal{H}$, and commuting $A_j = A_j^* \in
\mathcal{B}(\mathcal{H})$ be as in Definition~\ref{Djointspec}. Then:
\begin{enumerate}
\item The joint spectrum in the special case $d=1$ for the single
operator $A_j$, equals the spectrum $\sigma(A_j)$ as in~\eqref{Espec}.
More generally,
\[
\sigma_{joint}({\bf A}) = \left\{ \bla \in \mathbb{R}^d : 0 \in \sigma
\left( \textstyle\sum_{j=1}^d (A_j - \lambda_j {\rm id}_{\mathcal{H}})^2
\right) \right\}.
\]

\item $\sigma_{joint}({\bf A})$ is a compact subset of $\prod_{j=1}^d
\sigma(A_j) \subseteq \times_{j=1}^d [-\|A_j\|, \|A_j\|] \subset
\mathbb{R}^d$, and each $j$th projection map $: \mathbb{R}^d \to
\mathbb{R}$ sends $\sigma_{joint}({\bf A})$ onto $\sigma(A_j)$.

\item {\em (Gelfand--Naimark.)} There is an isometric $*$-algebra
monomorphism $\Psi : (C(\sigma_{joint}({\bf A})), \|~\cdot~\|_\infty) \to
(\mathcal{B}(\mathcal{H}), \| \cdot \|)$ satisfying:
$\Psi( {\bf z} \mapsto z_j ) = A_j$ and
$\Psi( {\bf z} \mapsto 1 ) = {\rm id}_{\mathcal{H}}$.

\item {\em (Spectral mapping.)} $\sigma(\Psi(f)) = f(\sigma_{joint}({\bf
A})) := \{ f(\bla) : \bla \in \sigma_{joint}({\bf A}) \}$.
\end{enumerate}
\end{theorem}

\begin{remark}\label{Rfnalcalc}
Given any continuous map $f$ on the spectrum, we define $f({\bf A}) :=
\Psi(f)$ via Theorem~\ref{Tfnalcalc}. In particular, if $f$ is the
uniform limit of polynomials $p_k$ on the compact set
$\sigma_{joint}({\bf A})$, then $f({\bf A}) = \lim_{k \to \infty}
p_k({\bf A})$ in supnorm, hence in the strong operator topology.
\end{remark}

With these preliminaries, we can now state the

\begin{theorem}[Spectral theorem]\label{Tspectral}
Let $d \geq 1, \mathcal{H}$, and commuting $A_j = A_j^* \in
\mathcal{B}(\mathcal{H})$ be as in Definition~\ref{Djointspec}. Then:
\begin{enumerate}
\item There exists a unique spectral measure $E$, defined on Borel
subsets of $\sigma_{joint}({\bf A})$, such that $\displaystyle A_j =
\int_{\sigma_{joint}({\bf A})} \lambda_j \, d E(\bla)$ for $1 \leq j \leq
d$. Moreover, if $U \subseteq \mathbb{R}^d$ is an open set intersecting
the joint spectrum nontrivially, then $E(U \cap \sigma_{joint}({\bf A}))
\neq 0$.

\item More strongly, for all continuous maps $f : \sigma_{joint}({\bf A})
\to \mathbb{C}$, one has
\begin{equation}
f({\bf A}) = \int_{\sigma_{joint}({\bf A})} f(\bla) \, d E(\bla).
\end{equation}
This is defined as meaning:
\[
\tangle{ g, f({\bf A}) h} = \int_{\sigma_{joint}({\bf A})} f(\bla) \, d
E_{g,h}(\bla), \quad \forall g,h \in \mathcal{H}.
\]

\item For all (not necessarily self-adjoint) operators $B \in
\mathcal{B}(\mathcal{H})$, we have $A_jB=BA_j$ for all $j$ if and only if
$E(F) B = B \, E(F)$ for all Borel subsets $F \subseteq
\sigma_{joint}({\bf A})$.
\end{enumerate}
\end{theorem}

\begin{example}
As a special case, suppose $\mathcal{H} = \mathbb{C}^n$ for some $n \geq
1$, and $A = A^*$ is an $n \times n$ matrix with eigenvalues $\lambda_1,
\dots, \lambda_k$ (for $k \leq n$) with eigenspaces $V_1, \dots, V_k$
respectively. Let $P_{\lambda_j}$ denote the projection onto $V_j$. Then
the spectral measure corresponds to the counting measure on $\sigma(A) =
\{ \lambda_1, \dots, \lambda_k \}$, and the Spectral theorem says that:
\[
E(\{ \lambda_j \}) = P_{\lambda_j}, \qquad A = \sum_{j=1}^k \lambda_j
P_{\lambda_j}, \qquad f(A) = \sum_{j=1}^k f(\lambda_j) P_{\lambda_j}.
\]

More generally, if $A_j = A_j^*$ are pairwise commuting $n \times n$
matrices for $1 \leq j \leq d$, then they are simultaneously
diagonalizable on $\mathbb{C}^n$, i.e.\ there is a basis of simultaneous
eigenvectors with finitely many joint eigenvalues $\bla$. These comprise
$\sigma_{joint}({\bf A})$, and $E(F) = \sum_{\bla \in F \cap
\sigma_{joint}({\bf A})} P_{\bla}$ for all Borel subsets $F \subset
\mathbb{R}^d$. \qed
\end{example}

\begin{example}\label{Excompact}
Even more generally, suppose $A_j = A_j^* \in \mathcal{B}(\mathcal{H})$
commute pairwise and are all compact -- so $\sigma(A_j) \setminus \{ 0
\}$ is countable, discrete, and consists of eigenvalues. Then all
$\sigma(A_j)$ and hence $\sigma_{joint}({\bf A})$ are countable sets, by
Theorem~\ref{Tfnalcalc}(2). \qed
\end{example}

We also state a consequence of the Spectral theorem, used above.

\begin{corollary}\label{Cmapping}
A tuple of commuting self-adjoint operators $A_j \in
\mathcal{B}(\mathcal{H})$ are all positive if and only if
$\sigma_{joint}({\bf A}) \subseteq \times_{j=1}^d [0, \|A_j\|]$. In this
case, $E_{v,v}$ is a nonnegative real measure on $\sigma(A)$, with total
mass $ \leq \|v\|^2$.
\end{corollary}

\subsection{Weyl functional calculus}\label{Sweylappendix}

We briefly describe here the Weyl functional calculus; see first the
opening paragraphs of Section~\ref{Sweylintro}.

Suppose $\mathcal{H}$ is a complex Hilbert space, $d \geq 2$ an integer,
and ${\bf A} = (A_1, \dots, A_d) \in \mathcal{B}(\mathcal{H})^d$ a tuple
of self-adjoint operators. Given a smooth map $f$ (even a polynomial) of
$d$ variables, as the $A_j$ need not commute pairwise, one cannot
directly apply the Spectral theorem to compute $f({\bf A})$. A workaround
is the Weyl calculus: note that $\tangle{\xi, {\bf A}} := \sum_{j=1}^d
\xi_j A_j$ is self-adjoint for all real tuples $\xi \in \mathbb{R}^d$,
and so $\exp( {\bf i} \tangle{\xi,{\bf A}} )$ is unitary. Now if $f \in
L^1(\mathbb{R}^d)$ is rapidly decreasing, then using its Fourier
transform we define the Bochner integral
\begin{equation}\label{Eweyl}
f({\bf A}) := \frac{1}{(2 \pi)^d} \int_{\mathbb{R}^d} \exp( {\bf i}
\tangle{\xi,{\bf A}} ) \widehat{f}(\xi) \, d \xi.
\end{equation}

This formula was originally suggested by Weyl~\cite{Weyl}, for ${\bf A} =
(Q,P)$, the coordinate and momentum operators in quantum mechanics. Today
it is known as the \textit{Weyl functional calculus}, and it extends to
smooth maps:

\begin{theorem}
Given a self-adjoint tuple ${\bf A} \in \mathcal{B}(\mathcal{H})^d$,
there exists a unique operator-valued distribution with compact support
$\mathcal{W}_{\bf A} : f \mapsto f({\bf A})$ for $f \in
C^\infty(\mathbb{R}^d)$, whose restriction to the Schwarz space of
rapidly decreasing $f$ is precisely the Weyl functional
calculus~\eqref{Eweyl} -- equivalently, $\mathcal{W}_{\bf A}$ is linear
and is given on the monomial ${\bf x}^{\bf m} := x_1^{m_1} \cdots
x_d^{m_d}$ via:
\[
\mathcal{W}_{\bf A}({\bf x}^{\bf m}) = \frac{m_1! \cdots m_d!}{|{\bf
m}|!} \sum_{W\, {\rm word}} W, \qquad \forall m_1, \dots, m_d \in
\mathbb{Z}_{\geq 0},
\]
where $|{\bf m}| := \sum_{j=1}^d m_j$, and $W$ runs over all distinct
words with $m_j$ copies of $A_j$ for each $j$.
\end{theorem}

In particular, for $d=2$ and integers $0 \leq l \leq m$, we have
\[
\mathcal{W}_{A,B}(x^{m-l}y^l) = \frac{1}{\binom{m}{l}} S_{m,l}(A,B).
\]

A comprehensive account of this calculus is provided in Jefferies'
monograph~\cite{Jeff}. We list here some basic properties:
\begin{enumerate}
\item For matrices $A_j$: The support of the above Weyl distribution is
within the convex hull of the joint numerical range
\[
\{ (\tangle{A_1 v, v}, \dots, \tangle{A_d v,v}) : v \in \mathcal{H}, \| v
\| = 1 \},
\]
and hence inside $\overline{D(0,\rho)}$, where $\rho^2 = \sum_{j=1}^d
\|A_j\|^2$. This is shown in~\cite{And1} via the Lie--Trotter formula and
the Paley--Wiener theorem.

\item For $d=2$, Anderson shows in \cite[Theorem~5.2]{And1} that for
operators, one needs to take the closure:
\begin{equation}\label{Enumrange}
{\rm supp} \, \mathcal{W}_{A,B} \subseteq \overline{\mathcal{C}}, \quad
\text{where} \quad \mathcal{C} := 
\{ (\tangle{A v, v}, \tangle{B v,v}) : v \in \mathcal{H}, \| v \| = 1 \}.
\end{equation}
(The closure is already necessary for one diagonal operator with infinite
discrete spectrum.)

\item If all $A_j$ are positive, then the support of the Weyl
distribution is within $\prod_{j=1}^d [0, \| A_j \|]$.

\item If the operators $A_j$ pairwise commute, the Weyl distribution
becomes the joint spectral measure, and the map $\mathcal{W}_{\bf A} : f
\mapsto f({\bf A})$ becomes multiplicative (i.e., an algebra map).
\end{enumerate}

This functional calculus was studied by Nelson (see his
article~\cite{Nelson} in a remarkable collection) and perhaps due to his
enthusiasm, by Anderson~\cite{And1} and Taylor~\cite{Taylor} --
coincidentally all at Princeton, as were Weyl and Wigner! -- and later as
well by Hein\"avaara~\cite{Heina} in going towards constructing joint
spectral measures.

\subsection{Proofs of some results above}\label{Appendix3}

Here we provide the proofs of some of the results above. We first show:

\begin{lemma}\label{Leigen}
Suppose $\mathcal{H}$ is a real or complex Hilbert space, and $A = A^*
\in \mathcal{B}(\mathcal{H})$ has discrete nonzero spectrum (possibly
empty). Then $A$ has a complete orthonormal eigenbasis for $\mathcal{H}$.
\end{lemma}

\begin{proof}
This is standard: first note that the compact set $\sigma(A) \subset
[-\|A\|, \|A\|]$ is countable, say $\sigma(A) \setminus \{ 0 \} = \{
\lambda_n : n \geq 1 \}$. Then ${\rm id}_{\mathcal{H}} = E(\mathbb{R}) =
E(\{ 0 \}) + \sum_{n \geq 1} E(\{ \lambda_n \})$. Since $\{ \lambda_n \}$
is discrete, there is an open interval $U_n$ with $U_n \cap \sigma(A) =
\{ \lambda_n \}$, so $E(\{ \lambda_n \}) \neq 0$. Hence it has nonzero
range $V_{\lambda_n}$, aka the $\lambda_n$-eigenspace of $A$. Finally,
$\mathcal{H} = V_0 \oplus \bigoplus_{n \geq 1} V_{\lambda_n}$. Now choose
a complete orthonormal basis of each (possibly infinite-dimensional)
subspace here.
\end{proof}

The following proof extends our results from complex Hilbert spaces to
real ones.

\begin{proof}[Proof of Theorem~\ref{Treal}]
A look at Theorems~\ref{Thilbert} and~\ref{Tnxn} for operators shows that
we need to understand how the following change under complexification:
$S_{m,l}(A,B)$; positivity; commutation; spectrum; and eigenspaces. This
is a matter of bookkeeping about complexification, and we refer to e.g.\
\cite[Section 5.1.1]{Theo} for the following standard facts.

Given a real Hilbert space $(\mathcal{H}_{\mathbb{R}},
\tangle{\cdot,\cdot}_{\mathbb{R}})$, define its complexification to be
$\mathcal{H} := \mathcal{H}_{\mathbb{R}} \otimes_{\mathbb{R}} \mathbb{C}
= \mathcal{H}_{\mathbb{R}} + {\bf i} \mathcal{H}_{\mathbb{R}}$, with
scalar multiplication and inner product given by:
\begin{align*}
(s+ {\bf i} t) \cdot (g+ {\bf i}h) := &\ (sg-th) + {\bf i} (sh+tg),\\
\tangle{g_1 + {\bf i} h_1, g_2 + {\bf i} h_2}
:= &\ ( \tangle{g_1, g_2}_{\mathbb{R}} + \tangle{h_1, h_2}_{\mathbb{R}})
+ {\bf i} ( \tangle{g_1, h_2}_{\mathbb{R}} - \tangle{h_1,
g_2}_{\mathbb{R}} ).
\end{align*}

Next, given a linear operator $A : \mathcal{H}_{\mathbb{R}} \to
\mathcal{H}_{\mathbb{R}}$, define its complexification via:
$A^{\mathbb{C}} (g + {\bf i} h) := Ag + {\bf i} Ah$. Then we have:
\begin{itemize}
\item If $A$ is bounded (resp.\ self-adjoint or positive) then so is
$A^{\mathbb{C}}$. The converse is easy.

\item $(A+B)^{\mathbb{C}} = A^{\mathbb{C}} + B^{\mathbb{C}}$ and
$(AB)^{\mathbb{C}} = A^{\mathbb{C}} B^{\mathbb{C}}$ for all $A,B$.

\item If $V \subseteq \mathcal{H}_{\mathbb{R}}$ is a closed subspace,
then so is $V^{\mathbb{C}} \subseteq \mathcal{H}$, and
$(P_V)^{\mathbb{C}} = P_{V^{\mathbb{C}} }$.

\item $\sigma(A^{\mathbb{C}}) = \sigma(A) \subseteq [-\|A\|, \|A\|]$ if
$A = A^*$.

\item If $A = A^*$ and $\mathcal{H}_{\mathbb{R}} = \bigoplus_{\lambda \in
\sigma(A)} V_\lambda$ is a direct sum of (closed) eigenspaces for the
real operator $A$, then $\mathcal{H} = \bigoplus_{\lambda \in \sigma(A)}
V_\lambda^{\mathbb{C}}$.
\end{itemize}

Now the proofs are along expected lines. For instance, if $S_{m,l}(A,B)
\geq 0$, then complexifying yields that the operator
\[
S_{m,l}(A,B)^{\mathbb{C}} = S_{m,l}(A^{\mathbb{C}}, B^{\mathbb{C}})
\]
is also positive, with the operators $A^{\mathbb{C}}, B^{\mathbb{C}}$
being self-adjoint. If this held for infinitely many $m \geq l$, then by
the complex case, $(A^{\mathbb{C}})^2 = (A^2)^{\mathbb{C}},
B^{\mathbb{C}}$ commute, whence so do their restrictions $A^2,B$ to
$\mathcal{H}_{\mathbb{R}}$.
Similarly and moreover, if $0 \supseteq \sigma(A) \cap -\sigma(A) =
\sigma(A^{\mathbb{C}}) \cap -\sigma(A^{\mathbb{C}})$, then
$A^{\mathbb{C}}, B^{\mathbb{C}}$ commute, whence so do $A,B$.
\end{proof}

Next, we elaborate on an example above.

\begin{proof}[Proof in Example~\ref{Exeventual}]
We show that with the matrices $A,\mathscr{B}$ as in the example,
$e_2(S_{m,1}(A,\mathscr{B})) < 0$ for all $m \geq 2$.
Indeed, applying~\eqref{Ecompute2x2} yields:
\[
S_{m,1}(A,\mathscr{B}) = \begin{pmatrix} 2 ma^{m-1} & \frac{a^m-1}{a-1}
\\ \frac{a^m-1}{a-1} & m \end{pmatrix} \oplus m {\rm Id}_{n-2}.
\]
Now $e_2(M)$ is the sum of the principal $2 \times 2$ minors of $M$, so
for all $m \geq 2$,
\begin{align*}
e_2(S_{m,1}(A,\mathscr{B}))
= &\ \binom{n-1}{2} m^2 + 2 (n-1) m^2 a^{m-1} - (1 + a + \cdots +
a^{m-1})^2\\
< &\ 2 \binom{n}{2} m^2 a^{m-1} - a^{2(m-1)} \leq a^{m-1} \left( 2
\binom{n}{2} m^2 - \binom{n}{2}^{m-1} 8^{m-1} \right)\\
\leq &\ a^{m-1} \binom{n}{2} (2 m^2 - 8^{m-1}).
\end{align*}
Using the elementary inequality $8^{m-1} \geq 2 m^2$ for all $m \geq 2$,
the final expression is $\leq 0$. Hence $e_2(S_{m,1}(A,\mathscr{B})) < 0$
for $m \geq 2$.
\end{proof}

\begin{proof}[Proof of Proposition~\ref{Pgenfn}]
We first show that $G_l$ exists on $D(0,1/\|A\|)$. Write
\[
G_l(z) = z^l \cdot \sum_{j_1 \geq 0} (zA)^{j_1} \cdot B \cdot \sum_{j_2
\geq 0} (zA)^{j_2} \cdot B \ \cdots \ B \cdot \sum_{j_{l+1} \geq 0}
(zA)^{j_{l+1}}.
\]

This multi-series has partial sums
\[
P_{\bf n}(z) := z^l \cdot \sum_{0 \leq j_1 \leq n_1} (zA)^{j_1} \cdot B
\cdot \sum_{0 \leq j_2 \leq n_2} (zA)^{j_2} \cdot B \ \cdots \ B \cdot
\sum_{0 \leq j_{l+1} \leq n_{l+1}} (zA)^{j_{l+1}}
\]
for all integers $n_1, \dots, n_{l+1} \geq 0$, where ${\bf n} = (n_1,
\dots, n_{l+1})$. The error/tail is
\[
E_{\bf n}(z) := G_l(z) - P_{\bf n}(z),
\]
and we \textbf{claim} that $E_{\bf n}(z) \to 0$ in operator norm, as
${\bf n} \to \infty^{l+1}$. This claim shows the result.

To show the claim, define
\[
\alpha_{\bf j} := (zA)^{j_1} B (zA)^{j_2} \cdots B (zA)^{j_{l+1}}, \qquad
{\bf j} := (j_1, \dots, j_{l+1}) \in \mathbb{Z}_{\geq 0}^{l+1}.
\]
Then $\| \alpha_{\bf j} \| \leq |z|^{|{\bf j}|} \| A \|^{|{\bf j}|}
\|B\|^l$, where $|{\bf j}| := j_1 + \cdots + j_{l+1}$. Now we upper bound
$\| E_{\bf n}(z) \|$:
\begin{equation}\label{Epreunion}
\| E_{\bf n}(z) \| = |z|^l \left\| \sum_{{\bf 0} \leq {\bf j}}
\alpha_{\bf j} - \sum_{{\bf 0} \leq {\bf j} \leq {\bf n}} \alpha_{\bf j}
\right\| \leq |z|^l \sum_{{\bf 0} \leq {\bf j}} \| \alpha_{\bf j} \| -
|z|^l \sum_{{\bf 0} \leq {\bf j} \leq {\bf n}} \| \alpha_{\bf j} \|,
\end{equation}
where the first (and bigger) term is bounded above by $|z|^l \|B\|^l
(1-\|zA\|)^{-(l+1)} < \infty$. Now apply a union bound: this difference
of sums (times $|z|^l$) is bounded above by
\[
|z|^l \sum_{j_1 > n_1} \| \alpha_{\bf j} \| + \cdots +
|z|^l \sum_{j_{l+1} > n_{l+1}} \| \alpha_{\bf j} \|.
\]

Combining this with the above upper-bound for each $\| \alpha_{\bf j}
\|$, we get:
\[
\| E_{\bf n}(z) \| \leq |z|^l \sum_{k=1}^{l+1} \sum_{{\bf j} \geq {\bf
0}, j_k > n_k} |z|^{|{\bf j}|} \| A \|^{|{\bf j}|} \|B\|^l = |z|^l
\|B\|^l \sum_{k=1}^{l+1} \frac{\|zA\|^{n_k+1}}{(1-\|zA\|)^{l+1}},
\]
This upper-bound goes to zero as ${\bf n} \to (\infty, \dots, \infty)$.

Similarly, we now show $F_l$ exists for the same $z$:
\begin{align*}
\left\| F_l(z) - \sum_{m=l}^N z^m S_{m,l}(A,B) \right\| \leq &\
\sum_{m=N+1}^\infty |z|^m \binom{m}{l} \|A\|^{m-l} \|B\|^l\\
= &\ \| zB \|^l \sum_{n := m-l = N+1-l}^\infty \binom{n+l}{n} \|zA\|^n,
\end{align*}
and the sum is the tail of the (convergent) negative binomial series
$(1-\|zA\|)^{-(l+1)}$. So this tail goes to zero as desired.

Finally, we claim that the above finite sums $P_{\bf n}(z)$ also converge
to $F_l(z)$ as ${\bf n} \to (\infty, \dots, \infty)$. Indeed, we compute
the error:
\begin{align*}
\| F_l(z) - P_{\bf n}(z) \| = &\ |z|^l \left\| \sum_{m \geq l}
\sum_{|{\bf j}|=m-l} \alpha_{\bf j} - \sum_{{\bf 0} \leq {\bf j} \leq
{\bf n}} \alpha_{\bf j} \right\|,
\end{align*}
and as above, the first sum ``includes'' the second, so this quantity is
at most
\[
\leq |z|^l \sum_{m \geq l} \sum_{|{\bf j}| = m-l} \| \alpha_{\bf j} \| - 
|z|^l \sum_{{\bf 0} \leq {\bf j} \leq {\bf n}} \| \alpha_{\bf j} \|.
\]
But this is exactly the upper bound in~\eqref{Epreunion}, so it also goes
to zero as ${\bf n} \to (\infty, \dots, \infty)$.
\end{proof}

Finally, we show:

\begin{proof}[Proof of Lemma~\ref{Ldecay}]
Let $f(x) = \sum_{k \geq 0} c_k x^k$ on $\mathbb{R}$, with all $c_k \geq
0$. We fix a degree $N>0$, and upper-bound the tail and the remaining
part of $f$. First, the tail: for $k > N$, we have $(c/a)^k < (c/a)^N$,
so
\[
\sum_{k=N+1}^\infty c_k (cx)^k < \sum_{k=N+1}^\infty c_k (ax)^k (c/a)^N
\leq (c/a)^N f(ax).
\]

For the remaining part: since $f$ is not a polynomial, there is some
$M>N$ such that $c_M>0$. Now $f(ax) > c_M (ax)^M$, so
\[
\frac{1}{f(ax)} \sum_{k \leq N} c_k (cx)^k \leq \frac{1}{c_M (ax)^M}
\sum_{k \leq N} c_k (cx)^k,
\]
and this approaches $0$ as $x \to \infty$. Combining these upper bounds,
\[
0 \leq \frac{f(cx)}{f(ax)} < \frac{c^N}{a^N} + \frac{1}{c_M (ax)^M}
\sum_{k \leq N} c_k (cx)^k.
\]
Taking $x \to \infty$ yields $\lim \sup_{x \to \infty} f(cx)/f(ax) \in
[0, (c/a)^N]$. As $N$ is arbitrary, letting $N \to \infty$ proves the
result.
\end{proof}



\end{document}